\documentclass{amsart}
\usepackage{amssymb,amscd,verbatim}
\usepackage[all]{xy}

\newcommand{\ie}{i\@.e\@.\ }
\newcommand{\eg}{e\@.g\@.\ }

\newcommand\Cinf{C^\infty}
\newcommand\dCinf{\dot C^\infty}
\newcommand\dCI{\dot{\mathcal{C}}^\infty}
\newcommand\CI{\mathcal{C}^\infty}
\newcommand\CmI{\mathcal{C}^{-\infty}}
\newcommand\Lap{\Delta}
\newcommand\Ham{\Delta+V}

\newcommand\pa{\partial}
\newcommand\Mand{\text{ and }}

\newcommand\Mif{\text{ if }}
\newcommand\Mforsome{\text{ for some }}
\newcommand\Mfor{\text{ for }}
\newcommand\Mwith{\text{ with }}
\newcommand\Msuchthat{\text{ such that }}

\newcommand\Mwhere{\text{ where }}

\newcommand\ep{\epsilon}
\newcommand\ev{\lambda}

\renewcommand\Sp{\operatorname{Sp}}
\newcommand\cS{\mathcal{S}}
\newcommand\cSp{\mathcal{S}'}
\newcommand\cM{\mathcal{M}}
\newcommand\tcM{\widetilde{\mathcal{M}}}

\renewcommand\inf{\operatorname{inf}}
\renewcommand\sup{\operatorname{sup}}
\newcommand\supp{\operatorname{supp}}
\newcommand\pp{\operatorname{pp}}
\newcommand\ess{\operatorname{ess}}
\newcommand\mic{\operatorname{mic}}
\newcommand\ff{\operatorname{ff}}
\newcommand\new{\operatorname{ff}}
\newcommand\old{\operatorname{o}}
\newcommand\Hess{\operatorname{Hess}}
\newcommand\WF{\operatorname{WF}}
\newcommand\WFsc{\WF_{\text{sc}}}
\newcommand\WFscp{\WF'_{\text{sc}}}
\newcommand\Hsc{H_{\text{sc}}}
\newcommand\Nat{\mathbb{N}}
\newcommand\Real{\mathbb{R}}
\newcommand\Cx{\mathbb{C}}
\newcommand\im{\operatorname{Im}}
\newcommand\re{\operatorname{Re}}
\newcommand\bbR{\mathbb{R}}
\newcommand\bbN{\mathbb{N}}

\newcommand\bbC{\mathbb{C}}
\newcommand\RR{\operatorname{RR}}
\newcommand\Rp{q}
\newcommand\RP{\operatorname{RP}}
\newcommand\RPinc{\operatorname{RP}_-}
\newcommand\RPout{\operatorname{RP}_+}

\newcommand\Cp{z}

\newcommand\Eve{e}
\newcommand\bbS{\mathbb{S}}

\newcommand\nin{\not\in}

\newcommand\Crit{\operatorname{Crit}}
\newcommand\Cv{\operatorname{Cv}}
\newcommand\Max{\operatorname{Max}}

\newcommand\Min{\operatorname{Min}}

\newcommand\Diffsc{\operatorname{Diff}_{\text{sc}}}
\newcommand\Diffb{\operatorname{Diff}_{\text{b}}}
\newcommand\Vb{{\mathcal V}_{\bl}}

\newcommand\bl{{\operatorname{b}}}
\newcommand\scT{{}^{\scl} T}
\newcommand\scbT{{}^{\scl}\overline{T}}

\newcommand\scl{{\operatorname{sc}}}
\newcommand\sccl{{\operatorname{scc}}}

\newcommand\Psisc{\Psi_{\scl}}
\newcommand\Psiscc{\Psi_{\sccl}}
\newcommand\interior{\operatorname{int}}
\newcommand\Dist{C^{-\infty}}
\newcommand\dist{C^{-\infty}}
\newcommand\CIdot{\dot C^{\infty}}

\newcommand\scH{\kern3pt{}^{\scl} H}
\newcommand\Isc{I_{\scl}}

\renewcommand\Re{\operatorname{Re}}
\renewcommand\Im{\operatorname{Im}}

\def\dbyd#1#2{\frac{ \partial #1}{\partial #2}}
\def\Id{\operatorname{Id}}

\newtheorem{lemma}{Lemma}[section]
\newtheorem{prop}[lemma]{Proposition}
\newtheorem{thm}[lemma]{Theorem}
\newtheorem{cor}[lemma]{Corollary}

\newtheorem*{thm*}{Main Results}
\newtheorem*{theorem*}{Theorem}
\newtheorem*{prop*}{Proposition}
\numberwithin{equation}{section}
\theoremstyle{remark}
\newtheorem{rem}[lemma]{Remark}
\theoremstyle{definition}

\theoremstyle{definition}
\newtheorem{Def}[lemma]{Definition}

\newcommand\ci{${\mathcal{C}}^\infty$}
\newcommand\cT{\mathcal{T}}

\newcommand\SEEF{E_{\ess}^{\infty}}
\newcommand\EEF{E_{\ess}^{-\infty}}
\newcommand\FEEF[1]{E_{\ess}^{#1}}
\newcommand\Moe{M}

\newcommand\psit{\tilde\psi}

\newcommand\Bt{\tilde B}
\newcommand\Et{\tilde E}
\newcommand\Pt{\tilde P}

\newcommand\Vt{\tilde V}

\newcommand\inc{-}
\newcommand\out{+}
\newcommand\Vy{{V_0}}
\newcommand\mub{\mu_{\bl}}
\newcommand\qminp{\Rp_{\min}^+}

\newcommand\qmaxp{\Rp_{\max}^+}

\newcommand\intro{%
\renewcommand{\theequation}{I.\arabic{equation}}%
\renewcommand{\thelemma}{I.\arabic{lemma}}}
\newcommand\paperbody{%
\renewcommand{\theequation}{\thesection.\arabic{equation}}%
\renewcommand{\thelemma}{\thesection.\arabic{lemma}}}

\begin{document}
\author[Andrew Hassell]{Andrew Hassell$^*$}
\address{Centre for Mathematics and its Applications, Australian National
   University, Canberra ACT 0200 Australia}
\email{hassell@maths.anu.edu.au}
\author[Richard Melrose]{Richard Melrose$^\dag$}
\address{Department of Mathematics, Massachusetts Institute of Technology,
Cambridge MA 02139}
\email{rbm@math.mit.edu}
\author[Andr\'as Vasy]{Andr\'as Vasy$^{\ddag}$}
\address{Department of Mathematics, Massachusetts Institute of Technology,
Cambridge MA 02139}
\email{andras@math.mit.edu}
\title[Symbolic potentials of order zero]
{Spectral and scattering theory for symbolic potentials of order zero}
\subjclass{35P25, 81Uxx}
\keywords{scattering metrics, degree zero potentials, asymptotics of generalized eigenfunctions, Microlocal Morse decomposition, asymptotic completeness}

\begin{abstract} The spectral and scattering theory is investigated for a
generalization, to scattering metrics on two-dimensional compact manifolds
with boundary, of the class of smooth potentials on $\bbR^2$ which are
homogeneous of degree zero near infinity. The most complete results
require the additional assumption that the restriction of the potential to the
circle(s) at infinity be Morse. Generalized eigenfunctions associated to
the essential spectrum at non-critical energies are shown to originate both
at minima and maxima, although the latter are not germane to the $L^2$
spectral theory. Asymptotic completeness is shown, both in the traditional
$L^2$ sense and in the sense of tempered distributions. This leads to a
definition of the scattering matrix, the structure of which will be
described in a future publication.
\end{abstract}

\thanks{$^*$ Supported in part by the Australian Research Council,
$^\dag$ Supported in part by the National Science Foundation under grant
\#DMS-9622870,
$^\ddag$ Supported in part by the National Science Foundation under grant
\#DMS-99-70607}
\maketitle

\intro
\section*{Introduction}

In \cite{Herbst1}, Herbst initiated the spectral analysis and scattering
theory for the Schr\"o\-dinger operator $P = \Lap+V$ on $\bbR^n$ where $V$
is a smooth, real-valued potential which is homogeneous of degree zero near
infinity. Further results for such potentials were obtained by Agmon, Cruz
and Herbst \cite{Agmon-Cruz-Herbst1} and Herbst and Skibsted
\cite{Herbst-Skibsted1}. Such a potential is so large that it deforms the
geometry near spatial infinity and consequently the scattering theory of
$P$ is quite different from that of a Schr\"odinger operator with a short-, or
even a long-, range potential. Herbst \cite{Herbst1} and Herbst-Skibsted
\cite{Herbst-Skibsted1} showed that solutions of the time-dependent
Schr\"odinger equation, with $L^2$ initial data, concentrate for large times
near critical directions of the potential restricted to the sphere at
infinity.

In this paper we study in some detail the tempered eigenspaces of $P,$
consisting for each fixed $\ev$ of those functions $u,$ of polynomial
growth, satisfying $(P - \ev)u = 0.$ We restrict attention to the two-
dimensional case but generalize from Euclidean space, $\bbR^2,$ to
scattering metrics (defined below) on arbitrary compact 2-manifolds $X$
with boundary. Such metrics give the interior of $X$ the structure of a
complete Riemannian manifold with curvature vanishing at infinity. In this
wider context the natural condition replacing, and generalizing,
homogeneity of the potential near infinity is its smoothness up to the
boundary of $\pa X.$ For the Euclidean case, now appearing through the
radial compactification to a ball, this allows potential which are
classical symbols of order $0.$ Generally we assume that $V_0,$ the
restriction of the potential to the boundary, is Morse, so $V_0$ has only
nondegenerate critical points.

The scattering wavefront set from \cite{Melrose43} allows the
eigenfunctions of $P$ to be analyzed microlocally and hence to be related
to the classical dynamics of $P.$ It is the asymptotic behaviour which is
relevant here and this is described through the contact geometry of the
cotangent bundle over the boundary; the classical dynamics of $P$ are then
given by a Legendre vector field over the boundary. The critical points of
this vector field correspond to the critical points of $V_0.$ We give
microlocal versions of the familiar notions of incoming and outgoing
eigenfunctions, depending on which integral curves of the Legendre vector
fields are permitted in the scattering wavefront set. By direct
construction we show that there is a non-trivial space of outgoing
microlocal eigenfunctions associated to each critical point $q$ of $V_0$,
at energies $\ev$ above $V_0(q).$ For a minimum this space is isomorphic to
the space of Schwartz functions on a line, for a maximum it is isomorphic
to the space of formal power series (\ie Taylor series at a point) in one
variable. These isomorphisms are realized directly in terms of asymptotic
expansions near the corresponding critical point of $V_0;$ the form of the
expansion reflects the behaviour of the classical trajectories near these
points.

For non-critical eigenvalues there is a well-defined space of `smooth'
eigenfunctions, which may be characterized either as the range of the
spectral projection on Schwartz functions or else as the space of tempered
eigenfunctions with no purely tangential oscillation at infinity (this
latter condition is conveniently expressible in terms of the scattering
wavefront set). We show this space to have a `microlocal Morse
decomposition' in terms of the spaces of microlocally outgoing
eigenfuncitons associated to the critical points. The space of all tempered
eigenfunctions is identified by a boundary (or Green) pairing with the dual
of the space of smooth eigenfunctions, so it has a similar decomposition in
terms of the duals of the microlocal spaces. The latter are isomorphic
either to Schwartz distributions on a line (for minima) or polynomials in
one variable (for maxima). The map from generalized eigenfunctions to the
duals of the microlocally outgoing eigenspaces thus correspond to
`generalized boundary data'; the collective boundary data fixes the
eigenfunction via the inverse map, which is an analogue of the Poisson
operator. There are two such families of boundary data (corresponding to an
underlying choice of orientation of the Legendre vector field) and the
isomorphism between them is the scattering matrix.

These results are new even in the Euclidean setting. Although for the most
part the results here are limited to the two-dimensional case there are
obvious generalizations to higher dimensions. These extensions will be
discussed in future publications as will the related microlocal description
of the Poisson operators and scattering matrix.

Superficially, the existence of smooth eigenfunctions associated to a
maximum might appear to be in conflict with the results of Herbst and
Skibsted on decay, near maxima, of solutions to the time-dependent
Schr\"odinger equation \cite{Herbst-Skibsted1}. This is not the case, since
these eigenfunctions are, near the maximum point, in a weighted $L^2$ space
which implies that they are too small to violate this decay estimate. In
fact an eigenfunction which is microlocally outgoing near a maximum, and is
non-trivial there, is necessarily non-trivial (but not outgoing) near the
minimum. Thus there are no eigenfunctions \emph{purely} associated to a
maximum, rather these outgoing eigenfunctions can be thought of as
corresponding to the completion of the Schwartz functions on a line to
smooth functions on a circle by the addition of the Taylor series at the
point at infinity. When completed in a norm equivalent to the $L^2$ norm on
initial data for the Schr\"odinger operator, these extra terms are already
in the closure of outgoing eigenfunctions at the minimum, thus they do not
`appear' in the $L^2$ theory.

\

Next we give a more precise description of our results. Thus, let $X$ be
a compact manifold with boundary where, for the moment, we do not restrict
the dimension. The boundary $Y=\pa X$ consists of a finite union of compact
manifolds without boundary, $Y=Y_1\cup\dots\cup Y_N.$ It is always possible
to find a boundary defining function on $X,$ $x\in\CI(X)$ such that
$x\ge0,$ $Y=\{x=0\}$ and $dx\not=0$ on $Y.$ A Riemannian metric on the
interior of $X$ is a scattering metric if, for some choice of defining
function, it takes the form
\begin{equation}
g=\frac{dx^2}{x^4}+\frac h{x^2},\ h\in\CI(X;S^2X),\ h_0=h\big|_Y
\text{ positive definite} .
\label{HMV11}\end{equation}
That is, $h=x^2(g-x^{-4}dx^2)$ is a smooth 2-cotensor on $X$ which restricts to a
metric on $Y.$ In this setting we consider a real potential $V\in\CI(X)$ and
examine the spectral and scattering theory of $\Lap+V$ where $\Lap$ is the
Laplace operator of a scattering metric. The Euclidean case is included
since the Euclidean metric is a scattering metric on the radial
compactification of $\bbR^n$ to a ball (or half-sphere) and a smooth
potential which is homogeneous of degree zero near infinity is a smooth
function up to the boundary of the radial compactification. The lower order
parts of the Taylor series of $V$ allow additional terms decaying at
Coulomb rate or faster.

Under these assumptions, for $\dim X\ge2,$ the Schr\"odinger operator
$P=\Lap+V$
is self-adjoint and has continuous spectrum of infinite multiplicity
occupying the interval $[\kappa ,\infty)$ where 
\begin{equation}
\kappa =\inf_YV.
\label{HMV12}\end{equation}
In addition there may be point spectrum in the interval $[m,K],$  
\begin{equation}
m=\inf_XV,\ K=\sup_YV
\label{HMV13}\end{equation}
which is discrete in the open set $(-\infty,K)\setminus \Cv(V),$  
\begin{equation}
\Crit(V)=\{p\in Y;d_YV(p)=0\},\ \Cv(V)=V(\Crit(V)).
\label{HMV14}\end{equation}
The main interest lies in the continuous spectrum which we analyse in
detail here under the assumption that $\dim X=2$ and $V_0 = V\big|_Y$ is
Morse. For simplicity in this Introduction we outline the results under the
additional assumption that $Y$ has only one component circle and $V_0$
is perfect Morse, so only has a global maximum, $\Cp_{\max},$ and minimum
$\Cp_{\min},$ forming $\Crit(V).$ These restrictions are removed in the
body of the paper.

Our central result is the parameterization of all tempered distributions
associated to the continuous spectrum; this constitutes a distributional
form of `asymptotic completeness'. Thus we examine 
\begin{equation}
E(\ev )=\left\{u\in\CmI(X);(\Ham-\ev )u=0\right\}.
\label{HMV15}\end{equation}
The space of extendible distributions, $\CmI(X),$ reduces precisely
to the space of tempered distributions in the sense of Schwartz in case $X$
is the radial compactification of $\bbR^n.$ For eigenfunctions this is
equivalent to a polynomial bound. For $\ev <\kappa,$ $E(\ev)$ is
finite-dimensional and consists of square-integrable eigenfunctions. In
fact, for any $\ev \notin \Cv(V)$, we show in
Proposition~\ref{prop:eigenvalues} that  
$E_{\pp}(\ev )=E(\ev )\cap L^2(X)$ is contained in $\dCI(X)$ and
finite dimensional 
(where $L^2(X)$ is computed with respect to the Riemannian volume form)
and is trivial for $\ev >K$.
Since $\CmI(X)$ is the dual of $\dCI(X)$ this allows us to consider 
\begin{equation}
\begin{gathered}
\EEF(\ev )=\left\{u\in\CmI(X);(\Ham-\ev )u=0,\ \langle u,v\rangle
=0\ \forall\ v\in E_{\pp}(\ev )\right\},\\
E(\ev )=\EEF(\ev )\oplus E_{\pp}(\ev ).
\end{gathered}
\label{HMV17}\end{equation}

The structure of the space $\EEF(\ev)$ depends on $\ev;$ there
are three distinct cases, corresponding to the values of $\kappa,$ $K$ and
the additional transition point
\begin{equation}
\ev _{\Hess}=\kappa +2V''(\Cp_{\min}),\ \kappa =V(\Cp_{\min})
\label{HMV18}\end{equation}
where the derivatives are with respect to boundary arclength. Clearly
$\ev_{\Hess}>\kappa$ but the three possibilities  $\ev_{\Hess}<K$,
$\ev_{\Hess}=K,$ $\ev_{\Hess}>K$
may all occur for different $V$ (or for the same $V$ but different
metrics). For any particular problem only two or three of the
following four intervals can occur
\begin{equation}
\begin{cases}
\kappa <\ev <\min(\ev _{\Hess},K)&\text{``Near minimum''}\\
\ev _{\Hess}<\ev <K&\text{``Hessian range''}\\
K<\ev <\ev _{\Hess}&\text{``Mixed range''}\\
\max(\ev_{\Hess},K)<\ev &\text{``Above thresholds."}
\end{cases}
\label{HMV.9}\end{equation}
The critical values $\kappa$ and $K$ of $V_0$ are thresholds corresponding to
changes in the geometry of the fixed energy (\ie characteristic)
surface of the classical problem. Above these thresholds both the maximum
and the minimum correspond to zeroes, which we also call radial
points, 
of the classical flow (at infinity); the Hessian transition corresponds to
an energy at which there is a change in the local geometry at the radial points
corresponding to $\Cp_{\min}$. 

To analyse $\EEF(\ev)$ we use the limit of the resolvent on the
spectrum. This exists as an operator for
$\ev\in(\kappa,+\infty)\setminus(\Cv(V)\cup\sigma_{\pp}(\Ham)),$ as a
consequence of an appropriate version of the Mourre estimate, or of
microlocal estimates closely related to it, see \cite{Herbst1} for the
proof in the Euclidean setting.

\begin{thm}[Herbst \cite{Herbst1}, see also Theorem~\ref{thm:lim-abs}]
The resolvent 
\begin{equation}
R(\ev\pm it)=(\Ham-(\ev\pm it))^{-1},\ t>0,\ \ev\nin\Cv(V)
\end{equation}
extends continuously to the real axis, \ie $R(\ev\pm i0)$ exist, as
bounded operators 
\begin{equation}
x^{1/2+\delta}L^2(X)\ominus E_{\pp}(\ev)\longrightarrow
x^{-1/2-\delta}L^2(X)\ominus E_{\pp}(\ev),\ \forall\ \delta>0.
\end{equation}
\end{thm}
\noindent Note that
\begin{equation}
f\in x^{1/2+\delta}L^2(X)\ominus E_{\pp}(\ev) \Longrightarrow
[R(\ev+i0)-R(\ev-i0)]f\in \EEF(\ev)
\end{equation}
and conversely (essentially by Stone's theorem) the range is dense.

For $\ev\notin\Cv(V)$ we can thus define spaces of
`smooth' eigenfunctions $\SEEF(\ev)$ by 
\begin{equation}
\SEEF(\ev)=
[R(\ev+i0)-R(\ev-i0)](\dCI(X)\ominus E_{\pp}(\ev))
\label{HMV.10}\end{equation}
and $\SEEF(\ev)$ inherits a Fr\'echet topology from $\dCI(X).$ We give a
microlocal characterization of these spaces below. In all of the
non-transition regions
\begin{equation}
\SEEF (\ev)\subset \EEF(\ev)\text{ is dense}
\label{HMV.11}\end{equation}
in the topology of $\CmI(X).$
In fact, approximating sequences
can be constructed rather explicitly by extending
$R(\ev\pm i0)$ to distributions satisfying a scattering
wave front set condition; see Proposition~\ref{prop:density}.

One of the aims of scattering theory is to give explicit, and geometric,
parameterizations of the continuous spectrum. The reality of the Laplacian
implies that there are two equivalent such representations, which are
interchanged by complex conjugation, and the scattering matrix gives the
relation 
between them. The main task here is to describe the structure of the smooth
eigenfunctions in the four non-transition regions in \eqref{HMV.9}; our
explicit parameterizations arise as the `leading terms' in the complete
asymptotic expansions that these eigenfunctions possess. This also
identifies $\SEEF(\ev)$ in terms of familiar Fr\'echet spaces. The
expansions are intimately connected to the corresponding classical problem,
which is described in detail in \S\ref{S.Classical}. For $\ev\notin\Cv(V)$
the classical system consists of a Legendre vector field $W$ (so defined up
to a conformal factor) on a compact hypersurface, the characteristic
variety $\Sigma(\ev),$ in a contact manifold. The dual variable 
to the variable $1/x$ (which is the radial variable $r =
|z|$ in the Euclidean case) is a function, $\nu$ on $\Sigma
(\ev)$ with $W\nu \geq 0$. 

The microlocal structure of this problem, near each radial point, when
transformed from the scattering to the traditional context (essentially by
Fourier transformation) is the problem considered by Guillemin and
Schaeffer in \cite{MR55:3504}. Although we do not use their work explicitly
here, several of our results could be proved by their techniques, provided
one made the additional assumption of non-resonance of the linearization of
$\Ham-\ev$ at the radial points.

The projections of the radial points (that is, the zeroes of $W$) are the
critical points of $V$ on the boundary. Each such critical point
$\Cp\in\Crit(V)$ 
corresponds to two radial points $q_{\pm}\in\Sigma (\ev)$ with
$\pm\nu(q_{\pm})>0$ if $\ev>V(\Cp);$ it
corresponds to a singular point on $\Sigma (\ev)$ if $\ev=V(\Cp).$ The two
points correspond to opposite, non-zero, values of $\nu.$ We shall denote
by $\RPout(\ev)\subset\Sigma (\ev)$, resp. $\RPinc(\ev)\subset\Sigma
(\ev)$, the set of radial points at which $\nu>0$, resp. $\nu < 0$. 
This can also be identified with the
subset $\left\{q\in\Crit(V);V(q)<\ev\right\}$ and we let
\begin{equation}
\RPout(\ev)=\Max_+(\ev)\cup\Min_+(\ev)
\label{HMV.183}\end{equation}
be the decomposition into points associated to maxima and to minima of
$V_0$. We shall describe points in $\RPout(\ev)$, resp. $\RPinc(\ev)$,
as outgoing, resp. incoming, radial points.

Our examination of the structure of the eigenfunctions is based on the
description of the microlocally outgoing eigenfunctions associated to each
of the 
outgoing radial points. The notion of microlocality here is with respect to
the scattering wavefront set, which is the notion of wavefront set
associated to the scattering calculus in \cite{Melrose43}. By microlocal
elliptic regularity, for any $u\in\CmI(X),$ and any open set $O$
\begin{equation*}
\WFsc\left((\Ham-\ev)u\right)\cap O=\emptyset
\Longrightarrow \WFsc(u)\cap O\subset\Sigma(\ev).
\end{equation*}
If $\Rp\in\RP_\pm(\ev)$ and $O$ is a sufficiently small open neighbourhood of
$\Rp,$ which meets $\Sigma(\ev)$ in a $W$-convex neighbourhood with each
$W$ curve meeting $\nu=\nu(\Rp)$ we set
\begin{multline}
\tilde E_{\mic,\out}(O,\ev)=
\big\{u\in\CmI(X);O\cap\WFsc\big((P-\ev)u\big)=\emptyset\\
\Mand \WFsc(u)\cap O\subset\{\nu \geq \nu(\Rp)\}\big\},
\label{HMV.42}\end{multline}
with $\tilde E_{\mic,\inc}(O,\ev)$ defined by reversing the inequality. 
We may consider this as a space of microfunctions, $E _{\mic,\out}(\Rp,\ev)$,
by identifying elements when they differ by functions with wavefront set
not meeting $O$. The result is then independent of the choice of $O.$
That is, such microlocal solutions
are determined by their behaviour in an arbitrarily small neighbourhood of
$\Rp.$ Our convention is to call elements of either $\tilde
E_{\mic,\out}(\Rp,\ev)$ 
or $E_{\mic,\out}(\Rp,\ev)$ `microlocally outgoing eigenfuncions at
$\Rp\in\RPout(\ev).$' 
The structure of the microlocal eigenfunctions at a radial point is
essentially determined by the linear part of $W$, the normalized
eigenvalues of which are given by \eqref{HMV.37} below. The character of
the radial point is determined by the quantity
\begin{equation}
r(\ev)= r(\Rp,\ev)=\frac{1}{2}-\sqrt{\frac{1}{4}-\frac{V''_0}{2\tilde\nu^2}},\ 
\tilde\nu=\nu(\Rp),
\label{HMV.186}\end{equation}
where $V''_0$ is the Hessian with respect to boundary arclength at the
critical point.

The case $\Rp\in\Max_+(\ev)$ is the most straightforward, with $W$ having a
saddle point at $\Rp$ at which $r(\Rp,\ev)<0.$ Elements $u \in
E_{\mic,+}(\Rp, \ev)$ are then (\ie have representatives which are) 
functions of the form
\begin{equation}
u=e^{i\phi_{\ev}/x}x^{\frac12+\epsilon}v,\ v\in L^\infty(X)\cap\Cinf(X^\circ),\
\epsilon=-\frac12 r(\Rp,\ev)>0,
\label{HMV.41}\end{equation}
\ie Legendre distributions (see \cite{MZ}) 
where $\phi_{\ev}\in\CI(X)$ is a real function parameterizing the unstable
manifold for $W$ at $\Rp$ and $v$ has a complete
asymptotic expansion with $\Cinf$ coefficients as $x \to 0$
in powers of $x$ starting with a term in which the power is imaginary.
Appropriately evaluating the terms in the expansion of $v$ at
$\Cp=\pi(\Rp)$ gives a parameterizing map,
\begin{equation} 
M_+(\Rp,\ev):E_{\mic,\out}(\Rp,\ev)\ni u\longrightarrow \bbC[[x]],
\label{HMV.44}\end{equation}
which is an isomorphism to the space of formal power series in one
variable, or equivalently the space of arbitrary sequences
$s:\bbN\longrightarrow\bbC.$ This result is proved in
Theorem~\ref{thm:saddle-smooth}. Notice that the factor
$x^\epsilon$ in \eqref{HMV.41} means that these microlocal eigenfunctions
are `small' near the maximum, as they lie in the space $x^{-\frac12}L^2(X)$
locally near $z$. Here, the exponent $-1/2$ is critical for
eigenfunctions; any eigenfunction which is {\it globally} in
$x^{-\frac12}L^2(X)$ is actually rapidly decreasing.

When $\Rp\in\Min_+(\ev)$ the microlocal space $E_{\mic,\out}(\Rp,\ev)$ can be
identified with a global space of eigenfunctions
\begin{equation} 
E_{\mic,\out}(\Rp,\ev) \simeq
\left\{u\in\EEF(\ev);\WFsc(u)\cap\{\nu >0\}\subset\{\Rp\}\right\}
\subset\SEEF(\ev)
\label{HMV.184}\end{equation}
(see Corollary~\ref{cor:min-ident-8}). 
In this case, the structure of
$E_{\mic,\out}(\Rp,\ev)$ undergoes a distinct change as $\ev$ crosses
the Hessian threhold. First consider the `near minimal' range $\kappa
<\ev<\ev_{\Hess}.$ In this case $W$ has a center at $\Rp$ with
$\re r(\Rp,\ev)=\frac12,$ and the square root in \eqref{HMV.186} is
imaginary.
It is convenient to introduce the space $[X;\Cp]_{\frac12}$
obtained by parabolic blow up, with respect to the boundary, of $X$ at $\Cp
= \pi(\Rp)$
(see section~\ref{sec:blow-up} for a discussion of blowups). 
Then the space $\CI_{\ff}([X;\Cp]_{\frac12})$ consists of the smooth
functions on this blown-up space which vanish rapidly at all boundary faces
other than the `front face' created by the blowup. Each microlocal 
eigenfunction has a representation in the form of an infinite sum
\begin{equation}\begin{gathered}
u=e^{i\tilde\nu/x}\sum\limits_{j\ge0}x^{\frac14+i((2j+1)\alpha + V_1(\Cp))
}a_j v_j,\ a_j\in
\cS(\bbN), \ v_j \in \CI_{\ff}([X;\Cp]_{\frac12}), \\
\tilde\nu=\sqrt{\ev - V_0(\Cp)},\ V_1(z) = (\pa_x V)(z),
\end{gathered}\label{HMV.47}\end{equation}
where the $a_j,$ as
indicated, form a Schwartz series in $j$, and the $v_j$ restricted to
$\ff([X;\Cp]_{\frac12})$ are $L^2$-normalized eigenfunctions for a 
self-adjoint globally elliptic operator on
$\ff([X;\Cp]_{\frac12})$, namely a harmonic oscillator with
eigenvalues $(j+1/2)\alpha$. The sequence of coefficients, $a_j,$ is an
arbitrary Schwartz sequence and the map 
obtained by normalization from \eqref{HMV.47} gives an isomormphism
\begin{equation}
\Moe_+(\Rp,\ev):E_{\mic,\out}(\Rp,\ev)\longrightarrow
\dCI(\ff[X;\Cp]_{\frac12}),\  
u\longmapsto \sum\limits_{j\ge0}(a_j v_j\big|_{\ff}), \kappa
<\ev<\ev_{\Hess},
\label{HMV.48}\end{equation}
onto the Schwartz functions on a line. We prove this statement
in Theorem~\ref{thm:center-smooth}. Note that the normalization of $v_j$ is
only well-defined up to a factor
$e^{is(j)},$ where $s(j)$ is a real sequence. When it is analysed
in a subsequent paper, it will convenient to introduce such terms so that
$\Moe_+(\Rp,\ev)$ (or more precisely its inverse) will be a Fourier
integral operator of an appropriate type. 

The case $\ev_{\Hess}\le\ev,$ at a minimum, is similar, except that the
discrete spectrum of the model problem has collapsed, and $\Rp$ has become
a sink for $W.$ At $\ev=\ev_{\Hess}$
the `local homogeneity' is still parabolic as it is for $\ev<\ev_{\Hess}.$
As $\ev$ increases the local homogeneity is determined by $r=r(\Rp,\ev)$ 
given by \eqref{HMV.186}. Corresponding to this
homogeneity we introduce the space $[X;\Cp]_r$ obtained by blow up with the
correct scaling of the tangential variable relative to the normal
variable. Now $r$ is not necessarily rational, hence we need to replace
$\CI_{\ff}([X;\Cp]_r)$ by a space of polyhomogeneous functions allowing
both homogeneities $x^n$ and $x^{nr}$ where $n$ is a non-negative integer.
This space, $\mathcal{A}^{I}_{\text{phg},\ff}([X;\Cp]_r)$, is described in more
detail in Section~\ref{S.Out.sink}, but here we mention that its
elements are in particular continuous on $[X;\Cp]_r,$ and vanish at the
boundary away from $\ff.$ If $1/r(\ev)$ is not an integer, then the
microlocal eigenfunctions $u \in E_{\mic,+}(\Rp,\ev)$ take the form
\begin{equation}
u=e^{i\phi_{\ev}/x}x^{\frac{1-r}2+ic}a,\
a\in\mathcal{A}^{I}_{\text{phg},\ff}([X;\Cp]_r),\
r=r(\ev),
\label{HMV.49}\end{equation}
where $\phi_{\ev}\in\CI(X)$ is now a real function parameterizing a smooth
(minimal) Legendre submanifold through $\Rp$ and $c=c(\ev)$ is a
real constant defined from the subprincipal symbol of $P.$ Again,
restriction to the front face of the blown-up space gives an isomorphism 
\begin{equation}
\Moe_+(\Rp,\ev):E_{\mic,\out}(\Rp,\ev)\longrightarrow\dCI(\ff[X;\Cp]_r),\
\Rp\in\Min_+(\ev),
\label{HMV.50}\end{equation}
which, after selection of coordinates, maps isomorphically onto the Schwartz
functions on an associated line. The passage of $\ev$ across the Hessian
threshold results in a change in the underlying parameterizing space, which
from then onwards depends on $\ev.$ This result is proved in
Theorem~\ref{thm:sink-smooth}, together with a statement of the minor
changes needed when $1/r(\ev)$ is an integer.

The space of smooth eigenfunctions can be directly related
to these microlocal spaces and this is especially simple in case the
potential is perfect Morse on the one boundary component. There is then a
restriction map `to the maximum'
\begin{equation}
\SEEF(\ev)\longrightarrow E_{\mic,\out}(\Rp,\ev),\
\Rp\in\Max_+(\ev)\Mfor K<\ev.
\label{HMV.51}
\end{equation}

\begin{thm}\label{HMV.52} [Microlocal Morse decomposition]
If $\dim X=2,$ $\pa X$ is connected, $V_0$ is perfect Morse
and $\ev\nin\Cv(V),$ the
map \eqref{HMV.184} is an isomorphism
\begin{equation}
E_{\mic,\out}(\Rp^+_{\min},\ev)\longrightarrow \SEEF(\ev),\Mfor\kappa <\ev<K,
\label{HMV.185}\end{equation}
whereas, for $\ev>K,$ it combines with the map \eqref{HMV.51} to give a
short exact sequence 
\begin{equation}
0 \longrightarrow E_{\mic,\out}(\Rp^+_{\min},\ev)\longrightarrow
\SEEF(\ev)\longrightarrow 
E_{\mic,\out}(\Rp^+_{\max},\ev) \longrightarrow 0.
\label{HMV.54}\end{equation}
\end{thm}
\noindent Whilst \eqref{HMV.54} does not split, there is a
left inverse on any finite dimensional subspace.  
Combined with the isomorphisms \eqref{HMV.44} and
\eqref{HMV.50} this result gives a quite complete description of the smooth
eigenfunctions, which we could call `smooth asymptotic completeness at
energy $\ev$' in view of \eqref{HMV.10}.

There is a non-degenerate sesquilinear `Green pairing' on $\SEEF(\ev)$
which extends to a continuous bilinear map (see section~\ref{sec:pairing})
\begin{equation}
B:\SEEF(\ev)\times\EEF(\ev)\longrightarrow \bbC
\label{HMV.177}\end{equation}
allowing $\EEF(\ev)$ to be identified with the dual space of
$\SEEF(\ev).$ This pairing also restricts to a pairing between
$E_{\mic,\out}(\Rp,\ev)$ and $E_{\mic,\inc}(\Rp,\ev)$ for each radial point
$\Rp$ and
$\ev>V(\pi(\Rp)).$ The adjoint of the map \eqref{HMV.44} becomes an
isomorphism 
\begin{equation}
\bbC[x]\longleftrightarrow E_{\mic,\inc}(\Rp,\ev),\ \Rp\in\Max_+(\ev),
\label{HMV.181}\end{equation}
where $\bbC[x]$ is the space of finite power series (\ie
polynomials) --- see Theorem~\ref{thm:saddle-big}. 
The range space in \eqref{HMV.181} again consists of functions
which are Legendre distributions near $\Rp$; only
the trivial elements of $E_{\mic,\inc}(\Rp,\ev)$ are in the space
$x^{-\frac12}L^2(X)$ near $\pi(\Rp)$ (this is dual to the property 
\eqref{HMV.41}). 
Similarly the adjoint of \eqref{HMV.50} extends to an isomorphism 
\begin{equation}
\cSp(\bbR_r)\longleftrightarrow E_{\mic,\inc}(\Rp,\ev),\ \Rp\in\Min_+(\ev).
\label{HMV.180}\end{equation}

In the perfect Morse case (with $\pa X$ connected) microlocal
decompositions of the distributional eigenspaces are given by the
identification 
\begin{equation}
E_{\mic,\inc}(\qminp,\ev)\longrightarrow \EEF(\ev),\ \Mfor\kappa<\ev<K
\label{HMV.178}\end{equation}
and for $\ev > K$ the short exact sequence 
\begin{equation}
0 \longrightarrow E_{\mic,\inc}(\qmaxp,\ev)\longrightarrow
\EEF(\ev)\longrightarrow E_{\mic,\inc}(\qminp,\ev) \longrightarrow 0. 
\label{HMV.179}\end{equation}
Combined with the isomorphism \eqref{HMV.181} and \eqref{HMV.180} this gives a
decomposition of the full distributional eigenspaces in terms of 
microlocal eigenspaces. 

One can compare the sequence in \eqref{HMV.54} with the familiar short
exact sequence
\begin{equation} \cS(\bbR)\hookrightarrow
\CI(\bbS)\longrightarrow \bbC[[x]]
\label{HMV.58}\end{equation}
where the second map is the passage to Taylor series at a point and the
first is by identification of the complement of this point with a
line. Then \eqref{HMV.179} corresponds to the dual sequence for
distributions 
\begin{equation}
\bbC[x]\longrightarrow \CmI(\bbS)\longrightarrow \cSp(\bbR)
\label{HMV.182}\end{equation}
where the polynomials correspond to the (Dirac delta) distributions at the
point.

In this way we may think of $\SEEF(\ev)$ as a space of test functions and
$\EEF(\ev)$ as the dual space of distributions. It is then natural to think
of the intermediate `Sobolev spaces' of eigenfunctions. In this case they
may be defined directly:
\begin{equation}
E^s_{\ess}(\ev)
=\big\{u\in\EEF(\ev);\ \WFsc^{0,s-1/2}(u)\cap\{\nu=0\}=
\emptyset\big\},\ s\in\bbR.
\label{HMV.157}\end{equation}
The case $s=0$ is particularly interesting since it is related to the notion
of asymptotic completeness. In Section~\ref{sec:pairing} we prove the
following theorem.

\begin{thm} [Asymptotic completeness at energy $\ev$] If $\ev\nin\Cv(V)$ then
\begin{equation}
\begin{gathered}
E_{\ess,+}(\qminp,\ev)\hookrightarrow E^0_{\ess}(\ev)\text{ is dense and}\\
\Moe_+(\qminp,\ev)\text{ extends to a unitary operator }
E^0_{\ess}(\ev)\longrightarrow L^2(\bbR_r),
\end{gathered}
\label{eq:Moe(ev)-0}\end{equation}
with respect to the metric induced by $B$ from \eqref{HMV.177} on the one
hand, and $L^2(\bbR_r)$ 
equipped with a translation-invariant measure induced by $h_0$ from
\eqref{HMV11} on the other.
\label{thm:AC-ev}\end{thm}

\noindent Integration with respect to $\ev$ gives the usual version of
asymptotic completeness which we state in Theorem~\ref{thm:AC}.

These two theorems, combined with \eqref{HMV.41}, have immediate
implications for the `size' of $u_\pm=R(\ev\pm i0)f,$ $f\in\dCI(X),$ away
from the minima of $V_0.$ Namely, for any pseudodifferential 
operator $Q$ with 
$\qminp\nin\WFscp(Q),$ in particular for multiplication operators
supported away from $\Cp_{\min},$ $Qu_+\in x^sL^2(X)$ for all
$s<(-r(\qminp,\ev)-1)/2.$ Since $r(\qminp,\ev)<0,$ this is an
improvement over the statement that $Qu_+\in x^{-1/2}L^2(X).$ The
statement $Qu_+\in x^{-1/2}L^2(X)$ can be interpreted as the absence of
channels at the maxima 
of $V_0$ and is closely related to the work of Herbst and Skibsted
\cite{Herbst-Skibsted1}; indeed this can be seen from the injectivity of
$\Moe_+(\ev)$ in \eqref{eq:Moe(ev)-0}. Thus, on the one hand, our results
strengthen theirs in this special case (\ie when $\dim X=2$) by showing
that $u_+$ has additional decay away from the minima of $V_0$ and
gives the precise asymptotic form of $u_+.$ On the other hand, it also
shows that the decay is not rapid. One interpretation of this phenomenon,
supported by the short exact sequence \eqref{HMV.54} and isomorphism
\eqref{eq:Moe(ev)-0}, is that the $L^2$-theory does not `see' the maxima of
the potential (because $u_+$ is too small there), while the smooth theory
(working modulo $\dCI(X))$ does.

As well as the identification \eqref{eq:Moe(ev)-0} arising from the
outgoing parameterization there is the corresponding incoming
parameterization, corresponding to $\nu<0$. The
scattering matrix is then the unitary operator given by the composite
\begin{equation}
S(\ev)=\Moe_+\Moe_-^{-1}:L^2(\bbR_r)\to L^2(\bbR_r).
\label{eq:S-def}
\end{equation}
In a future publication it will be shown to be a Fourier integral operator,
of an appropriate type, associated with the limiting scattering relation
arising from the vector field $W.$

The classical dynamical system underlying the eigenvalue problem is
analyzed in \S~\ref{S.Classical}, and \S~\ref{sec:blow-up} contains a brief
description of the blow-up procedure used later to discuss expansions. In
\S~\ref{sec:lim-abs} extension of results of Herbst on the limiting
absorption principle to the present more general setting are discussed;
they are proved by positive commutator methods in \S~\ref{S.Propagation}. The
basic properties of the (incoming and outgoing) microlocal eigenfunctions
associated to critical points of the classical system (\ie radial points)
are presented in \S~\ref{sec:micro} and an outline of the methods used to
analyse them is given in \S~\ref{S.Test}. The next four sections contain the
detailed analysis of incoming microlocal eigenfunctions at radial points
associated to minima and both incoming and outgoing eigenfunctions
associated to maxima. The decomposition of smooth eigenfunctions in terms
of these microlocal components is described in \S~\ref{sec:Morse} and this is
used, together with the natural pairing on eigenfunctions, to describe all
tempered eigenfunction, the scattering matrix and various forms of
asymptotic completeness in the final section.

The authors are delighted to acknowledge conversations with Ira Herbst and
Erik Skibsted, who were working on closely related problems at the time,
and Gunther Uhlmann for useful advice. 
They are also grateful for the hospitality and support of the Erwin
Schr\"odinger Institute in Vienna, Austria, where some of the work on this
paper was carried out.

\paperbody
\section{Classical problem}\label{S.Classical}

The classical system formally associated to the Schr\"odinger operator
$\Lap+V$ on $\bbR^n$ is generated by the Hamiltonian function 
\begin{equation}
p=|\zeta|^2+V(z)\in\CI(\bbR^n_z\times\bbR^n_\zeta).
\label{HMV.12}\end{equation}
In fact, only the `large momentum' (the usual semiclassical limit)
$\zeta\to\infty$ or the `large distance' limit $z\to\infty$ are relevant
to the behaviour of solutions. Since we are interested in eigenfunctions 
\begin{equation}
(\Lap+V-\ev )u=0,
\label{HMV.13}\end{equation}
for finite $\ev,$ only the latter limit is important and hence our
microlocal analysis takes place in the vicinity of the finite (fixed)
energy, or characteristic, surface which
can be written formally
\begin{equation}
\Sigma (\ev )=\left\{|\zeta |^2+V(z)=\ev,\ |z|=\infty \right\}.
\label{HMV.14}\end{equation}
Setting $|z|=\infty$ is restriction to the sphere at infinity
for the radial compactification $\overline{\bbR^n}$
of $\bbR^n$ so $\Sigma (\ev
)\subset\bbS^{n-1}_\omega \times\bbR^n,$ $z=|z|\omega.$ This phase space at
infinity is a contact manifold where the contact form is
\begin{equation}
\alpha =-d\nu +\mu \cdot d\omega ,\ \zeta =\mu \oplus \nu\omega,
\ \mu\in\omega ^\perp.
\label{HMV.15}\end{equation}

In the general case of a scattering metric on a compact manifold with
boundary (see \cite{Melrose43}) the corresponding phase space consists of 
\begin{equation}
C_{\pa}=\scT_Y^*X,
\label{HMV.16}\end{equation}
the restriction of the scattering cotangent bundle to the boundary.
Let $(x,y)$ be local coordinates near a boundary point, where $x$ is a
local defining function and $y$ restrict to coordinates on $Y$. Then the
forms $dx/x^2$ and $dy/x$ are a basis for the fibre $\scT^*_p X$ at each
point $p$ in the coordinate chart, and hence one may write an arbitrary
point $q \in \scT^* X$ in the form
\begin{equation}
q = -\nu \frac{dx}{x^2} + \sum_i \mu_i \frac{dy_i}{x},
\end{equation}
giving local coordinates $(x,y,\nu,\mu)$ on $\scT^* X$. 
The characteristic variety of $P$ at energy $\ev$ then is  
\begin{equation}
\Sigma (\ev )=\left\{|(\nu ,\mu )|_y^2+V(y)=\ev \right\}
\label{HMV.17}\end{equation}
where $|\cdot|^2_y$ is the metric function at the boundary. 

Again $C_{\pa}$ is naturally a contact manifold. The contact
structure arises as the `boundary value' of the singular symplectic
structure on $\scT^*X,$ corresponding to the fact that the associated
semiclassical model is the leading part, at the boundary, of the
Hamiltonian sytem defined by the energy. In terms of local
coordinates as above, the symplectic form,
arising from the identification $\scT^*X^\circ\simeq T^*X^\circ,$ is 
\begin{equation}
\widetilde\omega =d(\nu d\frac{1}{x} +\mu \cdot\frac{dy}x)
=(-d\nu+\mu \cdot dy)\wedge\frac{dx}{x^2}+\frac{d\mu \wedge dy}x,
\label{HMV.18}\end{equation}
and is hence singular at the boundary. Since the
vector field $x^2\pa_x$ is well-defined at the boundary, modulo multiples,
the contact form 
\begin{equation}
\alpha =\widetilde\omega (x^2\pa_x,\cdot)=-d\nu +\mu \cdot dy
\label{HMV.19}\end{equation}
defines a line subbundle of $T^*C_\pa$ which fixes the contact structure on
$C_{\pa}.$

The Hamilton vector field $H_p$ of $P$, regarded as a vector field on
$\scT^*X$, is of the form $x \tilde W$, where $\tilde W$ is tangent to $\pa
(\scT^*X)$. Here $W$ has the form 
\begin{equation*}
\tilde W = -2\nu x\pa_x  + W, \quad W \text{ a vector field on }\scT^*_{\pa
  X}X.
\end{equation*}
On $\Sigma(\ev) = \{ p=\ev \}$, we have $x^{-1}H_p = x^{-1}H_{p-\ev} =
H_{x^{-1}(p-\ev)} - (p-\ev)H_{x^{-1}} = H_{x^{-1}(p-\ev)}$. Since
$x^{-1}(p-\ev)$ is order $-1$, its Hamilton vector field $W_\ev$ restricted
to $\pa(\scT^*X)$ is the contact vector generated by the Hamiltonian $p-\ev$
on $\pa(\scT^*X)$. It is determined directly in terms of the contact form by 
\begin{equation}
d\alpha(\cdot,W_\ev)+\gamma \alpha=dp,\ \alpha (W_\ev)=p-\ev,
\label{HMV.20}\end{equation}
for some function $\gamma$.
The vector field $W_\ev$ is tangent to
$\Sigma(\ev)$ whenever the latter is smooth, as follows by pairing
\eqref{HMV.20} with $W_\ev.$ By the semiclassical model at energy $\ev$ we
mean the flow defined by $W = W_\ev$ on $\Sigma(\ev).$

\begin{lemma}\label{HMV.21} For any scattering metric and any real
$V\in\CI(X)$ the characteristic surface $\Sigma (\ev)$ 
is smooth whenever $\ev$ is not a critical value of $V_0.$ For regular
$\ev$ the critical set of the vector field $W$ on $\Sigma(\ev)$
is the union of the two radial sets 
\begin{equation}
\RP_\pm(\ev)=\{(y,\nu,\mu)\in\Sigma (\ev);
\nu =\pm\sqrt{\ev-V(p)},\ \mu =0,\ y=p\Mwhere d_YV(p)=0\}
\label{HMV.22}\end{equation}
and projection onto $Y$ gives bijections 
\begin{equation}
\pi:\RP_\pm(\ev)\longrightarrow \Crit(V)\cap\{y\in Y;V(y)<\ev\}
\label{HMV.63}\end{equation}
for each sign.
\end{lemma}

\noindent

\begin{proof} The vector field is fixed by \eqref{HMV.20} so in
any local coordinates $y$ in the boundary it follows from \eqref{HMV.19} that
\begin{equation}
d\mu \wedge dy(\cdot,W)=\gamma (d\nu -\mu dy)+dp.
\label{HMV.31}\end{equation}
Given a boundary point $q$ we may choose Riemannian normal
coordinates based at $q$ in the boundary. Thus the boundary metric is
Euclidean to second order, so $p=\nu ^2+|\mu|_y ^2+O(|y|^2|\mu|^2)+V(y)-\ev$
and we may easily invert \eqref{HMV.31} to find that at the fibre above $q$
\begin{equation}
W=2\mu\cdot\pa_y-(V'(y)+2\nu\mu)\cdot\pa_\mu+2|\mu|_y^2\pa_\nu +
O(|y||\mu|),\quad  \gamma =-2\nu.
\label{HMV.32}\end{equation}
This can only vanish when $\mu=0$ and $dV_0(\Cp)=0$ giving \eqref{HMV.22}.
\end{proof}

We remark that in the Euclidean setting, \eqref{HMV.22} amounts to
\begin{equation}
\RP_\pm(\ev )=\{(y,\zeta)\in\Sigma (\ev );\ y=\pm\zeta/|\zeta|,
\ y=\Cp\Mwhere d_YV(\Cp)=0\}. 
\label{HMV.138}\end{equation}

If $V_i=V\big|_{Y_i}$ is perfect Morse then the component of $\Sigma
(\ev)$ above $Y_i,$
\begin{equation}
\Sigma _i(\ev)\text{ is a }
\begin{cases}
\text{ sphere for }&
\kappa _i=\min(V\big|_{Y_i})<\ev <K_i=\max(V\big|_{Y_i}),\\
\text{ torus for }&\ev >K_i.
\end{cases}
\label{HMV.24}\end{equation}
In the general Morse case, $\Sigma _i(\ev)$ is a union of disjoint spheres for
non-critical $\ev<K_i$ and a torus for $\ev>K_i.$

In all cases if $\ev$ is not a critical value of $V\big|_{Y_i}$ then 
\begin{equation}
I_i(\ev )=\Sigma _i(\ev )\cap\{\nu =0\}
\label{HMV.25}\end{equation}
is a smooth curve (empty if $\ev<\kappa _i,$ generally with several
components) to which $W$ is transversal. The flow is symmetric under $\nu
\to-\nu,$ $\mu \to-\mu.$

\begin{prop}\label{HMV.29} Provided $\dim X=2$ and $V_0$ is Morse the
critical point $\Rp_{\pm}(\ev)=(\Cp,\pm\sqrt{\ev -V(\Cp)},0)$ in
\eqref{HMV.22} is
\begin{equation}
\begin{cases}
\text{a centre if $\Cp$ is a minimum and }\ev <\ev _{\Hess}(\Cp)\\
\text{a sink/source if $\Cp$ is a minimum and }\ev \ge\ev _{\Hess}(\Cp)\\
\text{a saddle if $\Cp$ is a maximum,}
\end{cases}
\label{HMV.30}\end{equation}
where $\ev_{\Hess}(\Cp)=V(\Cp)+2V''(\Cp).$
\end{prop}

\begin{proof} As already noted above, \eqref{HMV.32} is valid locally in
geodesic boundary coordinates. If $\Cp$ is a critical point for $V$ then $\mu=0$ at
the critical points $\Cp_\pm(\ev)$ of $W$ on $\Sigma (\ev)$ and we may use
local coordinates $y,\mu$ in the characteristic variety nearby, since $\nu
=\pm\sqrt{\ev-V(y)-\mu ^2}$ is smooth. In these coordinates the
linearization of $W$ inside $\Sigma(\ev)$ at the 
critical point is 
\begin{equation}
L=2[\mu \pa_y-(ay+\tilde\nu\mu)\pa_\mu],\ \text{where}\ V''(\Cp)=2a,\
\tilde\nu=\pm\sqrt{\ev -V(\Cp)}. 
\label{HMV.33}\end{equation}
The eigenvalues of $L$ are therefore the roots of
\begin{equation}
(\frac s 2)^2+\tilde\nu(\frac s 2)+a=0.
\label{HMV.34}\end{equation}
It is convenient for future reference to write the eigenvalues
$s_j$ as
\begin{equation}
s_j =-2\tilde\nu r_j;
\end{equation}
the $r_j$ thus satisfy
\begin{equation}
r^2-r+\frac{a}{\tilde\nu^2}=0.
\label{HMV.34p}\end{equation}

If $a<0,$ so $\Cp$ is a local maximum then the eigenvalues of $L$ are real,
and are given by 
\begin{equation}
r_1=\frac{1}{2}-\sqrt{\frac{1}{4}-\frac{a}{\tilde\nu^2}}
<0<1<r_2=\frac{1}{2}+\sqrt{\frac{1}{4}-\frac{a}{\tilde\nu^2}},
\label{HMV.36}\end{equation}
so the critical point is a saddle. If $a>0,$ so $\Cp$ is a local minimum,
then the discriminant
\begin{equation}
1-\frac{4a}{\tilde\nu^2} \text{ is }\begin{cases}\text{ negative }
&\Mif\ev<\ev_{\Hess}\\ 
\text{ zero }&\Mif \ev=\ev _{\Hess}\\
\text{ positive }&\Mif\ev>\ev _{\Hess}
\end{cases},\ \ev_{\Hess}=V(\Cp)+2V''(\Cp).
\label{HMV.35}\end{equation}
Correspondingly the eigenvalues are of the form $-2\tilde\nu r_j,$ $j=1,$
$2,$ where, respectively,
\begin{equation}\begin{split}
&r_1=\frac{1}{2}+ i\sqrt{\frac {a}{\tilde\nu^2}-\frac{1}{4}},
\ r_2=\frac{1}{2}- i\sqrt{\frac {a}{\tilde\nu^2}-\frac{1}{4}},\\\\
&r_1=r_2=1/2,\\
&0<r_1=\frac{1}{2}-\sqrt{\frac{1}{4}-\frac {a}{\tilde\nu^2}}<\frac{1}{2}<r_2
=\frac{1}{2}+\sqrt{\frac{1}{4}-\frac {a}{\tilde\nu^2}}<1
\end{split}\
\label{HMV.37}\end{equation}
and the critical point is a centre, a degenerate center, or a source/sink.
\end{proof}

\begin{rem}
If $\ev\neq\ev_{\Hess}(\Cp),$ a basis for the eigenvectors of $L$ is
\begin{equation}
\Eve_j=(s_j/2+\tilde\nu)\,dy+d\mu=\tilde\nu(1-r_j)\,dy+d\mu.
\label{eq:lin-evec}\end{equation}
If $\ev=\ev_{\Hess},$ $L$ has only one eigenvector,
$\frac 12 \tilde\nu\,dy+d\mu,$ but the generalized eigenspace (of eigenvalue
$-\tilde\nu$) is of
course two-dimensional.
\end{rem}

We also need a result about integral curves associated to the eigenvectors
at a sink.

\begin{prop}\label{prop:smooth-Leg}
Suppose that $\Rp \in \Min_\pm(\ev)$, and let $\Eve_j$
be the eigenvectors of $L$ as in \eqref{eq:lin-evec}. If $r_2/r_1>1$ is not an
integer, then there are two smooth one-dimensional submanifolds $L_j$
through $\Cp$ such that $W$ is tangent to $L_j$ in a neighbourhood of $\Rp$
and $\Eve_j\in N^*L_j.$ 
If $r_2/r_1=N$ is an integer, then there is one smooth Legendre curve
$L_1$ as above and another smooth curve, $L_2,$ such that $W=W_t+a \tilde
W,$ with $W_t$ tangent to $L_2,$ and where $a$ vanishes at $\Cp$ to order
$N$ \ie $a\in\mathcal{I}^{N},$ $\mathcal{I}$ being
the ideal of functions vanishing at $\Cp.$
\end{prop}

\begin{proof}
By a linear change of variables, $(y,\mu) \to (v_1,v_2)$, 
we may assume that $W$ takes the form
$$
W = r_1 v_1 \dbyd{}{v_1} + r_2 v_2 \dbyd{}{v_2} + O(|v|^2)
$$
near $v=0$. In these coordinates, the eigenvector $e_i$ corresponds to $dv_i$.
Suppose first that $r_2/r_1$ is not an integer. Then the eigenvalues are
nonresonant, and are both positive, so by the Sternberg linearization
theorem \cite{Sternberg} there is a smooth change of variables so that $W$
takes the form 
$$
W = r_1 v_1 \dbyd{}{v_1} + r_2 v_2 \dbyd{}{v_2}.
$$
Then we may take $L_i = \{ v_i = 0 \}$. In the resonant case, $r_2/r_1=N
\geq 2$, we still
obtain a normal form, but we are unable to remove resonant terms on the
right hand side. In this case, there is just one term, so we find
\cite{Sternberg}
that there is a smooth change of variables so that $W$ takes the form
$$
W = \frac1{N+1} v_1 \dbyd{}{v_1} + \frac{N}{N+1} v_2 \dbyd{}{v_2} + c v_1^{N}
\dbyd{}{v_2}. 
$$
Again, if we take $L_i = \{ v_i = 0 \}$ then we satisfy the conditions of
the theorem. 
\end{proof}

The corresponding result at saddle points is the stable/unstable
manifold theorem.

\begin{prop}\label{prop:smooth-Leg-saddle}
Suppose that $\Rp=\Cp_{\pm}(\ev)$ is a saddle, and let $\Eve_j$
be the eigenvectors of $L$ as in \eqref{eq:lin-evec}. Then
there are two smooth Legendre curves $L_j$ through $\Cp$
such that $W$ is tangent to $L_j$ and $\Eve_j\in N^*L_j.$
Moreover, if $\nu(\Rp)<0,$ then
$\nu\geq\nu(\Rp)$ on $L_1$ and $\nu\leq\nu(\Rp)$ on $L_2,$ while if
$\nu(\Rp)>0,$ then
$\nu\leq\nu(\Rp)$ on $L_1$ and $\nu\geq\nu(\Rp)$ on $L_2$. 
\end{prop}

\begin{proof}
The first part follows from the stable/unstable manifold theorem, while
the second part follows from \eqref{HMV.64}.  
\end{proof}

\begin{rem} In both cases, there are coordinates $(v_1,v_2)$ on
$\Sigma(\ev)$ near $q$ such that $W$ restricts to $r_2 v_2 \pa_{v_2}$ on
$L_1$ and $r_1 v_1 \pa_{v_1}$ on $L_2$ near $q$. 
\end{rem}

Below we make use of open neighbourhoods of the critical points which are
well-behaved in terms of $W.$

\begin{Def}\label{HMV.92} By a $W$-balanced neighbourhood of a
critical point $\Rp\in\Sigma(\ev),$ $\ev\nin\Cv(V),$ we shall mean an
open neighbourhood, $O$ of $\Rp$ in $C_{\pa}$ which contains no other
radial point, which meets $\Sigma(\ev)$ in a $W$-convex set (that is, each
integral curve of $W$ meets $\Sigma(\ev)$ in a single interval, possibly
empty) and is such
that the closure of each integral curve of $W$ in $O$ meets
$\nu=\nu(\Rp).$
\end{Def}

\begin{lemma}\label{HMV.187} Any critical point $\Rp\in\RP(\ev)$ for
$\ev>V(\pi(\Rp))$ has a neighbourhood basis of $W$-balanced open sets.
\end{lemma}

\begin{proof} In the case of sinks or sources this is clear, since
arbitrarily small balls around these that are convex with respect to $W$ are
$W$-balanced. In the case of saddle points we may use the stable 
manifold theorem. It suffices to suppose that $\nu(\Rp)>0.$ Taking a
neigbourhood of $L_2\cap\{\nu=\nu(\Rp)-\epsilon\},$
for $\epsilon>0$ but
small, a $W$-balanced neighbourhood is given by the $W$-flow-out of this set
in the direction of increasing $\nu$ strictly between $\nu=\nu(\Rp)-\epsilon$
and $\nu=\nu(\Rp)+\epsilon,$ together with the parts of $L_1$ and $L_2$
in $|\nu-\nu(\Rp)|<\epsilon.$
\end{proof}

Consider the global structure of the dynamics of $W.$ From \eqref{HMV.32}
it follows directly that 
\begin{equation}
W(\nu) = 2|\mu|_y^2 \geq 0, 
\text{ \ie $\nu$ is non-decreasing along integral curves of }W.
\label{HMV.64}\end{equation}
In fact, along maximally extended integral curves of $W,$ $\nu$ is only
constant if the curve reduces to a critical point. It is also easy to see
that $W$ has no non-trivial periodic orbits and every maximally extended
bicharacteristic $\gamma:\Real_t\to\Sigma(\ev)$ tends to a point in
$\RP_+(\ev) \cup \RP_-(\ev)$ as $t\to \pm\infty.$ Indeed,
$\lim_{t\to\pm\infty} \nu(\gamma(t))=\nu_\pm$ exists by the monotonicity
of $\nu,$ and any sequence $\gamma_k:[0,1]\to\Sigma(\ev)$,
$\gamma_k(t)=\gamma(t_k+t),$ $t_k\to +\infty,$ has a uniformly convergent
subsequence, which is then an integral curve $\tilde\gamma$ of $W.$ 
Then $\nu$ is constant along this bicharacteristic, hence $h$ is
identically $0.$ In view of the $\pa_\mu$ component of $W$ in
\eqref{HMV.32}, $\pa_y V$ is identically $0$ along $\tilde\gamma,$ hence
$y$ is a critical point of $V$ and the limit is a point in $\RP_+(\ev)\cup
\RP_-(\ev).$

We can now define the forward, resp.\ backward, bicharacteristic relation
as follows.

\begin{Def}\label{Def:Phi(K)-def}
The forward bicharacteristic relation $\Phi_+\subset\Sigma(\ev)\times
\Sigma(\ev)$ is the set of $(\xi',\xi)$ for which there exist
bicharacteristics $\gamma_j,$ $j=1,\ldots,N,$ $N\geq 1$,
and $T_-,T_+\in\Real,$ $T_-\leq T_+$,
such that $\gamma_1(T_-)=\xi',$ $\gamma_N(T_+)=\xi,$ and
$\lim_{t \to +\infty}\gamma_j(t)=\lim_{t\to-\infty}\gamma_{j+1}(t)$
for $j=1,2,\ldots,N-1$.

The backward bicharacteristic relation $\Phi_-$ is defined similarly, with
the role of $\xi$ and $\xi'$ reversed. The relations $\Phi_{\pm}$ depend on
$\ev$, but this is not indicated in notation. 
\end{Def}

In view of the previous observations, $\lim_{t\to\pm\infty}\gamma_j(t)
\in\RP_+(\ev)\cup\RP_-(\ev),$ so the content of the definition is to
allow $\xi$ and $\xi'$ to be connected by a chain of bicharacteristics
in the forward direction. Note that $\Phi_+$ is reflexive and transitive,
but not symmetric.

Also, if $\xi,\xi'\in\RP_+(\ev)\cup\RP_-(\ev),$ and there exists
a bicharacteristic $\gamma$ such that $\lim_{t\to+\infty}\gamma(t)=\xi$,
$\lim_{t\to-\infty}\gamma(t)=\xi',$ then $(\xi',\xi)\in\Phi_+.$ Indeed,
we can take $\gamma_1$ and $\gamma_3$ to be appropriately parameterized
constant bicharacteristics with image $\xi',$ resp.\ $\xi$,
and $\gamma_2=\gamma.$ Of course, $\gamma$ may be replaced by a chain
of bicharacteristics.

We recall
that the image of a set $K\subset\scT^*_Y X$ under $\Phi_+$ is
the set
\begin{equation}
\Phi_+(K)=\{\xi\in\scT^*_Y X;\exists\ \xi'\in K\ \text{s.t.}
\ (\xi',\xi)\in\Phi_+\}.
\end{equation}
We will be interested in the forward flow-out $\Phi_+(\RP_+(\ev))$ of the
outgoing radial set. Since $\nu>0$ on $\RP_+(\ev),$ and is increasing
along bicharacteristics, $\nu>0$ on $\Phi_+(\RP_+(\ev))$ as well.

\begin{lemma}
If $\Rp\in\RP_+(\ev)$ is a center or sink, then
$\Phi_+(\{\Rp\})=\{\Rp\}.$ On the other hand, if $\Rp$ is a saddle,
then with the notation of Proposition~\ref{prop:smooth-Leg-saddle},
$\Phi_+(\{\Rp\})$ is locally (in $\Sigma(\ev)$)
given by the local unstable curve $L_2,$ \ie by the unique smooth
Legendre manifold $L_2$ such that $\nu\geq \nu(\Rp)$ on $L_2$,
$\Rp\in L_2$ and $W$ is tangent to $L_2$.
\end{lemma}

\begin{rem}\label{rem:Phi(RP)}
Suppose that $\ev \notin \Cv(V)$. Then 
on $\Phi_+(\RP_+(\ev)),$ $\nu\geq a(\ev),$ and on 
$\Phi_-(\RP_-(\ev)),$ $\nu \leq -a(\ev),$ where
$a(\ev)=\min\{\nu(\Rp);\Rp\in\RP_+(\ev)\}>0.$
\end{rem}

\section{Blow-ups and resolution of singularities of
flows}\label{sec:blow-up} As indicated in the Introduction, to aid in the
description of the asymptotic behaviour of eigenfunctions it is very
convenient to modify the underlying manifold with boundary by blowing it up
at the minima of $V_0$ and with homogeneity $r$ with respect to the
boundary. We give a brief discussion of the blow-up of a stable critical
point of a vector field with respect to the vector field.

Suppose that $W$ is a real $\Cinf$ vector field on a manifold $M,$ and that
$o\in M,$ with $W(o)=0,$ is a linearly stable critical point. That is, the
eigenvalues of the linearization of $W$ at $o$ all have negative real
parts. We handle unstable critical points by changing the sign of the
vector field.

If $\Phi:M\times\Real_t\to M$ denotes the flow generated by $W,$ then $o$
has a neighbourhood $O$ such that for $p\in O,$ $\lim_{t\to+\infty}
\Phi(p,t)=o.$ There always exists a closed embedded submanifold,
diffeomorphic to a sphere, $S\subset O$ which is transversal to $W$. It
also necessarily satisfies
\begin{equation}
\lim_{t\to+\infty}\Phi(S,t)=\{o\}.
\label{HMV.190}\end{equation}
Then $\Phi$ restricts to a diffeomorphism of $S\times[0,+\infty)$ to
$O'\setminus\{o\}$ where $O'$ is the union of $S$ and the neighbourhood of
$o$ consisting of one of the components of $M\setminus S.$

We may compactify $S\times[0,+\infty)_t$ by embedding it as a dense subset 
\begin{equation*}
S\times[0,+\infty)_t\hookrightarrow S\times[0,1],\ (s,t)\longmapsto (s,e^{-t}).
\label{HMV.188}\end{equation*}
This makes $e^{-t}$ a defining function for a boundary hypsersurface. Using
the diffeomorphism $\Phi$ this compactifies $O'\setminus\{o\}$ to a
manifold with boundary. Finally then, replacing $O'$ as a subset of $M$ by
$S\times[0,1]$ using $\Phi$ to identify $O'\setminus\{o\}$ with
$S\times(0,1)$ we obtain a compactification of $M\setminus\{o\}$ which is
a compact manifold with boundary.

\begin{lemma}\label{HMV.189} The compact manifold with boundary
obtained by blow-up, with respect to a vector field $W,$ of a linearly
stable critical point, $\{o\},$ is independent of the transversal $S$ 
to $W,$ satisfying \eqref{HMV.190}, used to define it.
\end{lemma}
\noindent We may therefore denote the blown up manifold by $[M,\{o\}]_W.$

\begin{proof} Since there are spheres transversal to $W$ in any
neighbourhood of $\{o\}$ it suffices to show that the two manifolds obtained
using a transversal $S$ and a second transversal $\tilde S\subset O'$ are
naturally diffeomorphic where $O'$ is the component of $M\setminus S$
containing $o.$ Since it lies within $S,$ $S'$ is precisely the
image of $\Phi$ on $\{(s,T(s));s\in S\}$ for some smooth map $T:S\longrightarrow
(0,\infty).$ In terms of the compactification this replaces $e^{-t}$ by
$e^{-t}e^{-T(s)}$ which induces a diffeomorphism of the neighbourhoods of
the boundaries of the two compactifications.
\end{proof}

Thus the abstract manifold, $[M,\{o\}]_W,$ defined in this way is, as a
set, the union of $M\setminus\{o\}$ with an abstract sphere as
boundary. The boundary can be realized more concretely by using the flow of
$W$ to identify any two transversals to $W$ which are sufficiently close to
$\{o\}.$ The Lemma above shows that $[M,\{o\}]_W$ has a natural \ci\
structure as a compact manifold with boundary with interior which is
canonically diffeomorphic to $M\setminus\{o\}.$ The boundary hypersurface
is called the front face of the blow-up and is usually
denoted below as $\ff,$ or more precisely as $\ff([M;\{o\}]_W).$ Note that
it is \emph{not} generally the case that the natural map 
\begin{equation}
[M,\{o\}]_W\longrightarrow M,
\label{HMV.191}\end{equation}
under which $\ff\longmapsto \{o\},$ is smooth.

More generally, if $X$ is a smooth manifold with boundary and $W$ is a
smooth vector field on $X$ which is tangent to the boundary then we may use
the construction above to define $[X,\{o\}]_W$ even for a stable critical
point of $W$ on the boundary. Namely, we may simply extend $W$ across the
boundary to the double of $X$ and then observe that the closure of the
preimage of $X\setminus\{o\}$ in $[2X,\{o\}]_{\tilde W}$ is independent of
the extension $\tilde W$ of $W.$  

One particular use of this construction is to define inhomogeneous
blow-ups. Let $X$ be a 2-dimensional manifold with boundary, with $o\in\pa
X.$ Choose local coordinates $(x,y)$ near $o$ such that $x\geq 0$ is a
boundary defining function and $o$ is given by $x=0,$ $y=0.$ Let
$r\in(0,1)$ be a given homogeneity. We wish to blow up $o$ in $X$ so that
$y$ is homogeneous of degree $1$ and $x$ is homogeneous of degree $1/r$.
To do so consider the vector field $-W=r^{-1}x\pa_x+y\pa_y,$ and apply the
construction above, denoting the result $[X,\{o\}]_r;$ clearly it does
not depend on how $W$ is extended outside a neighbourhood of $o.$ Observe
that $y/x^r$ and $x/|y|^{1/r}$ are homogeneous of degree $0,$ \ie are
annihilated by $W,$ where they are bounded, so they can be regarded as
variables on the transversal $S.$ Thus, local coordinates on the blown-up
space, in the lift of the region $|y/x^r|<C$ are given by $y/x^r$ and
$x^r,$ while local coordinates in the lift of the region $x/|y|^{1/r}<C'$
are given by $x/|y|^{1/r}$ and $y.$

In fact the manifold $[X;\{o\}]_r$ does depend on the choice of $W,$ and so
on the choice of coordinates $(x,y).$ However the space of classical
conormal functions, with respect to the boundary, is defined independently
of choices. This is in essence because $x$ has homogeneity greater than
that of the boundary variable. More explicitly, any change of coordinates
takes the form $y'=a(x,y)y+b(x,y)x,$ where $x'=c(x,y)x,$ $c(0,0)>0$ and
$a(0,0)\neq 0.$ So, for example, $y'/(x')^r=(a/c^r)(y/x^r) +(b/c^r)x^{1-r}$
which is bounded (in $x\leq x_0,$ $x_0>0,$ $|y|\leq y_0,$ $y_0>0$) if and
only if $y/x^r$ is bounded (since $r\in(0,1)$). Such calculations show that
the blow-ups using $(x,y)$ and using $(x',y')$ coincide as topological
manifolds. Their $\Cinf$ structures are different as can be seen from the
appearance of terms involving $x^{1-r}$ in the transformation law. However,
it is easy to see that the class of conormal functions, and also that of
polyhomogeneous conormal functions, on $[X;\{o\}]_r,$ is well-defined,
independent of the choice of local coordinates used in the definition
(although the orders may need to be appropriately adjusted). Thus,
$[X;\{o\}]_r$ may be thought of as a `conormal manifold' rather than
as a $\Cinf$ manifold, \ie the algebra of smooth functions should be
replaced by the algebra of polyhomogeneous conormal functions as the basic
object of interest.

Such inhomogeneous blow-ups may be generalized to higher dimensions,
although there one needs extra structure to make it independent of the
choice of coordinates. Since we do not need this notion in the present
paper, we do not pursue it further.

\section{Limiting absorption principle}\label{sec:lim-abs}

In later sections we prove, by positive commutator methods, various results
of propagation of singularities type. Here we summarize, in four
theorems, the main conclusions concerning the limit of the resolvent on the
real axis and the behaviour of $L^2$ eigenfunctions. The scattering
wavefront set is discussed in \cite{Melrose43} as are the scattering
Sobolev spaces $\Hsc^{m,r}(X);$ the characteristic variety $\Sigma (\ev),$
scattering relation $\Phi_+$ and radial sets $\RP_{\pm}(\ev)$ are defined
in Section~\ref{S.Classical} above.

First, we recall H\"ormander's theorem on the propagation of singularities;
in the present setting this is \cite[Proposition~7]{Melrose43}.

\begin{thm}[Propagation of singularities]
\label{thm:prop-sing}
For $u\in\dist(X)$, and any real $m$ and $r$,
\begin{equation*}\begin{gathered}
\WFsc^{m,r-1}(u)\setminus\WFsc^{m,r}\left((P-\ev)u\right)
\text{ is a union of maximally extended } \\
\text{ bicharacteristics of $P-\ev$ in
$\Sigma(\ev)\setminus\WFsc^{m,r}\left((P-\ev)u\right)$}.
\end{gathered}\end{equation*}
Similarly, if $u_t$, for $t \in G \subset (0,1]$ is a bounded family in
$\Hsc^{M,R}(X)$, then 
\begin{equation*}\begin{gathered}
\WF_{\scl,L^\infty(G_t)}^{m,r-1}(u_t)\setminus
\WF_{\scl,L^\infty(G_t)}^{m,r}\left((P-(\ev+it))u_t\right)  
\text{ is a union of maximally } \\
\text{ forward-extended bicharacteristics of $P-\ev$ in }
\Sigma(\ev)\setminus\WF_{\scl,L^\infty(G_t)}^{m,r}\left((P-\ev)u\right).
\end{gathered}\end{equation*}
\end{thm}

\begin{thm}[Decay of outgoing eigenfunctions]\label{thm:unique}
If $\ev\nin\Cv(V)$ then the space $E_{\pp}(\ev)$ is finite dimensional, is
contained in $\dCinf(X)$, and may be characterized by 
\begin{multline}
E_{\pp}(\ev) = \{u\in\dist(X);(P-\ev)u=0\Mand\\
\WFsc^{m,l}(u)\cap \RP_+(\ev)=\emptyset\Mforsome l>-1/2\}.
\label{HMV.68}\end{multline}
The set of eigenvalues $\sigma_{\pp}(P) = \{ \ev; E_{\pp}(\ev) \neq \{ 0 \}
\}$ is discrete 
outside $\Cv(V)$ and is contained in $(-\inf_XV,\sup\Cv(V)]$. 
\end{thm}

As an immediate consequence of Theorem~\ref{thm:unique},

\begin{equation}
E_{\ess}^{-\infty}(\ev)=\left\{u\in\CmI(X);(P-\ev)u=0,\ \langle u,v\rangle
=0\ \forall\ v\in E_{\pp}(\ev)\right\}
\end{equation}
is well defined for $\ev\nin\Cv(V)$ as are spaces such as 
\begin{equation*}
\Hsc^{m,r}(X)\ominus E_{\pp}(\ev)=\{f\in\Hsc^{m,r}(X);
\langle f,\phi\rangle=0\ \forall\ \phi\in E_{\pp}(\ev)\}.
\end{equation*}

\begin{thm}[Limiting absorption principle]\label{thm:lim-abs}
The resolvent 
\begin{equation}
R(\ev+it)=(P-(\ev+it))^{-1},\ 0\ne t\in\bbR,\ \ev\nin\Cv(V)
\label{HMV.38}\end{equation}
extends continuously to the real axis, \ie $R(\ev\pm i0)$ exist, as
bounded operators 
\begin{multline}
\Hsc^{m,r}(X)\ominus E_{\pp}(\ev)\longrightarrow
\Hsc^{m+2,l}(X)\ominus E_{\pp}(\ev),\
\forall\ m\in\bbR,\ r>1/2,\ l<-1/2.
\label{HMV.65}\end{multline}
\end{thm}

\begin{thm}[Forward propagation]\label{thm:res-WF}
For $\ev\nin\Cv(V)$ and $f\in\dCinf(X)\ominus E_{\pp}(\ev),$ 
\begin{equation}
\WFsc(R(\ev+i0)f)\subset\Phi_+(\RP_+(\ev))\subset\{\nu\geq 0\}
\label{HMV.66}\end{equation}
and $R(\ev+i0)$ extends by continuity to
\begin{multline}
R(\ev+i0):\{v\in\dist(X)\ominus E_{\pp}(\ev);
\WFsc(v)\cap\Phi_-(\RP_-(\ev))=\emptyset\}
\longrightarrow \\
\dist(X),
\label{HMV.67}\end{multline}
with wavefront set bound given by 
\begin{equation}
\WFsc(R(\ev+i0)v)\subset\WFsc(v)\cup\Phi_+(\WFsc(v)\cap\Sigma(\ev))
\cup\Phi_+(\RP_+(\ev)).   
\label{eq:WF-bound}\end{equation}
\end{thm}

We next show that various statements in the Introduction are direct
consequences of these results.

\begin{prop}\label{P.Closed.range}
For $\ev\nin\Cv(V)$ the space $(P-\ev)\dCI(X)$ is closed
in $\dCI(X)$ and for $\ev\in\sigma_{\pp}(P)\setminus\Cv(V)$,
$(P-\ev)\dCI(X)$ is $L^2$-orthogonal to $E_{\pp}(\ev)$ and 
\begin{equation}
R(\ev+i0)f=R(\ev-i0)f\quad \forall\ f\in(P-\ev)\dCI(X).
\label{HMV.70}\end{equation}
\end{prop}

\begin{proof} By assumption $\ev\nin\Cv(V)$ so
$E_{\pp}(\ev)\subset\dCI(X).$ The self-adjointness of $P$ gives
\begin{equation*}
\langle (P-\ev)v,\phi\rangle=
\langle v,(P-\ev)\phi\rangle=0\Mif \phi 
\in E_{\pp}(\ev)
\end{equation*}
so $(P-\ev)\dCI(X)$ is $L^2$ orthogonal to $E_{\pp}(\ev).$

Suppose $u$ is in the closure of $(P-\ev)\dCI(X)$ in $\dCI(X).$ Thus
$u_j=(P-\ev)v_j\to u$ with $v_j\in\dCI(X).$  By the finite dimensionality of
$E_{\pp}(\ev)$ we may assume that $v_j\perp E_{\pp}(\ev).$ Now consider
$w_j=R(\ev+i0)u_j$ which exist by Theorem~\ref{thm:lim-abs} and by
Theorem~\ref{thm:res-WF} satisfy $\WFsc(w_j)\subset\Phi_+(\RP_+(\ev)).$ Thus, 
\begin{equation*}
(P-\ev)(v_j-w_j)=0 \Mand \WFsc(v_j-w_j)\subset \Phi_+(\RP_+(\ev))
\end{equation*}
so, by Theorem~\ref{thm:unique}, $v_j-w_j\in E_{\pp}(\ev),$ to which both
terms are orthogonal, so $v_j=w_j=R(\ev+i0)u_j.$
Since $u_j\to u$ in $\dCI(X),$ $v_j\to R(\ev+i0)u$ in $\CmI(X).$
The same argument shows that $v_j\to R(\ev-i0)u$ in $\CmI(X),$ so
$v=R(\ev+i0)u=R(\ev-i0)u.$ By \eqref{HMV.66} applied with
each sign, $\WFsc(v)$ is disjoint from both $\{\nu\le0\}$ and
$\{\nu\ge0\}$ and hence is empty. Thus $v\in\dCI(X),$ and since
$u=(P-\ev)v,$ $(P-\ev)\dCI(X)$ is closed in $\dCI(X);$ this also proves
\eqref{HMV.70}.
\end{proof}

\begin{cor}
For $\ev\nin\Cv(V),$ the space $\dCI(X)/(P-\ev)\dCI(X)$ is a Fr\'echet
space, as is $(\dCI(X)\ominus E_{\pp}(\ev))/(P-\ev)\dCI(X)$.
\end{cor}

An abstract parameterization of generalized eigenfunctions, which has
little specific to do with our particular problem, may be obtained via the
two terms of \eqref{HMV.10}. Namely, let
\begin{equation}
\RR_{\ess,\pm}^{\infty}(\ev)=
[R(\ev\pm i0)(\dCI(X)\ominus E_{\pp}(\ev ))]/(\dCI(X)\ominus E_{\pp}(\ev ))
\label{HMV.137}\end{equation}
be the range of the incoming ($-$), resp.\ outgoing ($+$),
boundary value of the resolvent
acting on Schwartz functions, modulo Schwartz functions.
It is now immediate that the map
\begin{equation*}
(P-\ev):\RR_{\ess,\pm}^{\infty}(\ev)\longrightarrow
(\dCI(X)\ominus E_{\pp}(\ev))/(P-\ev)\dCI(X)
\end{equation*}
is an isomorphism (for each sign). This induces a Fr\'echet topology
on $\RR_{\ess,\pm}^{\infty}(\ev)$.

\begin{prop}\label{prop:Sp}
For $\ev\nin\Cv(V)\cup\sigma_{\pp}(P)$ the nullspace of
$$
\Sp(\ev)=\frac1{2\pi i}\left(R(\ev+i0)-R(\ev-i0)\right):
\dCI(X)\longrightarrow\SEEF(\ev)
$$
is $(P-\ev)\dCI(X);$ if $\ev\in\sigma_{\pp}(P)\setminus\Cv(V),$ the same
remains true for the map
$(2\pi i)^{-1}\left(R(\ev+i0)-R(\ev-i0)\right):\dCI(X)\ominus
E_{\pp}(\ev)\to\SEEF(\ev)$, for $\ev \in \sigma_{\pp}(P) \setminus \Cv(V)$.
\end{prop}

\begin{proof}
By \eqref{HMV.70}, $(P-\ev)\dCI(X)$ lies in the nullspace of $\Sp(\ev).$
Conversely, suppose that $f\in\dCI(X),$ $\Sp(\ev)f=0,$ so
$v=R(\ev+i0)f=R(\ev-i0)f.$ Since
$\WFsc(R(\ev+i0)(f))$ and $\WFsc(R(\ev-i0)(f))$ are
disjoint ($\nu<0$ on one, $\nu>0$ on the other), both of these are
empty, so $v\in\dCI(X),$ and $(P-\ev)v=f$ is in $(P-\ev)\dCI(X)$ as claimed.
\end{proof}

\begin{cor}
If $\ev\nin\Cv(V)\cup\sigma_{\pp}(P)$ the spaces 
$\SEEF(\ev),$ $\RR_{\ess,+}^{\infty}(\ev)$,
$\RR_{\ess,-}^{\infty}(\ev)$ and $\dCI(X)/(P-\ev)\dCI(X)$
are all isomorphic; this remains true for
$\ev\in\sigma_{\pp}(P)\setminus\Cv(V)$ with $\dCI(X)/(P-\ev)\dCI(X)$
replaced by $(\dCI(X)\ominus E_{\pp}(\ev))/(P-\ev)\dCI(X)$.
\label{cor:SEEF-RR}\end{cor}

Explicitly, the isomorphisms
\begin{equation}
I_\pm(\ev):\SEEF(\ev)\to\RR_{\ess,\pm}^{\infty}(\ev)
\label{HMV.60}\end{equation}
are the maps
\begin{equation*}\begin{split}
&I_+(\ev):u=[R(\ev+i0)-R(\ev-i0)]f\mapsto R(\ev+i0)f,\\
&I_-(\ev):u=[R(\ev+i0)-R(\ev-i0)]f\mapsto -R(\ev-i0)f,\ f\in\dCI(X).
\end{split}\end{equation*}
Although $f$ is not well-defined in $\dCI(X),$ it is
in $\dCI(X)/(P-\ev)\dCI(X),$ hence the image of $R(\ev+i0)f$ in
$\RR_{\ess,+}^{\infty}(\ev)$ is also well defined.
By Theorem~\ref{thm:res-WF},
\begin{equation*}
\WFsc(I_\pm(\ev)u)\subset\Phi_\pm(\RP_\pm(\ev)).
\end{equation*}
Thus, $I_-(\ev)$ maps $u\in\SEEF(\ev)$ to its `incoming part',
while $I_+(\ev)$ maps $u$ to its outgoing part. We define the
abstract scattering matrix as the map
\begin{equation}
S(\ev):\RR_{\ess,-}^{\infty}(\ev)\to \RR_{\ess,+}^{\infty}(\ev),
\quad S(\ev)=I_+(\ev)I_-(\ev)^{-1}.
\end{equation}
In Section~\ref{sec:Morse} we identify $\RR_{\ess,\pm}^{\infty}(\ev)$
geometrically. We use this identification in Section~\ref{sec:pairing}
to extend $S(\ev)$ to a unitary operator on a Hilbert space, giving
rise to the definition of the S-matrix in the introduction in
\eqref{eq:S-def}.

We now turn to the relationship between $\SEEF(\ev)$ and $\EEF(\ev)$
which has two facets. On the one hand, as we show below,
$\SEEF(\ev)$ is dense in $\EEF(\ev)$ in the topology of $\dist(X)$,
indeed in a stronger sense. On the other hand they are dual to
each other with respect to a form of `Green's pairing' that we describe
in Section~\ref{sec:pairing}.

By Theorem~\ref{thm:res-WF}, $\Sp(\ev)$ may be applied to any distribution
which is $L^2$ orthogonal to $E_{\pp}(\ev),$ provided
its scattering wavefront set is contained in $\{|\nu| <\epsilon \}$ for
$\epsilon >0$ sufficiently small, namely $\ep<a(\ev),$
\begin{equation*}
a(\ev)=\min\{\nu(\Rp);\Rp\in\RP_+(\ev)\}>0.
\end{equation*}

\begin{prop}\label{prop:Sp-parameterization}
If $\ev\nin\Cv(V)$ and $\epsilon >0$ then every
$u\in \EEF(\ev)$ is of the form $\Sp(\ev)f$ for some
$f\in\Dist(X)$ with $\WFsc(f)\subset\{|\nu|<\epsilon\}$ and $f\perp
E_{pp}(\ev).$ 
\end{prop}

\noindent Of course if $\ev\nin\sigma_{\pp}(P)$ the orthogonality
condition is void.

\begin{proof}
Suppose first that $\ev\nin\sigma_{\pp}(P)$.

Let $A \in \Psisc^{-\infty,0}(X)$ be a pseudodifferential operator with
boundary 
symbol given by $\psi_1(\sigma_\pa (P - \ev)) \psi_2(\nu/a(\ev)),$ where
$\psi_1(t)$ is supported in $|t| < \ep$ and $\psi_2(t)$ is supported in $t
> 1/4,$ and equal to $1$ when $t > 1/2.$ 
Then, by Theorem~\ref{thm:res-WF}, we may apply $R(\ev - i0)$ to $(P -
\ev)(\Id - A)u$ and $R(\ev + i0)$ to $(P - \ev)Au.$ Since
\begin{equation*}
(P - \ev) \big( (\Id - A)u - R(\ev - i0)(P - \ev)(\Id -
A)u \big) = 0,
\end{equation*}
and this function has scattering wavefront contained in $\nu \leq a(\ev)/2,$ we
conclude using Theorem~\ref{thm:unique} that $(\Id - A)u = R(\ev - i0)(P -
\ev)(\Id-A)u.$ Similarly, $Au = R(\ev + i0)(P - \ev)Au$.
Since $(P-\ev)Au=[P,A]u=-(P-\ev)(\Id-A)u$,
we may take $f = (2\pi i) [P,A] u$.

If $\ev\in \sigma_{\pp}(P),$ only a minor modification is necessary.
Namely, we consider $\Pi A$ in place of $A$ where $\Pi$ is orthogonal
projection off $E_{\pp}(\ev).$ Note that
$\Id-\Pi\in\Psisc^{-\infty,\infty}(X).$ Then the rest of the argument
goes through.
\end{proof}

In fact, this proof shows more.

\begin{prop}\label{HMV.192}
Suppose $\ev\nin\Cv(V)$ and that $A\in\Psisc^{-\infty,0}(X)$ has
\begin{equation*}
\WFscp(A)\subset\{\nu>-a(\ev)\}\Mand \WFscp(\Id-A)\subset\{\nu<a(\ev)\}
\end{equation*}
where $0<a(\ev)=\min\{\nu(\Rp);\Rp\in\RP_+(\ev)\};$ if $\Pi$ is the
orthogonal projection off $E_{\pp}(\ev)$ then 
\begin{equation}
2\pi i\Sp(\ev)\Pi [P,A]:\dist(X)\longrightarrow \dist(X)
\label{eq:Sp-comm}\end{equation}
is continuous, with range in $E_{\ess}^{-\infty}(\ev)$ and the restriction
of \eqref{eq:Sp-comm} to $E_{\ess}^{-\infty}(\ev)$ is the identity
map.
\end{prop}

\begin{proof}
By assumption $-a(\ev)<\nu<a(\ev)$ on $\WFscp(A)\cap\WFscp(\Id-A)$ so this
is also holds on $\WFscp([P,A]).$ The continuity of \eqref{eq:Sp-comm} follows.
The final statement is a consequence of the proof of the preceeding
proposition, namely that, given $u\in E_{\ess}^{-\infty}(\ev)$ one can take
$f=2\pi i \Pi[P,A]u$ and conclude that $u=\Sp(\ev)f.$ 
\end{proof}

This result allows us to characterize the topology on $\SEEF(\ev)$
three different ways. Let $\cT_1$ be the topology on $\SEEF(\ev)$
as in Corollary~\ref{cor:SEEF-RR}, induced by
$\Sp(\ev):(\dCI(X)\ominus E_{\pp}(\ev))\to\SEEF(\ev).$ Fix any $A$ as in
Proposition~\ref{HMV.192} and let $\cT_2$ be the
topology on $\SEEF(\ev)$ induced by the topology of $\dCI(X)$ on
$\SEEF(\ev)$ via the map $[P,A]:\SEEF(\ev)\to\dCI(X).$ Finally let
$r<-1/2,$ and let $B\in\Psisc^{-\infty,0}(X)$ be such that 
\begin{equation}
B \text{ is elliptic at }\{\nu=0\}\cap\Sigma(\ev), \text{ with }
\WFscp(B)\subset\{|\nu|<a(\ev)\},
\label{eq:B}\end{equation} and let
$\cT_3$ be the weakest topology on $\SEEF(\ev)$ stronger than both
the $x^r L^2_{\scl}(X)$ topology on $\SEEF(\ev)$ and the
topology induced by the map $B:\SEEF(\ev)\to\dCI(X).$

\begin{prop}
The topologies $\cT_1, \cT_2, \cT_3$ are equivalent. In particular,
$\cT_2$ and $\cT_3$ are independent of choices.
\end{prop}

\begin{proof}
The equivalence of $\cT_1$ and $\cT_2$ follows from
$2\pi i\Sp(\ev)\Pi[P,A]=\Id$ on $\SEEF(\ev).$
By the propagation of singularities, $\cT_3$ is stronger that $\cT_2.$
Finally, by the forward propagation, $\cT_1$ is stronger than $\cT_3.$
\end{proof}

We can also topologize
\begin{equation}
E^s_{\ess}(\ev)
=\big\{u\in\EEF(\ev);\ \WFsc^{0,s-1/2}(u)\cap\{\nu=0\}=
\emptyset\big\},\ s\in\bbR,
\label{HMV.157p}\end{equation}
similarly. Note first that
\begin{equation*}
E^s_{\ess}(\ev)\subset x^rL^2_\scl(X)\Mfor r<\min(-\frac{1}{2},s-\frac{1}{2}).
\end{equation*}
For $s\in\bbR$, we let $\cT_2^s$ be the topology
on $E^s_{\ess}(\ev)$ induced by $[P,A]:E^s_{\ess}(\ev)\to x^{s+1/2}
L^2_{\scl}(X),$ so
\begin{equation*}
\|u\|_{\cT_2^s}=\|u\|_2=\|[P,A]u\|_{x^{s+1/2}L^2_{\scl}(X)}.
\end{equation*}
In the same setting, let $\cT_3^s$ be the topology
on $E^s_{\ess}(\ev)$ induced by the $x^r L^2_{\scl}(X)$ topology on
$E^s_{\ess}(\ev)$, for $r$ with $r < -1/2$ and $r<s-1/2,$
and by the operator $B$ from
\eqref{eq:B}, that is by the norm
\begin{equation*}
\|u\|_{\cT_3^s}=\|u\|_3=\left(\|u\|^2_{x^rL^2_{\scl}(X)}
+\|Bu\|^2_{x^{s-1/2}L^2_{\scl}(X)}\right)^{1/2}.
\end{equation*}
For $s>0$, we deduce the following as above.

\begin{prop}\label{prop:topologies-s}
For each $s>0$, the topologies $\cT_2^s$ and $\cT_3^s$ are equivalent
and are independent of all choices, making
$E^s_{\ess}(\ev)$ into a Hilbert space.

The inclusion map $E^s_{\ess}(\ev)\hookrightarrow E^{s'}_{\ess}(\ev),$
with norm $\|.\|_3,$ is bounded for $s'\leq s.$
\label{prop:equiv-norms}\end{prop}

We consider $s=0$ when we study Green's pairing in Section~\ref{sec:pairing}.

From the density of $\dCI(X)$ in the subspace of $\Dist(X)$ with wave front
set in $\{|\nu|< \epsilon\}\cap\Sigma(\ev),$ with the topology of
H\"ormander, and the continuity of $\Sp(\ev)$ from this space to
$\Dist(X)$ we conclude (using Proposition~\ref{prop:Sp-parameterization})
that the smooth eigenfunctions, $E_{\ess}^\infty(\ev)$, 
are dense in $\EEF(\ev)$.

More explicitly, this can be seen by considering $u_r=\phi(x/r)u,$ $r>0$
where $\phi\in\Cinf(\bbR)$ is identically $1$ on $[1,+\infty)$ and $0$ on
$[0,1/2).$ Then $\Sp(\ev)\Pi [P,A]u_r\in\SEEF(\ev)$ for $r>0,$ and $u_r\to
u$ in $\dist(X),$ hence
\begin{equation*}
\SEEF(\ev)\ni 2\pi i\Sp(\ev)\Pi [P,A]u_r\to u
\end{equation*}
in
$\dist(X).$ If $u\in E_{\ess}^s(\ev),$ $s>0,$ then $[P,A]u_r\to[P,A]u$ in
$x^{s+1/2}L^2_{\scl}(X),$ hence
\begin{equation*}
[P,A]2\pi i\Sp(\ev)\Pi [P,A]u_r\to [P,A]2\pi i\Sp(\ev)\Pi [P,A]u=[P,A]u
\end{equation*}
in $x^{s+1/2}L^2_{\scl}(X).$
More generally, suppose that $u\in E_{\ess}^s(\ev),$ $s\in\bbR,$ and let
$Q\in\Psisc^{0,0}(X)$ be such that $\WFscp(Q)\subset\{|\nu|<a(\ev)\}$,
$\WFscp(\Id-Q)\cap\WFscp([P,A])=\emptyset$. Let $p$ be such that
$p<s-1/2$, $p<-1/2$.
Then $[P,A]u_r\to[P,A]u$ in $x^{s+1/2}L^2_{\scl}(X)$ shows that
$2\pi i\Sp(\ev)\Pi Q [P,A]u_r\to 2\pi i\Sp(\ev)\Pi Q [P,A]u$ in
$x^pL^2_\scl(X)$ and $(\Id-Q)[P,A]u_r\to(\Id-Q)[P,A]u$ in $\dCI(X)$ so
$2\pi i\Sp(\ev)\Pi [P,A]u_r\to u$ in $x^pL^2_\scl(X)$, and, with $B$
from \eqref{eq:B}, $2\pi iB\Sp(\ev)\Pi [P,A]u_r\to Bu$ in
$x^{s-1/2}L^2_\scl(X)$. We have thus proved the following.

\begin{cor}\label{prop:density} For each $\ev\nin\Cv(V),$ $\SEEF(\ev)$
is dense in $\EEF(\ev)$ in the topology of $\Dist(X).$
Moreover, for each $\ev\nin\Cv(V),$ $s\in\bbR,$ $\SEEF(\ev)$
is dense in $E_{\ess}^s(\ev)$ in the topology $\cT_3^s$.
\end{cor}

\section{Propagation of singularities}\label{S.Propagation}
In this section, we derive wavefront set bounds on eigenfunctions,
or more generally of solutions of $(P-\ev)u=f$,
and use these to prove the results of the previous section. Let us first
consider $\ev\in\Cx\setminus[\inf V_0,+\infty)$. Then the symbol $p-\ev$ of $\sigma_\pa(P-\ev)$ never
vanishes, so $P - \ev$ has a 
parametrix in the scattering calculus, \ie there is $G(\ev) \in
\Psisc^{-2,0}(X)$ such that $E(\ev) = (P - \ev) G(\ev) - \Id \in
\Psisc^{-1,1}(X).$ Thus $E(\ev)$ is compact on $L^2.$ By analytic
Fredholm theory we conclude that $P$ has discrete spectrum in
$(-\infty,\inf V_0).$ Also ellipticity of $p-\ev$ at the boundary
implies that $\WFsc(u)$ is empty if $(P - \ev)u = 0$, so $u \in
\CIdot(X)$. 

In the rest of the section we consider $\ev > \inf V_0$, and use
positive commutator estimates to derive wavefront set bounds. First we
recall some standard results in the scattering calculus. 
Since $P \in \Psisc^{*,0}(X)$, the Hamilton vector field
$H_p$ of $P$ vanishes to first order at $x=0$, we define the scattering
Hamilton vector field $\scH_p = x^{-1} H_p$, which is smooth and
nonvanishing up to the boundary. In fact, $\scH_p = W -2\nu x \pa_x$ in
local coordinates near the boundary, where $W$ is the contact vector field
from Section~\ref{S.Classical}. If $A \in
\Psisc^{*,-l-1}(X)$ is a scattering pseudodifferential operator, with
symbol $x^{-l-1}a$, then the symbol of the commutator $i[P, A]$ is 
\begin{equation}
\sigma_{\pa}(i[P,A]) =  H_p(x^{-l-1}a) = x^{-l} \big( Wa + 2(l+1)\nu a
\big).
\end{equation}
We also recall the notion of wavefront set with respect to a family of
functions, $\WF_{\scl,L^\infty(G_t)}^{m,l}(u_t)$, where $t \in G$ is a
parameter. For a single function $v$, the statement $q \notin
\WFsc^{m,l}(v)$ is equivalent to the existence of $A \in \Psisc^{m,-l}(X)$
which is elliptic at $q$ and such that $Av \in L^2(X)$, while the statement 
$q \notin \WF_{\scl,L^\infty(G_t)}^{m,l}(u_t)$ is equivalent to the
existence of $A$ as above such that $\| Au_t \|_2$ is {\it uniformly}
bounded.

We begin with a Lemma describing the basic structure of
positive commutator estimates.  

\begin{lemma}\label{lemma:pos-comm}
Suppose $a \geq 0 \in x^{-l-1}\Cinf(\scT^*X)$ satisfies
\begin{equation}
\scH_p a=-x^{-l-1}b^2+x^{-l-1}e,
\label{eq:commutator}\end{equation}
with
$a^{1/2}b,$ $b$ and $e\in\Cinf(\scT^*X).$ Suppose
\begin{equation*}
A\in\Psisc^{-\infty,-l-1}(X),\ B\in\Psisc^{-\infty,-l-1/2}(X)\Mand
E\in\Psisc^{-\infty,-2l-1}(X)
\end{equation*}
have
principal symbols $a,$ $a^{1/2}b$ and $ae$ respectively, with
$\WFscp(E)\subset\supp e$, and  assume that 
\begin{equation}
a \leq C b^2 \text{ on } \{ e = 0 \} \text{ for some
  } C>0. 
\label{eq:abe-cond}\end{equation}
Then there exists $F\in \Psisc^{-\infty,-2l}(X)$
such that for all $u\in \Hsc^{*,l+1/2}(X)$,
\begin{equation}
\|Bu\|^2+2t\|Au\|^2\leq |\langle u,Eu\rangle|+|\langle u,Fu\rangle|
+2|\langle u,A^*A(P-(\ev+it))u\rangle|.
\label{eq:pos-comm-16}\end{equation}
In particular, if $G\subset(0,1]_t,$ and $u_t,$ $t\in G$,
is a (not necessarily bounded) family in $\Hsc^{*,l+1/2}(X)$ such that 
\begin{multline*}
\begin{gathered}
\WF_{\scl,L^\infty(G_t)}^{*,l}(u_t)\cap\supp a=\emptyset,\\
\WF_{\scl,L^\infty(G_t)}^{*,l+3/2}((P-(\ev+it))u_t)\cap\supp a
=\emptyset
\\
\Mand
\WF_{\scl,L^\infty(G_t)}^{*,l+1/2}(u_t)\cap\supp e=\emptyset,
\end{gathered}\\
\text{ then }\WF_{\scl,L^\infty(G_t)}^{*,l+1/2}(u_t) \subset \{ b = 0 \}.
\end{multline*}
The assumptions $t>0,$ $G\subset (0,1]$ can be replaced by $t<0,$ $G\subset
[-1,0)$ if we arrange that
$0\leq a\in x^{-l-1}\psi_0(p)\Cinf(\scT^*X)$ satisfies
$$
\scH_p a=+
x^{-l-1}b^2+x^{-l-1}e,
$$
with
$a^{1/2}b,$ $b$ and $ae\in\Cinf(\scT^*X)$.
\end{lemma}

\begin{proof}
Let $A\in
\Psisc^{-\infty,-l-1}(X)$ be
a quantization of $a,$ \ie be such that $\sigma_{\pa,-l-1}(A)=x^{l+1}a$.
From the symbol calculus, and \eqref{eq:commutator},
\begin{equation}
ix^{l+1/2}[A^*A,P]x^{l+1/2}= 2x^{l+1/2}(B^*B+E+F)
x^{l+1/2},
\label{eq:pos-comm-7}\end{equation}
where
\begin{equation}
\begin{gathered}
F\in\Psisc^{-\infty,-2l}(X),\ \WFscp(F)\subset\supp (x^{l+1}a),\\
B\in\Psisc^{-\infty,-l-1/2}(X)\Mand\sigma_{\pa,-l-1/2}(B)=ba^{1/2}.
\end{gathered}
\label{eq:pos-comm-9}\end{equation}

For $v\in\Hsc^{*,l+1}(X),$ $t\in\Real$,
\begin{equation}
\langle v,i[A^*A,P]v\rangle
=2\im\langle v,A^*A(P-(\ev+it))v
\rangle-2t\|A v\|^2.
\end{equation}
Now if $u\in\Hsc^{*,l+1/2}(X),$ consider this for $v=u_r=\phi(x/r)u$ where
$\phi\in \Cinf(\Real)$ is identically $1$ on $[2,+\infty)$ and identically
$0$ on $[0,1].$ Writing 
$(P-(\ev+it))u_r=\phi(x/r)(P-(\ev+it))u+[P,\phi(x/r)]u,$ observe that
$A^*A[P,\phi(x/r)]\in\Psisc^{-\infty,\infty}(X)$ is, for $r\in(0,1]$,
uniformly bounded in $\Psisc^{-\infty,-2l-1}(X),$ and indeed converges
to $0$ strongly as an operator
$\Hsc^{*,l+1/2}(X)\longrightarrow\Hsc^{*,-l-1/2}(X),$ while $u_r\to u$ in
$\Hsc^{*,l+1/2}(X).$ Taking the limit $r\to 0,$ we deduce that
\begin{equation}
\langle u,i[A^*A,P]u\rangle
=2\im\langle u,A^*A(P-(\ev+it))u
\rangle-2t\|A u\|^2.
\label{eq:pos-comm-12}\end{equation}
Combining this with \eqref{eq:pos-comm-7}, we conclude that
\eqref{eq:pos-comm-16} holds.

Now consider the uniform statement for the family. Suppose first that 
\begin{equation*}
\WF_{\scl,L^\infty(G_t)}^{*,l+2}((P-\ev+it)u_t)
\cap\supp a=\emptyset.
\end{equation*}
Since $t>0,$ the second term on the left in \eqref{eq:pos-comm-16} can be
dropped. Since $u_t$ is microlocally uniformly bounded in $\Hsc^{0,l}(X)$
on $\supp a$ and $f_t=(P-(\ev+it))u_t$ is microlocally uniformly bounded in
$\Hsc^{0,l+2}(X)$ by our assumption, it follows that the right hand side of
\eqref{eq:pos-comm-16} stays bounded as $t\to 0.$ Thus, $B u_t$ is
uniformly bounded in $L^2_\scl(X).$ 
This proves that $\supp b$
is disjoint from $\WF^{*,l+1/2}_{\scl,L^\infty(G_t)}(u_t).$

In the general case, when we only assume that
$\WF_{\scl,L^\infty(G_t)}^{*,l+3/2}((P-\ev+it)u_t)
\cap\supp a=\emptyset,$ we need to remove $u$ from the last term
of \eqref{eq:pos-comm-16}.
To do so, we apply Cauchy-Schwarz to the last term and estimate it by
\begin{multline*}
2\|x^{1/2}Au_t\|_2\,\ \|x^{1/2}A(P-(\ev+it))u_t\|_2\leq\\
\ep \|x^{1/2}Au_t\|_2^2+\ep^{-1}\|x^{-1/2}A(P-(\ev+it))u_t\|_2^2,
\end{multline*}
for any $\ep>0$. By the assumption
$\WF_{\scl,L^\infty(G_t)}^{*,l+1/2}(u_t)\cap\supp e=\emptyset,$
we can find $Q \in
\Psisc^{0,0}(X)$ such that $\WFsc'(Q) \cap
\WF_{\scl,L^\infty(G_t)}^{*,l+1/2}(u_t) = \emptyset$
(hence $Qu_t$ is uniformly bounded in $\Hsc^{*,l+1/2}(X)$), while
$\supp e\cap\WFscp(\Id-Q)=\emptyset,$ so $e=0$ on
$\WFsc'(\Id - Q)$. By assumption \eqref{eq:abe-cond}, we have $a \leq C
b^2$ on $\WFsc'(\Id - Q)$. Then we may estimate
\begin{multline}
\ep \langle x^{1/2}Au_t, x^{1/2}Au_t \rangle = 
\ep \Big( \langle x^{1/2}AQ^2u_t, x^{1/2}Au_t \rangle  + \langle
x^{1/2}A(\Id-Q^2)u_t, x^{1/2}Au_t \rangle \Big) \\
\leq \ep C' \Big( \| x^{1/2} AQu_t \|_2^2 + \| Bu_t \|_2^2 + |\langle \tilde
Fu_t,u_t \rangle| \Big),
\end{multline}
for some $C'>0$ and $\tilde F \in \Psisc^{*,-2l}(X)$. Choose $\epsilon =
(2C')^{-1}$. Then we may remove the $\| Bu_t \|_2^2$ term from the right hand
side of \eqref{eq:pos-comm-16} by replacing $\|Bu_t\|^2$ by
$1/2 \|Bu_t\|^2$ on the left hand side, obtaining
\begin{equation}\begin{split}
&\frac{\|Bu_t\|^2}{2}+t\|Au_t\|^2\\
&\leq C'' \Big( \| x^{1/2} AQu_t \|_2^2  +
|\langle u_t,Fu_t\rangle| +\ep^{-1}\|x^{-1/2}A(P-(\ev+it))u_t\|^2 \Big).
\label{eq:pos-comm-32}\end{split}\end{equation}    
The argument of the preceeding paragraph now applies and completes the proof.
\end{proof}

\begin{rem} The conclusions of Lemma~\ref{lemma:pos-comm} hold if we
weaken the assumption that $u \in \Hsc^{*,l+1/2}(X)$ to 
$u \in \Hsc^{*,k}(X)$, for some $k\in\Real$, and
$\WFsc^{*,l+1/2}(u) \cap \supp a = \emptyset$. We may then apply operators
$C \in \Psisc^{*,-l-1/2}(X)$ with essential support contained in $\supp a$
to $u$. This is done by the prescription
$$
Cu = CQu + C(\Id-Q)u,
$$
where $Q\in \Psisc^{0,0}(X)$ is such that $\WFsc'(Q) \cap \WFsc'(C) =
\emptyset$, and $\WFsc'(\Id-Q) \cap \WFsc^{*,l+1/2}(u) = \emptyset$. This
is independent of the choice of pseudodifferential operator
$Q$. This observation allows the iterative use of the Lemma to gain
regularity, as is done below. 
\end{rem}

If $\nu$ has a fixed sign, then so does the coefficient $-2\nu x$ of
$\pa_x$ in the Hamilton vector field of $p.$ This enables us to prove the
second part of the conclusion of the lemma above without requiring the
existence of a $t$-dependent family $u_t$ satisfying various conditions,
\ie we can work directly with $t=0.$

\begin{lemma}\label{lemma:pos-comm-reg}
Suppose that $a,b,e,$ etc., are as in Lemma~\ref{lemma:pos-comm}, and suppose
that $\nu>0$ on $\supp a$ and that $u\in \Hsc^{*,l}(X),$ with
\begin{multline*}
\begin{gathered}
\WF_{\scl}^{*,l+3/2}((P-\ev)u)\cap\supp a=\emptyset\Mand\\
\WF_{\scl}^{*,l+1/2}(u)\cap\supp e=\emptyset.
\end{gathered}\\
\text{ Then } \WF_{\scl}^{*,l+1/2}(u)\subset\{ b=0\}.
\end{multline*}
The same conclusion holds if $a\geq 0,$ $\nu<0$ on $\supp a$ and
$$
\scH_p a=x^{-l-1}\psi_0(p)b^2+x^{-l-1}\psi_0(p)e,
$$
with
$a^{1/2}b,$ $b$ and $e\in\Cinf(\scT^*X).$
\end{lemma}

\begin{rem}
This theorem will be used to analyze functions whose wavefront set is
concentrated near a radial point of $\scH_p$. 
Away from the radial points, the sign of $\nu$ is irrelevant
and Theorem~\ref{thm:prop-sing} tells us how the wavefront set of $u$
propagates. 
\end{rem}

\begin{proof}
For $r>0$ set $\chi_r=(1+r/x)^{-1/2}=(\frac{r+x}{x})^{-1/2}.$ Let $A$, $B$,
$E$ and $F$ be as before, and let 
$A_r=A\chi_r,$ $B_r=B\chi_r$ and $E=E\chi_r$.
We repeat the arguments of the previous proposition with
$A_r$ in place of $A.$
Now,
$\scH_p\chi_r^2=-2\nu (r/x)(1+r/x)^{-2}=-c_r^2$ is negative in $\nu>0$
and so on $\supp a.$ Thus, for $r>0,$
\begin{equation}
ix^{l+1/2}[A^*_rA_r,P]x^{l+1/2}= x^{l+1/2}(B_r^*B_r+A^*C_r^*C_rA+ E_r+F_r)
x^{l+1/2},
\label{eq:pos-comm-7p}\end{equation}
where $C_r$ is the operator multiplication by $c_r$,
$F_r\in\Psisc^{-\infty,-2l+1}(X)$ for $r>0,$ is uniformly bounded
in $\Psisc^{-\infty,-2l}(X)$ as $r\to 0$ and $\WFscp(F_r)\subset\supp (x^{l+1}a)$
uniformly as $r\to 0.$

Note that $A_r\in\Psisc^{-l-1/2}(X)$ for $r>0$,
and it is uniformly bounded in $\Psisc^{-l-1}(X)$ as $r\to 0.$ Since
$u\in\Hsc^{*,l}(X)$ by our assumption, the proof of
\eqref{eq:pos-comm-12} (with $l$ replaced by $l+1/2$) is applicable. 
Combining this with \eqref{eq:pos-comm-7p}, we deduce that for $r>0$
\begin{equation}
\|B_r u\|^2+\|C_r Au\|^2\leq
|\langle u,F_r u\rangle|+|\langle u,E_r u\rangle|
+2|\langle u,A_r^*A_r(P-\ev)u\rangle|.
\label{eq:pos-comm-16p}\end{equation}
We now drop the term $\|C_r Au\|^2$, and apply Cauchy-Schwarz to the last
term as in the previous proof, obtaining  
\begin{equation}
\frac{\|B_ru\|^2}{2} \leq C'' \Big( |\langle E_ru,u\rangle| +
|\langle u,\tilde F_r u\rangle| +\ep^{-1}\|x^{-1/2}A_r(P-\ev)u\|^2 \Big).
\label{eq:pos-comm-16pp}\end{equation}
The right hand side of
\eqref{eq:pos-comm-16pp} remains uniformly bounded as $r\to 0,$ so
we conclude that $Bu\in \Hsc^{0,0}(X)$, which implies the conclusion of the
Lemma. 
\end{proof}

For basic spectral and scattering theory, partially microlocal estimates
suffice, 
\ie one does not need full microlocal results. These partially microlocal
results (namely, microlocal in the $\nu$ variable, global in $y$ and $\mu$)
are closely related to the two-body resolvent estimates of Isozaki and Kitada
\cite{Isozaki-Kitada:Microlocal},
and the corresponding many-body estimates of G\'erard, Isozaki and
Skibsted \cite{GIS}. 

\begin{prop}\label{prop:fine-reg-t-sg}
Suppose that, for $t\in G\subset (0,1],$ some $l > -1,$ $\nu_0<0$ 
and $r>-1/2$, the function $(P-(\ev+it))u_t=f_t$ satisfies
\begin{gather*}
\WF^{*,l}_{\scl,L^\infty(G_t)}(u_t)\cap\{\nu<\nu_0\}=\emptyset\Mand\\
\WF^{*,r+1}_{\scl,L^\infty(G_t)}(f_t)\cap\{\nu<\nu_0\}=\emptyset,
\label{HMV.193}\end{gather*}
then $\WF^{*,r}_{\scl,L^\infty(G_t)}(u_t)\cap\{\nu<\nu_0\}=\emptyset.$
\end{prop}

\begin{proof}
The key point is that
$\scH_p$ is radial at $\RP(\ev)=\RP_+(\ev)\cup \RP_-(\ev)$. That is, the
tangential component $W$ vanishes there, so we need to exploit the normal
component, which is $-2\nu x \partial_x$. We obtain a positive
commutator by using a weight $x^{-l-1}$, since $-l-1 < 0$. 
Fix $\ep>0,$ $\delta>0,$ $M>0,$
and let
\begin{equation}
a=x^{-l-1}\chi(\nu)\psit(x)\psi_0(p)\geq 0
\label{eq:fine-reg-t-sg-8}\end{equation}
where $\psit\in\Cinf_c(\Real)$ is identically $1$ near $0$ and
is supported in a bigger neighbourhood of $0$ (it is simply a cutoff
near $\pa X),$ $\psi_0\in\Cinf_c(\Real;[0,1])$ supported in $(\ev-\delta,
\ev+\delta),$ $\chi\in\Cinf_c(\Real;[0,1])$ is identically one on 
$(-\infty,\nu_0-3\ep)$, vanishes on $(\nu_0-\ep, \infty)$,
$\chi'\leq 0$,
and $\chi(t)=e^{M/(t-(\nu_0-\ep))}$ on $(\nu_0-2\ep,\nu_0-\ep).$ Then
\begin{equation}\begin{split}
\scH_p a=2((l+1)\nu\chi(\nu)+|\mu|^2_y \chi'(\nu))x^{-l-1}\psi_0(p)
=-b^2x^{-l-1},
\end{split}\label{eq:fine-reg-t-sg-16}\end{equation}
and $b$ is $\Cinf$ by the construction of $\chi$. Thus, if $r<l+3/2$, then 
Lemma~\ref{lemma:pos-comm} gives us the result (with $\nu_0$ replaced by
$\nu_0-\ep$). Otherwise, from the Lemma
we gain a power of $x^{1/2}$ in the wavefront set estimate for $u_t$, \ie,
we have 
$$
\WF^{*,l+1/2}_{\scl,L^\infty(G_t)}(u_t)\cap\{\nu<\nu_0 - \ep\}=\emptyset.
$$
We may now repeat the argument with $l$ replaced by $l+1/2$. 
Iterating the argument a finite number of times (cf. the remark after the
proof of Lemma~\ref{lemma:pos-comm}), then sending $\ep \to 0$,
we obtain our estimate. 
\end{proof}

A similar argument works in $\nu \geq 0,$ but this time we only need to
assume boundedness in $\Hsc^{*,l}(X)$ for an arbitrary $l\in\Real$ in
this region, \ie we do not need $l\in(-1,-1/2)$. On the other hand, 
we are only able to deduce
boundedness in $\Hsc^{*,r}(X)$ for $r<-1/2.$ In addition, we need to assume
that the desired boundedness, \ie in $\Hsc^{*,r}(X)$ for $r<-1/2,$ holds
for $\nu$ near $\nu_0>0.$

\begin{prop}\label{prop:rough-reg-t-sg}
Suppose $f_t=(P-(\ev+it))u_t,$ where $t\in G\subset (0,1],$ 
$\nu_0>0,$ $l\in\Real$ and $r<-1/2$. If we have 
\begin{multline*}
\WF^{*,l}_{\scl,L^\infty(G_t)}(u_t)\cap\{\nu\geq\nu_0\}=\emptyset,\\
\WF^{*,r}_{\scl,L^\infty(G_t)}(u_t)\cap\{\nu=\nu_0\}=\emptyset \Mand\\
\WF^{*,r+1}_{\scl,L^\infty(G_t)}(f_t)\cap\{\nu\geq\nu_0\}=\emptyset,\\
\text{ then }
\WF^{*,r}_{\scl,L^\infty(G_t)}(u_t)\cap\{\nu\geq\nu_0\}=\emptyset.
\label{HMV.194}\end{multline*}
\end{prop}

\begin{proof} The proof is nearly identical to that of
Proposition~\ref{prop:fine-reg-t-sg}. We may assume $l<-1$ and $l<r$ 
without loss of
generality. Also, the estimate is clear away from $\Sigma(\ev)$ by elliptic
regularity, so we may restrict our attention to a small neighbourhood $O$ of
$\Sigma(\ev)$ whose closure $\overline{O}$ is compact. Take $a$ as
in \eqref{eq:fine-reg-t-sg-8}, where now 
$\chi'\geq 0,$ $\chi$ vanishes on
$(-\infty,\nu_0-\ep]$, is identically $1$ on $[\nu_0+\ep,+\infty)$, and
given by $e^{M/(\nu-(\nu_0+\ep))}$ on $(\nu_0, \nu_0+\ep)$. We
assume that
$\ep>0$ is so small that $\WF^{*,r}_{\scl,L^\infty(G_t)}(u_t) \cap \{ \chi'
\neq 0 \} \cap \overline{O} = \emptyset$ and 
$\WF^{*,r+1}_{\scl,L^\infty(G_t)}(f_t) \cap \{ \chi \neq 0 \} \cap
\overline{O} = \emptyset$. Such a choice is possible since the wave front
set is 
closed and $\overline{O}$ is compact. Now $l+1<0$ and $\nu>0,$ so the sign
of the first term in \eqref{eq:fine-reg-t-sg-16} is unchanged, while the
sign of the second term is reversed. Hence, the second term becomes an
error term, $e,$ which satisfies \eqref{eq:abe-cond}. By construction we
have  
$\WF^{*,r}_{\scl,L^\infty(G_t)}(u_t) \cap \supp e = \emptyset$. Thus, we
can apply Lemma~\ref{lemma:pos-comm}, following the argument of the previous
Proposition, to complete the proof. 
\end{proof}

Versions of these results with no parameter $t$ can also be obtained. In
this case, there is no distinction between $\nu > 0$ and $\nu < 0$.

\begin{prop}\label{prop:t-ind} Let $\nu_0 < 0$ and $l >-1/2$. Suppose that
$Pu = f$, where  
\begin{equation*}
\WF_{\scl}^{*,l}(u)\cap \{ \nu \leq \nu_0 \}=\emptyset \Mand 
\WF_{\scl}^{*,r+1}(f)\cap \{ \nu \leq \nu_0 \}=\emptyset. 
\end{equation*}
Then 
\begin{equation*}\WF_{\scl}^{*,r}(u)\cap \{ \nu \leq \nu_0 \} = \emptyset.
\end{equation*}
A similar result holds with $\nu_0>0$ with all
inequalities $\nu\leq\nu_0$ replaced by $\nu\geq\nu_0.$
\end{prop}

\begin{proof} The proof
uses the symbol $a$ from Proposition~\ref{prop:fine-reg-t-sg}, and
proceeds via a regularization argument similar to
Lemma~\ref{lemma:pos-comm-reg}. Now, however, $\scH_p a$ and $\scH_p\chi_r^2$
have opposite signs, so the two terms cannot be treated separately.
Instead, we use the argument of the lemma with $a$ replaced by $x^{1/2}a$
and $\chi_r=(1+r/x)^{-1/2}$ replaced by $x^{-1/2}\chi_r=(x+r)^{-1/2}.$
Now $\scH_p (x^{1/2}a)$ and $\scH_p(x\chi_r^2)$ have the same sign, so
the proof of Lemma~\ref{lemma:pos-comm-reg} is applicable, and completes
the proof of the proposition. Note that the power of the factor of $x$ in
$x^{1/2}a$ is $x^{-l-1/2};$ this is where we use that $l>-1/2.$
\end{proof}

The combination of Theorem~\ref{thm:prop-sing} and 
Proposition~\ref{prop:t-ind} leads to 

\begin{prop}\label{HMV.195} Suppose that $l>-1/2,$ $\ev_0\nin\Cv(V).$
Then there exist operators $F_j\in\Psisc^{-1,1}(X),$ $j=1,\ldots,N$,
and constants $C>0,$ $\delta>0,$ such that for $u\in\Hsc^{m,l}(X)$,
$(P-\ev)u\in\Hsc^{m,l+1}(X),$ $\ev\in(\ev_0-\delta,\ev_0+\delta)$,
$$
\|u\|_{\Hsc^{m,l}(X)}\leq \sum_{j=1}^N\|F_j u\|_{\Hsc^{m,l}(X)}
+C\|(P-\ev) u\|_{\Hsc^{m,l+1}(X)}.
$$
\end{prop}

\begin{proof} Choose a small $\nu_0 > 0$ so that the Hamiltonian flow never
vanishes 
for $|\nu| \leq 2\nu_0$, on the chacteristic surfaces $\Sigma(\ev)$, $\ev
\in(\ev_0-\delta,\ev_0+\delta)$. 
We then decompose the identity operator in the form $\Id = Q_0 + Q_1 + B_+ +
B_-$, where all operators are in $\Psisc^{0,0}(X)$, and are microlocalized
as follows: $Q_0$ is microsupported in $|\nu| \leq 2\nu_0$, $B_\pm$ are of
the form that arise in the proof of Proposition~\ref{prop:t-ind}, thus
microsupported in $\{ \pm \nu > \nu_0 \}$, and $Q_1$
is microsupported away from the characteristic sets $\Sigma(\ev)$. 
Since the Hamilton flow does not vanish on the microsupport of $Q_0$, we
may write $Q_0$ as a commutator with $P-\ev$, modulo lower order terms, and
hence obtain an estimate
$$
\|Qu\|_{\Hsc^{m,l}(X)}\leq C \Big( \|(P-\ev) u\|_{\Hsc^{m,l+1}(X)} +
|\langle Eu, u \rangle| \Big),
$$
where $E \in \Psisc^{2m,-2l}(X)$ satisfies $\WFsc'(E)
\subset \{ |\nu| \geq \nu_0 \}$. We add to this a sufficiently
large multiple of the estimates for $\| B_\pm u\|_{\Hsc^{m,l}(X)}$,
obtained from the proof of Proposition~\ref{prop:t-ind}, and the estimate
for $\| Q_1u \|_{\Hsc^{m,l}(X)}$ obtained from ellipticity. This allows one
to absorb the error term $|\langle Eu, u \rangle|$ above. 
\end{proof}

This immediately yields the usual results on the spectrum of $P.$ For
the last statement we also need a unique continuation result.

\begin{prop}\label{prop:eigenvalues}(See Theorem~\ref{thm:unique}) For
$\ev\nin\Cv(V)$ the 
space $W$ of $L^2$ eigenfunctions with eigenvalue $\ev$ is a finite
dimensional subspace of $\CIdot(X)$. The pure point spectrum, 
$\sigma_{pp}(P)$, is disjoint from
$(\sup\Cv(V),+\infty)$ and can only accumulate at $\Cv(V).$
\end{prop}

\begin{proof}
Let $F_k$ be as in Proposition~\ref{HMV.195}. Since these operators are
compact, there exists a finite dimensional subspace $W_0$ of $W$ such that
for $u\in W$ orthogonal to $W_0,$ $\|F_k u\|_{\Hsc^{0,0}(X)}\leq
\|u\|_{\Hsc^{0,0}(X)}/2N.$ Suppose that $u\in W,$ $\|u\|_{\Hsc^{0,0}(X)}=1.$
Applying the proposition with $\ev_0=\ev$ we deduce that
$$
1=\|u\|_{\Hsc^{0,0}(X)}\leq \sum_{k=1}^N\|F_k u\|_{\Hsc^{0,0}(X)}
\leq 1/2,
$$
which is a contradiction. Hence $W=W_0,$ so $W$ is finite dimensional.

More generally, let $\tilde W$ be the sum of the $L^2$ eigenspaces
of $P$ with eigenvalue $\tilde\ev\in[\ev_0-\delta/2,\ev_0+\delta/2]$.
By compactness of the $F_k,$ there is a finite dimensional subspace $\tilde
W_0$ of $\tilde W$ such that for $u\in \tilde W$ orthogonal
to $\tilde W_0,$ $\|F_k u\|_{\Hsc^{0,0}(X)}\leq \|u\|_{\Hsc^{0,0}(X)}/2N.$
Again applying Proposition~\ref{HMV.195} with $\tilde\ev$ in place of $\ev$
in the last term, $\|(P-\tilde\ev) u\|_{\Hsc^{m,l+1}(X)}=0,$
shows that $\tilde W=\tilde W_0,$ so $\tilde W$ is also finite dimensional.

That an $L^2$ eigenfunction $u$ is in $\CIdot(X)$ follows from the
following Proposition. The statement that $\sigma_{\pp}(P)$ is disjoint from 
$(\sup\Cv(V),+\infty)$ follows from the unique continuation theorem
of Froese and Herbst \cite{FroExp}, as generalized to scattering metrics
in \cite{Vasy:Propagation-2}.
\end{proof}

\begin{prop}\label{prop:unique}(Theorem~\ref{thm:unique})
Suppose that $u\in\dist(X),$ $(P-\ev)u=0,$ $\WFsc^{m,l}(u)\cap
\RP_-(\ev)=\emptyset$ for some $l>-1/2$. Then $u\in\dCinf(X).$
\end{prop}

\begin{proof}
This is essentially Isozaki's uniqueness theorem from \cite{Isozaki1},
Lemma~4.5, proved the same way in this setting. From
Theorem~\ref{thm:prop-sing} and Proposition~\ref{prop:t-ind}, we have 
\begin{equation}
\WFsc(u)\subset \{ \nu \geq a(\ev) > 0 \}
\Mand u\in\Hsc^{m',l'}(X)\ \forall\ m'\in\Real,\ l'<-1/2,
\label{eq:unique-1}\end{equation}
where
$a(\ev)$ is as in Remark~\ref{rem:Phi(RP)}. 
Let $l\in(-1/2,0)$ and let $\phi\in\Cinf(\Real)$ be
$0$ on $(-\infty,1]$ and $1$ on $[2,\infty).$ For $r>0$ let
\begin{equation}
\chi_r(x)=r^{-2l-1}\int_0^{x/r}\phi^2(s)s^{-2l-2}\,ds.
\end{equation}
Thus, $\chi_r\in\Cinf_c(\interior(X))$ and
\begin{equation}
x^2\partial_x\chi_r(x)=x^{-2l}\phi^2(x/r).
\end{equation}
Note that $\chi_r$ is not uniformly bounded in $\Psisc^{m',l'}(X)$ for any
$m'$ and $l',$ but $x^2\partial_x\chi_r$ is. Since $\chi_r$ enters the
commutator $[\chi_r(x),P]$ only via $x^2\partial_x\chi_r,$ that boundedness
will suffice for us.
Now, by \cite[Equation 3.7]{Melrose43}
\begin{equation}
\Delta=(x^2 D_x)^2+i(n-1)x^3 D_x+x^2\Delta_h+x^3\Diffb^2(X),
\end{equation}
so
\begin{equation}
-i[\chi_r(x),P]=2x^{-2l}\phi^2(x/r)(x^2 D_x)+F'_r
\end{equation}
where $F'_r$ is bounded in $\Psiscc^{1,-2l+1}(X).$ Let $\psi\in\Cinf_c
(\Real)$ be supported close to $\ev$ and be identically $1$ near $\ev.$
Let $\rho\in\Cinf(\Real)$
be $0$ on $(-\infty,a(\ev)/4)$ and be $1$ on $(a(\ev)/2,\infty)$. Let
$B\in\Psisc^{-\infty,0}(X)$ and $E\in\Psisc^{-\infty,0}(X)$ have symbols
$\sigma_{\pa,0}(B)=\sqrt{\nu} \rho(\nu) \psi(p) \in \CI$ and
$\sigma_{\pa,0}(E)=\nu(1-\rho^2(\nu))\psi^2(p)$. Then
\begin{equation}
-i[\chi_r(x)\psi^2(P),P]=2x^{-l}\phi(x/r)(B^2+E)\phi(x/r)x^{-l}+F_r
\end{equation}
where $\WFscp(E)$ is disjoint from $\{ \nu \geq a(\ev)/2 \}$ 
and $F_r$ is bounded in $\Psisc^{-\infty,-2l+1}(X).$ Now, for $r>0$
\begin{equation}
\langle u,[\chi_r(x),P]u\rangle
=-2i\im\langle u,\chi_r(x)(P-\ev)u\rangle=0.
\label{eq:unique-4}\end{equation}
Hence,
\begin{equation}
\|x^{-l}\phi(x/r)Bu\|^2\leq |\langle x^{-l}\phi(x/r)u,E x^{-l}\phi(x/r)u
\rangle|+|\langle u,F_r u\rangle|.
\label{eq:unique-5}\end{equation}
In view of \eqref{eq:unique-1} the right hand side stays bounded as $r\to 0,$ 
so we conclude that $x^{-l}Bu\in L^2_{\scl}(X).$ Then, by \eqref{eq:unique-1},
it follows that $u\in\Hsc^{\infty,l}(X).$ Since $l>-1/2,$ we 
conclude from Theorem~\ref{thm:prop-sing} and Proposition~\ref{prop:t-ind}
that $u\in\dCinf(X).$
\end{proof}
As mentioned above, $\chi_r(x)$ is not bounded in $\Psisc^{m',l'}(X)$ for any
$m'$ and $l',$ so the place where we have used the assumption
$(P-\ev)u=0$ is really in the elimination of the term on the right hand side of
\eqref{eq:unique-4} from the right hand side of \eqref{eq:unique-5}.

\begin{thm}(Theorem~\ref{thm:lim-abs}) The resolvent $R(\ev+it),$ $t>0,$
$\ev\nin\Lambda =\Cv(V)\cup\sigma_{pp}(P)$ extends continuously to the real
axis, \ie $R(\ev+i0)$ exists, as a bounded operator
$\Hsc^{m,r}(X)\longrightarrow  \Hsc^{m+2,l}(X)$ for any $r>1/2$ and $l<-1/2.$
\end{thm}

\begin{proof} Let $f \in \Hsc^{m,r}(X)$, for $r>1/2$, let $-1 < l < -1/2$,
and suppose that $u_t=R(\ev+it)f$ is not bounded in $\Hsc^{m+2,l}(X)$ as
$t\to 0.$ We can then take a sequence $t_j\to 0$ with
$\|u_{t_j}\|_{\Hsc^{m+2,l}(X)}\to\infty.$ Now consider
$v_{t_j}=u_{t_j}/\|u_{t_j}\|_{\Hsc^{m+2,l}(X)},$ so $v_{t_j}$ remains
bounded in $\Hsc^{m+2,l}(X),$ and $(P-(\ev+it))v_{t_j}\to 0$ in $\Hsc^{m,r}(X).$
We may then pass to a convergent subsequence in $\Hsc^{m+2-\delta,l-\delta}(X)$,
$\delta>0,$ with limit
$v.$ Thus $(P-\ev)v=0,$ and by Proposition~\ref{prop:fine-reg-t-sg},
$\RP_-(\ev)\cap\WF^{*,l'}_{\scl,L^\infty(G_t)}(v_t)=\emptyset$ for
some $l'>-1/2$. By the $t$-dependent part of Theorem~\ref{thm:prop-sing} on
the propagation of singularities
$\WF^{*,l'}_{\scl,L^\infty(G_t)}(v_t)\subset\{\nu>0\},$ hence the same
holds for $v.$
The uniqueness result, Theorem~\ref{thm:unique} then shows that $v=0.$
However this contradicts the fact that
$\|v_{t_j}\|_{\Hsc^{m+2,l'}(X)}=1.$ Hence, there exists $C>0$
such that $\|u_t\|_{\Hsc^{m+2,l}(X)}\leq C$ for all $t.$
Then, by Proposition~\ref{prop:fine-reg-t-sg} and the $t$-dependent part
of Theorem~\ref{thm:prop-sing},
$\WF^{*,l'}_{\scl,L^\infty(G_t)}(u_t)\subset\{\nu>0\}.$

Now suppose that $t_j\to 0.$ By decreasing $m+2$ and $l,$ we can pass to a
convergent subsequence with limit $u,$ which then satisfies
$(P-\ev)u=f,$ and $\WFsc^{*,l'}(u)\subset\{\nu>0\}.$ Taking any other
sequence $t'_j$ and any convergent subsequence, we obtain another distribution
$u'$ with the same properties, hence by considering their difference
and using Theorem~\ref{thm:unique},
$u=u'.$ If now $t''_j$ is a sequence such that $t''_j\to 0$ and
$\|u_{t''_j}-u\|_{\Hsc^{m+2,l}(X)}\geq \ep>0,$ then taking a subsequence
of $t''_j$ with $u_{t''_j}$ converging in $\Hsc^{m+2-\delta,l-\delta}(X)$
to some $u'',$ we deduce that $u=u'',$ which is a contradiction.
This shows convergence of $R(\ev+it)f$ as $t\to 0,$ and similar arguments
easily give boundedness as an operator, and continuity in $\ev$ as well.
\end{proof}

We next refine this theorem by combining it with the propagation
estimates. To do so we need a microlocal version of
Proposition~\ref{prop:fine-reg-t-sg}, which we state afterwards.

\begin{thm}(Theorem~\ref{thm:res-WF})
If $\ev\nin\Lambda$ and $f\in\dCinf(X)$ then $\WFsc(R(\ev+i0)f)\subset
\Phi_+(\RP_+(\ev)).$ Moreover $R(\ev+i0)$ extends by continuity 
to $v\in\dist(X)$ with $\WFsc(v)\cap\Phi_-(\RP_-(\ev))=\emptyset$,
and for such $v,$ $R(\ev+i0)v$ satisfies \eqref{eq:WF-bound}.
\end{thm}

\begin{proof} Let $u = R(\ev+i0)f$. 
It follows from Proposition~\ref{prop:fine-reg-t-sg} that
$\WFsc(u)\subset \{\nu \geq 0\}.$ Then Theorem~\ref{thm:prop-sing} shows
that $\WFsc(u)\subset\Phi_+(\RP_+(\ev)).$ Indeed, given $q\in \Sigma(\ev)$,
let $\gamma$ be the bicharacteristic through $q.$ Then
$q_0=\lim_{t\to-\infty} \gamma(t)\in \RP_+(\ev)\cup \RP_-(\ev).$ If $q_0\in
\RP_+(\ev),$ there is nothing to prove, since then $q\in\Phi_+(\RP_+(\ev)).$ If
$q_0\in \RP_-(\ev),$ there exist points on $\gamma$ where $\nu$ is
negative, hence these are not in $\WFsc(u).$ But then by
Theorem~\ref{thm:prop-sing}, $q\nin\WFsc(u).$

A duality argument immediately gives the extension of $R(\ev+i0)$ to the
class of distributions described above. The wave front set bound is a
consequence of Theorem~\ref{thm:prop-sing} and the following fully microlocal
family version of Proposition~\ref{prop:fine-reg-t-sg}, stated below.
\end{proof}

\begin{prop}\label{prop:micro-fine-prop-t}
Suppose that $q\in \RP_-(\ev),$ $f_t=(P-(\ev+it))u_t,$ for $t\in
G\subset(0,1],$ and $O$ is a $W$-balanced neighbourhood (see
Definition~\ref{HMV.92}) of $q$ then for
$r>-1/2,$ and $l>-1$,
\begin{multline*}
O\cap\WF^{*,r+1}_{\scl,L^\infty(G_t)}(f_t)=\emptyset,\
O\cap\WF^{*,l}_{\scl,L^\infty(G_t)}(u_t)=\emptyset\Mand\\
(O\setminus\Phi_+(\{q\}))\cap\WF^{*,r}_{\scl,L^\infty(G_t)}(u_t)=\emptyset\\
\text{ implies } O\cap\WF^{*,r}_{\scl,L^\infty(G_t)}(u_t)=\emptyset.
\end{multline*}
\end{prop}

\begin{proof} Let $q=(y_0,\nu_0,0),$ $\nu_0<0,$ in our local
coordinates.  Due to our assumptions, and Theorem~\ref{thm:prop-sing}, 
we only need to make a commutator positive on $\Phi_+(\{q\})$.

For $l>-1$ consider a symbol $a$ of the form
\begin{equation}
a=x^{-l-1}\chi(\nu)\phi(y,\mu)\psit(x)\psi_0^2(p)\geq 0
\label{eq:fine-reg-t-ml-8}\end{equation}
where $\psit\in\Cinf_c(\Real)$ is identically $1$ near $0$ and
is supported in a bigger neighbourhood of $0$ (it is simply a cutoff
near $\pa X$), $\psi_0\in\Cinf_c(\Real;[0,1])$ supported in $(\ev-\delta,
\ev+\delta),$ $\phi\geq 0$ is $\Cinf,$ supported near $(y_0,0),$ identically
$1$ in a neighbourhood of this point,
$\chi\in\Cinf_c(\Real;[0,1])$ supported in $[\nu_0-\ep,
\nu_0+\ep],$ identically $1$ near $\nu_0,$ $\chi'\leq 0$ on
$[\nu_0-\ep/2,\nu_0+\ep],$ $\chi'\geq 0$ on $[\nu_0-\ep,\nu_0-\ep/2]$
and $\chi(\nu)=e^{M/(\nu-(\nu_0+\ep))}$ on $(\nu_0+\ep/2,\nu_0+\ep)$.
In addition, we choose $\ep>0$ so that $[\nu_0-2\ep,
\nu_0+2\ep]\times\supp d\phi\subset O\setminus\Phi_+(\{q\})$.
This can be achieved by fixing $\phi$ and making $\ep>0$ small since
$\{\nu_0\}\times\supp d\phi$ is disjoint from the compact set
$\Phi_+(\{q\})\cap \overline{O}$.

Let $\rho\in\Cinf_c(\Real)$ be identically $0$ on $(-\infty,\nu_0-3\ep/8)$,
identically $1$ on $(\nu_0-\ep/4,+\infty).$ Notice that $\chi' \leq 0$ on
the support of $\rho$. Then, near $x=0$,
\begin{equation}\begin{split}
&\scH_p a=2((l+1)\nu\chi(\nu)+|\mu|^2_y\chi'(\nu))\phi(y,\mu)
x^{-l-1}\psi_0^2(p) +x^{-l-1}\chi(\nu)\psi_0^2(p)\scH_p\phi\\
&\quad=-b^2x^{-l-1}+ex^{-l-1},\\
&b^2=\Big(-2\nu(l+1)\chi(\nu) -2\rho^2(\nu)
|\mu|^2_y\chi'(\nu) \Big)\phi(y,\mu)\psi_0^2(p),\\
&e=\Big(2(1-\rho^2(\nu))|\mu|^2_y\chi'(\nu))\phi(y,\mu)
+\chi(\nu)\scH_p\phi \Big)\psi_0^2(p),
\end{split}
\label{eq:fine-reg-t-ml-16}\end{equation}
and $b$ is $\Cinf$ by the construction of $\chi$.
Note that $\supp e$ is disjoint from $\WFsc^{*,r}(u)$, and $a$, $b$, $e$
satisfy \eqref{eq:abe-cond}. 
Thus, by Lemma~\ref{lemma:pos-comm}, we gain half a power of $x$ in our
regularity of $u_t$ near $\Rp$, \ie, we get
$\WF^{*,l+1/2}_{\scl,L^\infty(G_t)}(u_t) = \emptyset$ near $q$. By
Theorem~\ref{thm:prop-sing} this gives $O \cap
\WF^{*,l+1/2}_{\scl,L^\infty(G_t)}(u_t) = \emptyset$. Iterating the
argument, we get $O \cap \WF^{*,r}_{\scl,L^\infty(G_t)}(u_t) = \emptyset$
after a finite number of steps. 
\end{proof}

\section{Microlocal eigenfunctions}\label{sec:micro}

With each zero $\Rp$ of the vector field $W,$ when
$\ev>V(\pi(\Rp)),$ we associate spaces of microlocally incoming, respectively
outgoing  eigenfunctions. These are microfunctions near $\Rp$ with respect
to the scattering wavefront set, satisfy the eigenfunction equation
microlocally near
$\Rp,$ and have wave front set in $\nu<\nu(\Rp)$, resp.\ 
$\nu>\nu(\Rp)$.

For a $W$-balanced neighbourhood $O\ni\Rp,$ as in Definition~\ref{HMV.92},
the spaces of incoming and outgoing microlocal eigenfunctions, $\tilde
E_{\mic,\pm}(\Rp,\ev),$ are defined in \eqref{HMV.42}. The propagation of
singularities in regions of real principal type (\ie, where $W\neq 0$) 
shows that stronger restrictions on $\WFsc(u)$ follow directly.

\begin{lemma}\label{lemma:E-mic-structure-8}
If $O\ni\Rp$ is, for $\Rp\in \RP(\ev),$ a $W$-balanced neighbourhood then 
for every $u\in \tilde E_{\mic,\pm}(O,\ev)$ 
\begin{gather*}
\WFsc(u)\cap O\subset\Phi_\pm(\{\Rp\})\Mand\\
\WFsc(u)\cap O=\emptyset\Longleftrightarrow \Rp\nin\WFsc(u).
\label{HMV.196}\end{gather*}
\end{lemma}

\noindent Thus, we could have defined $\tilde E_{\mic,\pm}(O,\ev)$ by
strengthening the restriction on the wavefront set to $\WFsc(u)\cap
O\subset\Phi_\pm(\{\Rp\}).$ With such a definition there is no need for $O$
to be $W$-balanced; the only relevant bicharacteristics would be those
contained in $\Phi_\pm(\{\Rp\}).$

\begin{proof}
We may assume that $\Rp\in\RP_+(\ev)$.
For the sake of definiteness also assume that $u\in \tilde
E_{\mic,\out}(O,\ev)$; the 
other case follows similarly. Consider $\zeta\in
O\setminus\{\Rp\}.$ If $\nu(\zeta)<\nu(\Rp),$ then $\zeta\nin\WFsc(u)$ by
the definition of $\tilde E_{\mic,\out}(O,\ev),$ so we may suppose that
$\nu(\zeta)\geq\nu(\Rp).$ Since $\Rp\in \Phi_+(\{\Rp\})$ we may also
suppose that $\zeta\neq\Rp.$

Let $\gamma:\Real\to\Sigma(\ev)$ be the bicharacteristic through $\zeta$
with $\gamma(0)=\zeta$.
As $O$ is $W$-convex, and $\WFsc((P-\ev)u)\cap O=\emptyset$,
Theorem~\ref{thm:prop-sing} shows that
\begin{equation*}
\zeta\in\WFsc(u)\Rightarrow \gamma(\Real)\cap O\subset\WFsc(u).
\end{equation*}
As $O$ is $W$-balanced,
there exists $\zeta'\in \overline{\gamma(\Real)}\cap O$ such that
$\nu(\zeta')=\nu(\Rp)$.

If $\zeta'=\gamma(t_0)$ for some $t_0\in\Real$,
then for $t<t_0,$ $\nu(\gamma(t))<\nu(\gamma(t_0))=\nu(\Rp),$ and for
sufficiently small $|t-t_0|,$ $\gamma(t)\in O$ as $O$ is open. Thus,
$\gamma(t)\nin\WFsc(u)$ by the definition of $\tilde E_{\mic,\out}(O,\ev)$,
and hence we deduce that $\zeta\nin\WFsc(u)$.

On the other hand, if $\zeta'=\lim_{t\to-\infty}\gamma(t),$ then
$\zeta'\in \RP(\ev),$ and $\zeta'\in O,$ so the fact that
$O$ is $W$-balanced shows that $\zeta'=\Rp,$ hence $\zeta\in\Phi_+(\{\Rp\})$.

The last statement of the lemma follows from Theorem~\ref{thm:prop-sing},
and the fact that $\WFsc(u)$ is closed.
\end{proof}

\begin{cor} \label{HMV.197}
Let $O$ be a $W$-balanced neighbourhood of $\Rp\in\RP_+(\ev).$ If
$\pi(\Rp)$ is a minimum of $V|_Y$ and $u\in \tilde E_{\mic,\out}(O,\ev)$ then 
$\WFsc(u)\cap O\subset\{\Rp\}$.
If $\pi(\Rp)$ is a maximum of $V|_Y$ and $u\in \tilde E_{\mic,\pm}(O,\ev)$
then $\WFsc(u)\cap O\subset L_\pm$. Here $L_+$ is that
one of the Legendrians $L_1$, $L_2$ from 
Proposition~\ref{prop:smooth-Leg-saddle} which satsifies
$L_+ \subset \{ \nu \geq
\nu(\Rp) \},$ and $L_-$ is the other one of these.
\end{cor}

To avoid truly microlocal arguments in the analysis of these spaces, such
as microlocal solvability along the lines of H\"ormander \cite{MR48:9458},
we introduce, as a technical device, 
an operator $\tilde P$ which arises from $P$ by altering $V$
appropriately, and use global analysis of $\tilde P.$

\begin{lemma}\label{lemma:Vt} Fix $\ev > \min \Vy$, and $\tilde\nu>0$, and set
$K=\Vy^{-1}((-\infty,\ev-\tilde\nu^2]).$ There exists $\tilde V\in\Cinf(X)$
with $\tilde \Vy$ Morse such that the minima of $\Vy$ and $\tilde \Vy$
coincide, $V$ and $\tilde V$ themselves coincide in a neighbourhood of
their minima, $\tilde \Vy\geq \Vy,$ $\tilde V|_K=V|_K$ and no maximum value
of $\tilde V$ lies in the interval $(\ev-\tilde\nu^2,\ev).$

Let $\tilde P=\Lap+\tilde V$ be the corresponding Hamiltonian,
$\tilde\Sigma(\ev)$ its characteristic variety at eigenvalue $\ev$ and
$\widetilde{\RPout}(\ev)$ the set of its radial points in $\nu>0$.
If $\Rp\in\widetilde{\RPout}(\ev)$ with $\nu(\Rp)<\tilde\nu,$ 
then $\pi(\Rp)$ is a minimum of $\tilde V$ and the full symbols of $P$ and 
$\tilde P$ coincide in a neighbourhood of $\Rp$; in particular,
\begin{equation}
\Sigma(\ev)\cap\{\nu\geq\tilde\nu\}
=\tilde\Sigma(\ev)\cap\{\nu\geq\tilde\nu\}.
\label{eq:Sigma=Sigmat}\end{equation}
\end{lemma}

\begin{proof} It suffices to modify $V$ in small neighbourhood $U_{z}$ of
each maximum $z$ of $\Vy$ with value in the range
$(\ev-\tilde\nu^2,\ev]$ (if there are any). We do this so that the new
potential has $\tilde V>\ev-\tilde\nu^2$ on these sets, has only one
critical point in each and at that maximum takes a value greater than
$\ev.$

Now $\Rp\in \widetilde{\RPout}(\ev)$ implies that $z=\pi(\Rp)$ is
a critical point of $\tilde \Vy$ and $\nu(\Rp)^2+\Vt(z)=\ev.$ Suppose that
$z$ is a maximum of $\tilde V$ and $0<\nu(\Rp)<\tilde\nu.$ Then
$\nu(\Rp)^2+\Vt(z)=\ev$ shows that $\ev-\tilde\nu^2<\Vt(z)<\ev$,
which is contradicted by the construction of $\Vt.$ Thus, $z$ is
a minimum of $\Vt.$

Next, on $\Sigma(\ev)\cap\{\nu\geq\tilde\nu\},$ $\nu^2+|\mu|_y^2+\Vy=\ev$,
hence $\Vy\leq \ev-\tilde\nu^2,$ so $\Vy=\tilde \Vy$, and therefore
$\Sigma(\ev)\cap\{\nu\geq\tilde\nu\}\subset\tilde\Sigma(\ev).$  With the
converse direction proved similarly, \eqref{eq:Sigma=Sigmat} follows.
\end{proof}

\begin{rem}\label{rem:P-tilde} The point of this Lemma is that it allows
one to assume, in any argument concerning $q \in \RP_+(\ev)$, that there is
no $q' \in \Max_+(\ev)$ with $\nu(q') < \nu(q)$. The virtue of this is
illustrated in the proof of the following continuation result.
\end{rem}

\begin{lemma}\label{HMV.136} Suppose $u\in\CmI(X)$
satisfies 
\begin{equation*}
\WFsc(u)\subset\{\nu\geq\nu_1\}\Mand\WFsc((P-\ev)u)\subset\{\nu\geq\nu_2\},
\end{equation*}
for some $0<\nu_1<\nu_2$, then there exists $\tilde u\in\CmI(X)$ with
$\WFsc(u-\tilde u)\subset\{\nu\geq\nu_2\}$ and $(P-\ev)\tilde u\in\dCI(X).$
\end{lemma}

\begin{proof}
One is tempted to try to solve this by adding $R(\ev+i0)((P-\ev)u)$ to
$u$. This does not quite work, however, since the wavefront set of this is
contained in $\WFsc((P-\ev)u) \cup \Phi_+(\WFsc((P-\ev)u \cap \Sigma(\ev)))
\cup \Phi_+(\RP_+(\ev))$, and only the first two sets are contained in 
$\{ \nu \geq \nu_2 \}$. The third is contained in the set $\Min_+(\ev)$,
together with the flowouts along $L_2$ from $\Max_+(\ev)$. The flowouts are
a nuisance here because they `bump into' the rest of the wavefront set at 
$\{ \nu \geq \nu_2 \}$. To avoid dealing with this problem, we use the
operator $\tilde P$ where such points are eliminated. 

Let $\tilde\nu$ be such that $\nu_1<\tilde\nu<\nu_2,$ and $\tilde\nu$
is sufficiently
close to $\nu_1$ so that there are no critical points in $\Sigma(\ev)$ with
$\nu\in [\tilde \nu, \nu_2)$. Choose an operator $A\in\Psisc^0(X)$ with
$\WFscp(\Id-A)\cap\Sigma(\ev)\subset \{\nu\geq\tilde\nu\}$ and
$\WFsc(A)\subset\{\nu \leq \nu_2\}.$
Then $f=(P-\ev)Au$ has wavefront set confined to $\Sigma
(\ev)\cap\{\tilde \nu < \nu < \nu_2 \}.$ 

Let $\Vt$ be as in Lemma~\ref{lemma:Vt}, and let $\tilde\Pi$ be the
orthogonal projection off the $L^2$-nullspace of $\Pt-\ev$; $\Id-\tilde\Pi$
is a finite
rank projection onto a subspace of $\dCI(X)$. By \eqref{eq:Sigma=Sigmat},
$\WFsc(f)\subset\tilde\Sigma(\ev)\cap\{\tilde\nu<\nu<\nu_2\}.$ Now consider
$v=\tilde R(\ev+i0)\tilde\Pi f$ where
$\tilde R(\ev)$ is the resolvent for $\tilde
P.$ This has wavefront set in $\tilde\Sigma(\ev)\cap\{\nu>0\}.$
Furthermore in $\nu<\tilde\nu,$ where $f$ has no wavefront set, its only
wavefront set is associated to critical points. These are all minima, by
construction, so in $\nu<\tilde\nu,$ $\WFsc(v)$ only has these as
isolated points. Using microlocal cut-offs near these points we may excise
the wavefront set and so obtain a solution of $(\tilde P-\ev)v'=f+f',$
where $f'\in\dCI(X)$ and $\WFsc(v')$ is confined to $\nu\geq\tilde\nu$ so
\begin{equation*}
\WFsc(v')\subset \tilde\Sigma(\ev)\cap\{\nu\geq\tilde\nu\}
=\Sigma(\ev)\cap\{\nu\geq\tilde\nu\}.
\end{equation*}
Since
\begin{equation*}
\WFsc((\tilde P-P)v')\subset \Sigma(\ev)\cap\{\nu\geq\tilde\nu\}\cap
\WFscp(\tilde P-P)=\emptyset,
\end{equation*}
 it follows that $(P-\ev)v'=f+f''$
with $f''\in\dCI(X).$
Finally then set $\tilde u=Au-v'.$ By construction
this reduces to $u$ microlocally in $\nu<\tilde\nu$ and satisfies
$(P-\ev)\tilde u\in\dCI(X).$ Since there are no critical points
with $\nu\in[\tilde\nu, \nu_2)$ it follows from Theorem~\ref{thm:prop-sing} 
that $\WFsc(u-\tilde u)$ is contained in $\{\nu\geq\nu_2\}.$
\end{proof}

If $O_1$ and $O_2$ are two $W$-balanced neighbourhoods of $\Rp$ then
\begin{equation}
O_1\subset O_2\Longrightarrow
\tilde E_{\mic,\pm}(O_2,\ev)\subset\tilde E_{\mic,\pm}(O_1,\ev).
\label{HMV.84}\end{equation}
Since $\{u\in\CmI(X);\WFsc(u)\cap O=\emptyset\}\subset \tilde
E_{\mic,\pm}(O,\ev)$ for all $O$ and this linear space decreases with $O,$
the inclusions \eqref{HMV.84} induce similar maps on the quotients  
\begin{equation}
\begin{gathered}
E_{\mic,\pm}(\Rp,\ev)=\tilde E_{\mic,\pm}(O,\ev)/
\{u\in\CmI(X);\WFsc(u)\cap O=\emptyset\},\\
O_1\subset O_2\Longrightarrow
E_{\mic,\pm}(O_2,\ev)\longrightarrow E_{\mic,\pm}(O_1,\ev).
\end{gathered}
\label{HMV.85}\end{equation}

\begin{lemma}\label{HMV.86} Provided $O_i,$ for $i=1,$ $2,$ are $W$-balanced
neighbourhoods of $\Rp$ and the closure of $O_2$ contains no radial points
other than $\Rp,$ the map in \eqref{HMV.85} is an isomorphism.
\end{lemma}

\begin{proof}
We may assume $\Rp\in\RP_+(\ev)$. We prove the lemma first
for $E_{\mic,+}(O_i,\ev)$; we need a different method for
$E_{\mic,-}(O_i,\ev)$.

If $\Rp \in \Min_\pm(\ev)$, the lemma follows directly
from Corollary~\ref{HMV.197}, so we may assume that $\Rp \in \Max_\pm(\ev)$.
Let $L_+$ be that
one of the Legendrians $L_1$, $L_2$ from 
Proposition~\ref{prop:smooth-Leg-saddle} with $L_+ \subset \{ \nu \geq
\nu(\Rp) \}$.  

The map in \eqref{HMV.85} is injective since any element $u$
of its kernel has a representative $\tilde u\in\tilde E_{\mic,+}(O_2,\ev)$
which satisfies $\Rp\nin\WFsc(\tilde u),$ hence $\WFsc(\tilde u)\cap O_2=\emptyset$ by
the second half of Lemma~\ref{lemma:E-mic-structure-8}, so $u=0$ in
$E_{\mic,+}(O_2,\ev).$

The surjectivity follows from Lemma~\ref{HMV.136}. 
The assumption that $\overline{O_2}$ contains no radial points
other than $\Rp$ allows us to find an open set $O\supset\overline{O_2}$
which is $W$-balanced.
Now let $A\in\Psisc^{-\infty,0}(X)$ be microlocally the identity
on $L_+\cap \overline{O_1}$ and supported in a small neighbourhood of $L_+$
inside $O$. Then there exists $\nu_2 > \nu(\Rp)$ such that $\nu > \nu_2$ on 
$L_+\setminus O_1$, and $\WFscp(A)\setminus O_1 \subset \{ \nu \geq \nu_2 \}$.
Since $\WFsc(u)\cap O_1\subset L_+,$ $\WFsc(Au-u)\cap O_1=\emptyset$.
In addition, $\WFsc(Au)\subset\WFscp(A)\cap\WFsc(u),$ hence
$\nu\geq\nu(\Rp)$ on $\WFsc(Au).$ Moreover, $\WFsc(Au-u)\cap O_1=\emptyset$
implies that
\begin{equation*}
\WFsc((P-\ev)Au)\cap O_1=
\WFsc((P-\ev)Au-(P-\ev)u)\cap O_1=\emptyset,
\end{equation*}
so $\WFsc((P-\ev)Au)\subset\WFscp(A)\setminus O_1,$ hence is contained in 
$\{ \nu \geq \nu_2 \}$. 
Hence by Lemma~\ref{HMV.136}, there exists $\tilde u\in\dist(X)$
such that $\nu\geq\nu_2$ on $\WFsc(\tilde u-Au)$ and
$(P-\ev)\tilde u\in\dCI(X)$.
In particular, $\nu\geq\nu(\Rp)$ in $\WFsc(\tilde u),$ so $\tilde u
\in \tilde E_{\mic,\out}(O_2,\ev).$ Moreover,
$\nu\geq\nu_2$ on $\WFsc(\tilde u-u)\cap O_1,$ hence
by the second half of Lemma~\ref{lemma:E-mic-structure-8} $\WFsc(\tilde u-u)
\cap O_1=\emptyset,$ so $\tilde u$ and $u$ have the same image in
$E_{\mic,\out}(O_1,\ev).$ This shows surjectivity.

Rather applying Lemma~\ref{HMV.136}, we could have used H\"ormander's
existence theorem in the real principal type region \cite{MR48:9458},
to find a microlocal
solution $v$ of $(P-\ev)v=(P-\ev)Au$ on $O$, with
$\WFsc(v)\subset \Phi_+(\WFsc((P-\ev)Au))$, and then taken $\tilde u=Au-v$.
This method also allows us to deal with $E_{\mic,-}(O_i,\ev).$
\end{proof}

It follows from this Lemma that the quotient space $E_{\mic,\pm}(\Rp,\ev)$
in \eqref{HMV.85} is well-defined, as the notation already indicates, and
each element is determined by the behaviour microlocally `at' $\Rp.$

\begin{cor}\label{cor:min-ident-8}
If $\Rp \in \Min_+(\ev)$ then 
\begin{equation}
E_{\mic,\out}(\Rp,\ev)\simeq
\left\{u\in\CmI(X);(P-\ev)u\in\dCI(X),\ \WFsc(u)\subset\{\Rp\}\right\}/
\dCI(X).
\label{HMV.89}\end{equation}

If $\Rp \in \Max_+(\ev)$ then every element of $E_{\mic,\out}(\Rp,\ev)$
has a representative $u$ such that
$(P-\ev)u\in\dCI(X),$ and $\WFsc(u)\subset \Phi_+(\{\Rp \})$.
\end{cor}

\begin{proof} Suppose first that $\Rp \in \Min_+(\ev)$ then any microlocal
eigenfunction, $u\in\tilde E_{\mic,\out}(O,\ev),$ has $\Rp$ as (at most)
the only point in its wavefront set within $O.$ Using a microlocal cutoff
around $\Rp$ with support in $O$ gives a map which realizes the isomorphism
in \eqref{HMV.89}.

If $\Rp \in \Max_+(\ev)$ then $\tilde u$ constructed in the proof
of Lemma~\ref{HMV.86} provides a representative such that
$(P-\ev)u\in\dCI(X)$, with $\WFsc(\tilde u)$ contained in the union of 
$\Phi_+(\{ \Rp \})$ and $\Phi_+(\{\Rp'\})$, for $\Rp' \in \RP_+(\ev)$ with
$\nu(\Rp') > \nu(\Rp)$. If we choose $\Rp'$ from the set
\begin{equation}
\{ \Rp' \in \RP_+(\ev) \cap \WFsc(\tilde u) \mid \nu(\Rp') > \nu(\Rp), \
\Rp' \notin \Phi_+(\{ \Rp \}) \},
\label{eq:removing}\end{equation}
with $\nu(\Rp')$ minimal, then by localizing $\tilde u$ near $\Rp'$ we have
an element $v$ of $E_{\mic,+}(\Rp')$. By subtracting from $\tilde u$ a
representative for $v$ given by Lemma~\ref{HMV.86}, we remove the wavefront
set near $\Rp'$. Inductively choosing radial points from \eqref{eq:removing}
and performing this procedure, we remove all wavefront set from $\tilde u$
except that contained in $\Phi_+(\{ \Rp \})$. 
\end{proof}

At radial points corresponding to maxima these results can be strengthened
by using a positive commutator estimate which is dual, in some sense, to
Proposition~\ref{prop:micro-fine-prop-t}. 

\begin{prop}\label{prop:micro-rough-prop}
Assume that $\Rp\in \RP_+(\ev)$, and let 
$O$ be a $W$-balanced neighbourhood of
$\Rp$. If, for some $r<-1/2,$ $O\cap\WFsc^{*,r+1}((P-\ev)u)=\emptyset$ and
$O \cap \WFsc^{*,r}(u)\subset \Phi_+(\{\Rp\})$, then $O \cap \WFsc^{*,r}(u) =
\emptyset$.  The same
holds true if $\Rp\in \RP_-(\ev)$, $u \in E_{\mic,-}(\Rp,\ev)$, and 
$O\cap \WFsc^{*,r}(u)\subset \Phi_-(\{\Rp\})$. 
\end{prop}

\begin{proof} Assume that $\Rp\in\Max_+(\ev)$ for the sake of
definiteness, and let $L_+$ be that one of the Legendrians $L_1$, $L_2$ from 
Proposition~\ref{prop:smooth-Leg-saddle} with $L_+ \subset \{ \nu \geq
\nu(\Rp) \}$. The argument for $\Rp \in \Min_+(\ev)$ is similar, and slightly
easier since in that case $\Phi_+(\{\Rp\})$ consists of the single point
$\Rp$. The case $\Rp \in \RP_-(\ev)$ follows by taking the complex conjugate.

Let $l<-1$, $l<r$, be such that $u \in \Hsc^{*,l}(X)$, 
let $\nu_0 = \nu(\Rp) > 0$ and let $y_0=y(\Rp)$ in some local coordinates.
Let $a$, $\chi$, $\phi$, $\tilde \psi$ and $\psi_0$
be as in Proposition~\ref{prop:micro-fine-prop-t} (except that here $\nu_0
> 0$). Then we may write 
\begin{equation}
\scH_p a=-b^2x^{-l-1}+ex^{-l-1}
\label{eq:fine-reg-t-ml-16-2}
\end{equation}
as in Proposition~\ref{prop:micro-fine-prop-t}, since now $l+1$ has changed
signs, but so has $\nu$ (from negative to positive) on the support of
$\chi$. If we use Lemma~\ref{lemma:pos-comm-reg} instead of
Lemma~\ref{lemma:pos-comm} (since we are dealing
here with a single function rather than a family), the argument from
Proposition~\ref{prop:micro-fine-prop-t} 
applies, except we can only reach values of $r$ less than $-1/2$ since we
must have $l<-1$. 
\end{proof}

\begin{cor}\label{cor:micro-rough-prop}
If $\Rp\in \RP_+(\ev),$ $O$ is a $W$-balanced neighbourhood of $\Rp$
and $u\in \tilde E_{\mic,\out}(O,\ev),$ then $\WFsc^{*,r}(u)\cap O=\emptyset$
for all $r<-1/2$.
\end{cor}

\section{Test modules and iterative regularity}\label{S.Test}

To investigate the regularity of microlocal eigenfunctions we use `test
modules' of scattering pseudodifferential operators. Since we work
microlocally we consider scattering pseudodifferential operators
microlocally supported in an open set $O\subset C=\pa (\scbT^*X).$ In
particular we shall use the notation 
\begin{equation}
\Psisc^{m,l}(O)=\left\{A\in \Psisc^{m,l}(X); \WFscp(A)\subset O\right\}.
\label{HMV.93}\end{equation}
Since we are only interested in boundary regularity here we shall suppose
that $O\subset C_{\pa}=\scT^*_{\pa X}X$ (see \eqref{HMV.16}) and so we may
take $m=-\infty.$ 

\begin{Def} A \emph{test module} in an open set $O\subset C_{\pa}$ is a
linear subspace $\cM\subset\Psisc^{-\infty,-1}(O)$ which contains and is a
module over $\Psisc^{-\infty,0}(O),$ which is closed under commutators and which
is finitely generated in the sense that there exist finitely many
$A_i\in\Psisc^{*,-1}(X),$ $i=0,1,\dots,N,$ $A_0=\Id,$ such that each $A\in
\cM$ can be written  
\begin{equation}
A=\sum\limits_{j=0}^N Q_iA_i,\ Q_i\in\Psisc^{-\infty,0}(O).
\label{HMV.94}\end{equation}
\end{Def}

Since we have assumed that $\cM \subset \Psisc^{-\infty,-1}(X),$ $A\in\cM,$
$C\in \Psisc^{-\infty,0}(X)$ implies that
$[A,C]\in\Psi^{-\infty,0}(O)\subset\cM.$ The structure of the module is
thus determined by the commutators of the generators. It certainly suffices
to have
\begin{multline}
[A_i,A_j]=\sum\limits_{k=0}^NC_{ijk}A_k+E'_{jk},\
C_{ijk}\in\Psisc^{0,0}(X),\\
E'_{jk}\in\Psisc^{0,-1}(X),\ \WFscp(E'_{jk})\cap
O=\emptyset.
\label{HMV.95}\end{multline}
This in turn follows from the purely symbolic condition on principal
symbols: 
\begin{equation}
\{a_i,a_j\}=\sum\limits_{k=0}^Nc_{ijk}a_k+e_{ij}
\text{ in a neighbourhood of }\overline{O}
\label{HMV.99}\end{equation}
where the $e_{ij}$ are symbols of order $0.$ 

In our arguments below, we frequently shrink $O,$ so the condition
\eqref{HMV.95} is in essence necessary. When $O'\subset O$ is another open
set we consider the restricted module over $O'$
\begin{equation}
\cM(O')=\sum\limits_{i=0}^N\Psi^{-\infty,0}(O')\cdot A_i.
\label{HMV.97}\end{equation}
which is easily seen to be independent of the choice of generators.

\begin{Def}\label{Test.module}
Let $\cM$ be a test module. For $u\in\CmI(X)$ we say $u\in\Isc^{(s)}(O,\cM)$ if
$\cM^mu\subset\Hsc^{\infty,s}(X)$ for all $m.$ That is, 
for all $m$ and for all $B_i\in\cM,$ $i=1,\ldots,m$, we have
$\prod_{i=1}^m B_i u\in\Hsc^{\infty,s}(X).$ We may define the finite
regularity spaces by saying $u\in\Isc^{(s),M}(O,\cM)$ if
$\cM^mu\subset\Hsc^{\infty,s}(X)$ for all $m\leq M.$
\end{Def}

Note that $\Isc^{(s)}(O,\cM)$ is a well-defined space of microfunctions
over $O$ since if $u\in\CmI(X)$ and $\WFsc(u)\cap O=\emptyset$ then
$Au\in\dCI(X)$ for all $A\in\Psisc^{*,*}(X)$ such that $\WFscp(A)\subset
O;$ in particular this is the case for all $A\in\cM.$ Clearly 
\begin{equation}
\Isc^{(s)}(O,\cM)=\bigcap_{M=0}^\infty\Isc^{(s),M}(O,\cM).
\label{HMV.96}\end{equation}
From the definition of the restriction in \eqref{HMV.97} it also follows
that if $O'\subset O$ then
\begin{equation}
\Isc^{(s),k}(O,\cM)\subset \Isc^{(s),k}(O',\cM).
\label{HMV.98}\end{equation}

Since $\Id\in\cM,$
\begin{equation*}
\Psisc^{-\infty,0}(O)=\cM^0\subset\cM\subset\cM^2\subset\ldots.
\end{equation*}

Definition~\ref{Test.module} above then reduces to 
\begin{equation}
u\in\Isc^{(s),M}(O,\cM)\Longleftrightarrow\cM^{M}(O)u\subset\Hsc^{\infty,s}(X).
\label{HMV.109}\end{equation}
If two test modules $\cM_1$ and $\cM_2$ define equivalent filtrations
of the same set of operators they define the same space of
iteratively-regular distributions. That is, 
\begin{multline}
\cM_1(O)\subset(\cM_2(O))^{k_1}\Mand \cM_2(O)\subset(\cM_1(O))^{k_2}
\Longrightarrow\\ \Isc^{(s)}(O,\cM_1)=\Isc^{(s)}(O,\cM_2),\ \forall\ s\in\bbR.
\label{HMV.110}\end{multline}

\begin{lemma}\label{HMV.101} If $A_i,$ $0\le i\le N,$ are generators for
$\cM$ in the sense of Definition~\ref{Test.module} with $A_0=\Id,$ then 
\begin{equation}
\cM^{k}=\left\{\sum\limits_{|\alpha|\le k} Q_\alpha
\prod\limits_{i=1}^{N}A_i^{\alpha _i},\ Q_\alpha \in\Psi^{-\infty,0}(O)\right\}
\label{HMV.102}\end{equation}
where $\alpha$ runs over multiindices $\alpha
:\{1,\dots,N\}\longrightarrow\bbN_0$ and $|\alpha |=\alpha _1+\dots+\alpha _N.$ 
\end{lemma}

\begin{proof} This is just a consequence of the fact that $\cM$ is a Lie
algebra and module. Thus, we proceed by induction over $k$ since
\eqref{HMV.102} certainly holds for $k=1,$ and by definition for $k=0.$ It
therefore suffice to prove the equality \eqref{HMV.102} modulo
$\cM^{k-1}.$ Then \eqref{HMV.102} just corresponds to
$[\cM^{r},\cM^{s}]\subset\cM^{r+s-1}$, applied with $r+s\leq k,$
which allows the factors to be freely rearranged.
\end{proof}

In the modules we consider, one of the generating elements is
$A_N=x^{-1}(P-\ev)$ and we deal with approximate eigenfunctions, so that
typically $(P-\ev)u\in\dCinf(X).$ We may reorder the basis so that this
element comes last. In this case it is enough to consider the action on $u$
of the remaining generators in an inductive proof of regularity since
if an operator as on the right hand side of
\eqref{HMV.102} with $\alpha_N\neq 0$ is applied
to $u$ then we obtain an element of $\dCI(X).$
We therefore use the `reduced' multiindex notation, with an additional
index to denote the weight, and set 
\begin{equation}
A_{\alpha,s}= x^{-s}\prod_{i=1}^{N-1}A_i^{\alpha_i},\ \alpha_i\in\bbN_0,\ 
1 \leq i \leq N-1.
\label{eq:red-multi}\end{equation}

\begin{cor}\label{cor.reduced.basis} Suppose $\cM$ is a test module,
$A_N=x^{-1}(P-\ev)$ is a generator, $u\in\CmI$ satisfies 
$(P-\ev)u\in\dCinf(X)$ and $u\in\Isc^{(s),m-1}(O,\cM).$ Then for $O'\subset
O,$ $u\in\Isc^{(s),m}(O',\cM)$ if for each multiindex $\alpha,$ with
$|\alpha|=m,$ $\alpha_N= 0$,
there exists $Q_\alpha\in\Psi^{-\infty,0}(X),$ elliptic on
$O'$ such that $Q_\alpha A_{\alpha,s}u\in L^2_\scl(X).$
\end{cor}

In the case of interest here, $O$ will be a $W$-balanced neighbourhood of
$\Rp \in \RP(\ev)$ and the generators $A_i$ will be certain operators,
related to the geometry of the classical dynamics near $\Rp$, with
principal symbols vanishing at $\Rp$. The last operator, as mentioned
above, will be $A_N=x^{-1}(P-\ev)$. 
Our main tool in the demonstration of regularity with respect to such a
test module is the use of positive commutator techniques. These depend on
finer structure of the module.
In particular the basic source of positivity
is the commutation relation 
\begin{equation}
\begin{gathered}
i[x^{-s-\frac12},P-\ev]=x^{-s+\frac12}C_0, \\
\text{ where } C_0\in\Psi^{1,0}_{\scl}(X) \text{ has principal symbol }
-(2s + 1) \nu,
\end{gathered}
\label{HMV.103}\end{equation}
so has a definite sign in a neighbourhood of a radial point $\Rp \in
\RP_\pm(\ev)$ if $s \neq -1/2$ and $\ev \notin \Cv(V)$. This strict sign of
$\sigma_0(C_0)$ allows us to deal with error terms. To fix signs we shall
assume that $\Rp \in \RP_+(\ev)$ and $\ev \notin \Cv(V)$ in the remainder
of this section. 

So suppose that $f(\nu,y,\mu)$ vanishes at $\Rp\in\RP_+(\ev)$,
and $df$ is
an eigenvector for the linearization of $W$ at $\Rp$ with eigenvalue
$-2\nu r$, where $\nu = \nu(\Rp) > 0$. Then
\begin{equation}
\scH_p(x^{m} f)=x^m[-2\nu(m+r)f + f_2  +{\mathcal{O}}(x)],
\label{HMV.46}\end{equation}
where $f_2$ is smooth and vanishes quadratically at the critical point
$\Rp.$ Thus, 
\begin{equation*}
\scH_p(x^{2m} f^2)=2x^{2m}[-2\nu(m+r)f^2+ f_2 f +{\mathcal{O}}(x)].
\end{equation*}
If $m+r \neq 0$ and $f_2$ is of the form $f \tilde f_2,$ then this
has a fixed sign near $\Rp.$ In general, $f_2$ does not factor this way,
so we need to consider a sum of such commutators corresponding
to functions $f,$ $f'$ such that $df$ and $df'$ span $T^*_{\Rp}\Sigma(\ev);$
then $x^{-2m}\scH_p(x^{2m} f^2)+x^{-2m}\scH_p(x^{2m} (f')^2)$ has a fixed
sign near $\Rp$ on $\Sigma(\ev).$
We use this as a basis for the construction of positive commutators, though
it is more convenient to describe our conditions on the generators $A_i,$
and then spell out in detail in the following sections why these
conditions hold.

The basic condition we require is that the remaining generators of the
test module, $A_i,$ $i=1,\dots,l=N-1,$ satisfy
\begin{equation}
\begin{gathered}
x^{-1}i[ A_i,P-\ev]=\sum\limits_{j=0}^N C_{ij} A_j,\
C_{ij}\in\Psi^{*,0}_{\scl}(X),\\
\Mwhere \sigma_{\pa}(C_{ij})(\Rp)=0,\ j\neq i,\\
\Mand \re\sigma_{\pa}(C_{jj})(\Rp)\geq 0,\ 0<j\le l,\ \Rp \in \RP_+(\ev),
\ s<-\frac12.
\end{gathered}
\label{eq:A_i-comm}\end{equation}

However, we will occasionally need slightly weaker conditions than
\eqref{eq:A_i-comm}. These weaker conditions allow 
$\sigma_{\pa}(C_{ij})(\Rp)\neq 0$ even for $i\neq j$, but they require
that the matrix $\sigma_{\pa}(C_{ij})(\Rp)$ is lower triangular. Hence,
in this case, the ordering of the $A_i$ (by their index) is not arbitrary.
Our weaker condition is then
\begin{equation}
\begin{gathered}
x^{-1}i[ A_i,P-\ev]=\sum\limits_{j=0}^N C_{ij} A_j,\
C_{ij}\in\Psi^{*,0}_{\scl}(X),\\
\Mwhere \sigma_{\pa}(C_{ij})(\Rp)=0,\ j>i,\\
\Mand \re\sigma_{\pa}(C_{jj})(\Rp)\geq 0,\ 0<j\le l,\ \Rp \in \RP_+(\ev),
\ s<-\frac12..
\end{gathered}
\label{eq:A_i-comm-p}\end{equation}
We need the following notation.
Suppose that $\alpha,$ $\beta$ are multiindices. We say that $\alpha\sim
\beta$ if
there exist $i\neq j$ such that for $k\nin\{i,j\},$ $\alpha_k=\beta_k$,
while $\alpha_i=\beta_i+1$, $\alpha_j=\beta_j-1.$ This means that
$A_\alpha$ has an additional factor of $A_i,$ and one fewer factor
of $A_j$ than $A_\beta.$ If $\alpha\sim\beta,$ we write $\alpha<\beta$
provided that $i<j.$

The following is a technical Lemma we need for Proposition~\ref{HMV.108} below,
which gives a method for showing membership of $u \in \Isc^{(s),m}(O,\cM)$.

\begin{lemma}\label{HMV.104} Let $\Rp \in \RP_+(\ev)$ and $s<-\frac12,$
where $\ev \notin
\Cv(V)$. Suppose $Q \in \Psisc^{-\infty,0}(O)$, let $K_\alpha > 0$ be arbitrary
positive constants, and let $C_0$ and $C_{ij}$ be given by
\eqref{HMV.103}, resp. \eqref{eq:A_i-comm}. Then, 
assuming \eqref{eq:A_i-comm}, and using the notation of
\eqref{eq:red-multi}, we have 
\begin{multline}
\sum_{|\alpha|=m}i K_\alpha[A_{\alpha,s+1/2}^*Q^*QA_{\alpha,s+1/2},P-\ev]\\
=\sum_{|\alpha|,|\beta|=m}A_{\alpha,s}^*Q^*C'_{\alpha\beta}QA_{\beta,s}\\
+\sum_{|\alpha|=m}
\left(A_{\alpha,s}^*Q^*E_{\alpha,s}+E^*_{\alpha,s}Q A_{\alpha,s}\right)
+\sum_{|\alpha|=m} A_{\alpha,s}^* i[Q^*Q,P-\ev]A_{\alpha,s}, 
\label{HMV.105}\end{multline}
where
\begin{gather}
E_{\alpha,s}=x^{-s}E_{\alpha},\ E_{\alpha}\in \cM^{m-1}+\cM^{m-1}A_N,
\ \WFscp(E_{\alpha})\subset\WFscp(Q),
\label{HMV.106}\\
\begin{gathered}
C'_{\alpha\beta} \in \Psisc^{1,0}(X), \ 
\sigma_{\pa}(C'_{\alpha\alpha})(\Rp)=2K_\alpha\left(-(2s+1)\nu
+\sum_j\alpha_j\Re\sigma_{\pa}(C_{jj})(\Rp)\right)>0 \\
\Mand \sigma_{\pa}(C'_{\alpha\beta})(\Rp)=0\Mfor\alpha\neq\beta.
\end{gathered}
\label{eq:C'>0}\end{gather}
In addition, $C'_{\alpha\beta}=0$ unless either
$\alpha=\beta$ or $\alpha\sim\beta.$

Assuming instead \eqref{eq:A_i-comm-p} instead of \eqref{eq:A_i-comm}, the
conclusions hold with the last equation of \eqref{eq:C'>0} replaced by
\begin{equation}\begin{split}
&\sigma_{\pa}(C'_{\alpha\beta})(\Rp)=K_\beta\beta_j
\overline{\sigma_\pa(C_{ji})(\Rp)}
\Mfor\alpha\sim\beta,\ \alpha<\beta,\\
&\sigma_{\pa}(C'_{\alpha\beta})(\Rp)=K_\alpha\alpha_i \sigma_\pa(C_{ij})(\Rp)
\Mfor\alpha\sim\beta,\ \alpha>\beta.
\label{eq:C'>0p}\end{split}\end{equation}
\end{lemma}

\begin{rem} The first term on the right hand side of \eqref{HMV.105} is the
principal term, in terms of $\cM$-order, since both the $A^*_{\alpha,s}$
and $A_{\beta,s}$ terms have $\cM$-order equal to $m$. Note that this term
has nonnegative principal symbol. In the second term,
the $E$ terms have $\cM$-order $m-1$ which allows these to be treated as
error terms. To deal with the third term, involving the commutator $[Q^*Q,
P-\ev]$, we need additional information about $Q$, discussed below. 
\end{rem}

\begin{proof}
We obtain the identity \eqref{HMV.105} by computing
\begin{multline}
[QA_{\alpha,s+1/2},P-\ev]=
[Qx^{-s-1/2}A_1^{\alpha_1}\ldots A_{N-1}^{\alpha_{N-1}},P-\ev],\\
Q\in\Psisc^{-\infty,0}(O),\ |\alpha|=m.
\end{multline}
The commutator with $P-\ev$ distributes over the product, in each term
replacing a single factor $A_i$ by $\sum xC_{ij}A_j.$ After
rearrangement of the order of the factors, giving an error included in
$E_\alpha$ below, it becomes
\begin{multline}
(C_0+\sum_j \alpha_jC_{jj})Qx^{-s+1/2}
A_\alpha+\sum_{|\beta|=m,\ \beta\neq\alpha}
C_{\alpha\beta}Q x^{-s+1/2}A_\beta\\
+x^{-s+1/2}E_{\alpha}+i[Q,P-\ev]A_{\alpha,s+1/2},
\end{multline}
where $\sigma_{\pa}(C_{\alpha\beta})(\Rp)=0$, $C_{\alpha\beta}=0$
unless $\alpha=\beta$ or $\alpha\sim\beta$, and the $E_\alpha$ are as in
\eqref{HMV.106} -- using in particular the inclusion
$x\cM^{m}\subset \cM^{m-1}.$ In addition, for $\alpha\sim\beta,$
\begin{equation*}
\sigma_{\pa}(C_{\alpha\beta})(\Rp)=\alpha_i\sigma_{\pa}(C_{ij})(\Rp).
\end{equation*}
This implies that
\begin{multline}
i[A_{\alpha,s+1/2}^*Q^*Q A_{\alpha,s+1/2},P-\ev]\\
=A_{\alpha,s}^*Q^*(C_0+C_0^*+\sum_j \alpha_j(C_{jj}+C_{jj}^*))Q
A_{\alpha,s}\\
+\sum_{\beta}A_{\alpha,s}^*Q^*C_{\alpha\beta} QA_{\beta,s}
+\sum_{\beta}A_{\beta,s}^*Q^*C_{\alpha\beta}^*
QA_{\alpha,s}+A_{\alpha,s}^*Q^*E_{\alpha,s}\\
+E^*_{\alpha,s} QA_{\alpha,s}
+A_{\alpha,s}^* i[Q^*Q,P-\ev]A_{\alpha,s};
\end{multline}
here $E_{\alpha,s}=x^{-s}E_{\alpha}$ and all sums over $\beta$ are understood
as $|\beta|=m,$ $\beta\neq\alpha.$
Let
\begin{equation}\begin{split}
&C'_{\alpha\alpha}
=K_\alpha\left(C_0+C_0^*+\sum_j \alpha_j(C_{jj}+C_{jj}^*)\right)\\
&C'_{\alpha\beta}=K_\alpha \alpha_i C_{ij}+K_\beta \beta_j C_{ji}^*\Mif
\alpha\sim\beta.
\end{split}\end{equation}
If we assume \eqref{eq:A_i-comm}, we deduce that \eqref{eq:C'>0} holds,
while if we assume \eqref{eq:A_i-comm-p}, \eqref{eq:C'>0p} follows. This
finishes the proof.
\end{proof}

We now consider $C'=(C'_{\alpha\beta}),$ $|\alpha|=m=|\beta|$
as a matrix of operators, or rather as an
operator on a trivial vector bundle with fiber $\Real^{|M_m|}$ over a
neighbourhood of $\Cp=\pi(\Rp)$ in $X,$ where
$|M_m|$ denotes the number of elements of the set $M_m$ of multiindices
$\alpha$ with $|\alpha|=m.$ Let $c'=\sigma_{\pa}(C')(\Rp),$
$c'_{\alpha\beta}=\sigma_{\pa}(C'_{\alpha\beta})(\Rp).$

Suppose first that
\eqref{eq:A_i-comm} holds. Then for any choice of $K_\alpha>0,$ \eg
$K_\alpha=1$ for all $\alpha,$ $c'=\sigma_{\pa}(C')(\Rp)$ is
positive or negative definite with the sign of $\sigma_{\pa}(C_0)(\Rp)$
(whose sign for $\Rp \in \RP_+(\ev)$ is that of $-s-1/2$). 
The same is therefore true microlocally near $\Rp.$ For the sake of
definiteness, suppose that $\sigma_{\pa}(C_0)(\Rp)>0.$ Then there exist a
neighbourhood $O_m$ of $\Rp,$ depending on $|\alpha|=m,$ and
$B\in\Psisc^{-\infty,0}(X),$ $G\in\Psisc^{-\infty,1}(X),$ with
$\sigma_{\pa}(B)>0$ on $O_m$ such that
\begin{equation}\label{eq:Q-O_m}
Q\in\Psisc^{-\infty,0}(X),\ \WFscp(Q)\subset O_m\Rightarrow
Q^*C'Q=Q^*(B^*B + G)Q.
\end{equation}

Assume now that \eqref{eq:A_i-comm-p} holds. There is a natural partial
order on $M_m$ given as follows: $\alpha<\beta$ if there exist
$\alpha^{(j)}$, $j=0,1,\ldots,n,$ $n\geq 1,$ $\alpha^{(0)}=\alpha,$  
$\alpha^{(n)}=\beta$, $\alpha^{(j)}\sim\alpha^{(j+1)}$ and $\alpha^{(j)}
<\alpha^{(j+1)}$ (which has already been defined if
$\alpha^{(j)}\sim\alpha^{(j+1)}$) for
$j=0,1,\ldots,n-1.$ This is clearly transitive, and if $\alpha<\beta$
then $\sum_j j\alpha_j<\sum_j j\beta_j,$ so $\alpha<\beta$ and $\beta<\alpha$
cannot hold at the same time.
We place a total order $<'$ on $M_m$ which is compatible with $<,$ so
$\alpha<\beta$ implies $\alpha<'\beta.$
For the sake of definiteness, suppose that $\sigma_{\pa}(C_0)(\Rp)>0.$
We then choose $K_\alpha$ inductively
starting at the maximal element of $<',$ making $\sigma_{\pa}(C')(\Rp)$
positive definite on $\oplus_{\alpha\leq' \beta}\Real_\beta.$
We can for example take $K_\alpha=1$ for the maximal element $\alpha.$
Once $K_\beta$ has been chosen for all $\beta$ with $\alpha<'\beta,$
so that $c'$ is positive definite on
$\oplus_{\alpha<' \beta}\Real_\beta,$ we remark that of all entries
$c'_{\beta\gamma}$ with $\alpha\leq \beta,\gamma,$
only $c'_{\alpha\alpha}$ depends on $K_\alpha.$
Expanding $(v,c'v)$ for $v=(v',v'')\in\Real_\alpha\oplus \left(
\oplus_{\alpha<' \beta}\Real_\beta\right),$ and using Cauchy-Schwarz for
the cross terms involving both $v'$ and $v'',$ we deduce that for
any $\ep>0,$
\begin{equation*}
(v,c'v)\geq c'_{\alpha\alpha}\| v' \|^2+(v'',c' v'')-\ep \|c'v''\|^2
-\ep^{-1}\| v' \|^2.
\end{equation*}
Since $c'$ is positive definite on $\oplus_{\alpha<' \beta}\Real_\beta,$
we can choose first $\ep>0$ sufficiently small, then $K_\alpha>0$
sufficiently large (depending on $\ep$) such that $c'$ is positive
definite on $\oplus_{\alpha\leq' \beta}\Real_\beta.$ Thus, for suitable
constants $K_\alpha,$ $c'=\sigma_{\pa}(C')(\Rp)$ is positive definite,
and we again conclude that there exist a
neighbourhood $O_m$ of $\Rp,$ depending on $|\alpha|=m,$ and
$B\in\Psisc^{-\infty,0}(X),$ $G\in\Psisc^{-\infty,1}(X),$ with
$\sigma_{\pa}(B)>0$ on $O_m$ such that
\begin{equation}\label{eq:Q-O_m-p}
Q\in\Psisc^{-\infty,0}(X),\ \WFscp(Q)\subset O_m\Rightarrow
Q^*C'Q=Q^*(B^*B + G)Q.
\end{equation}

Equations \eqref{HMV.103} and \eqref{eq:A_i-comm}
(or \eqref{HMV.103} and \eqref{eq:A_i-comm-p}) are all that is needed in the
two cases of centers and sinks, where we know {\em a priori} that
$\WFsc(u)\subset\{\Rp\}.$ In these cases we take $Q$ such that
$\Rp\nin\WFscp(\Id-Q),$ which implies that 
\begin{equation}
i[Q^*Q,P-\ev]= x^{1/2}\tilde Fx^{1/2}, \text{ where } \tilde
F\in\Psisc^{0,0}(O), \ \Rp \notin \WFscp(\tilde F).
\label{eq:[Q,P]-min-point}\end{equation}
In the case of a saddle point, \ie $\Rp \in \Max_+(\ev)$,
we need a more general setup and we shall arrange that
$\Rp\nin\WFscp(\Id-Q),$
\begin{equation}\begin{split}
&i[Q^*Q,P-\ev]= x^{1/2} (\tilde B^*\tilde B + \tilde G)x^{1/2}
+x^{1/2}\tilde Fx^{1/2}, \text{ where }\\
&\tilde B,\ \tilde F\in\Psisc^{0,0}(O), \ \tilde G \in \Psisc^{0,1}(X),
\ \Rp\nin\WFscp(\tilde F),
\end{split}
\label{eq:[Q,P]}\end{equation}
and in addition, $\tilde F$ satisfies
$\WFscp(\tilde F) \subset \{ \nu < \nu(\Rp) \}$. This condition on $\tilde F$
will ensure that $\WFscp(\tilde F) \cap \WFsc(u) = \emptyset$ for the
application in section~\ref{sect:saddle}. 

\begin{prop} \label{HMV.108p} Suppose that $m>0,$ $s<-1/2,$
$\Rp\in\RP_+(\ev),$ $\ev\nin\Cv(V),$ either
\eqref{eq:A_i-comm} or \eqref{eq:A_i-comm-p} hold, and let
$O_m$ be as in \eqref{eq:Q-O_m} (or \eqref{eq:Q-O_m-p}).
Suppose that $u\in\Isc^{(s),m-1}(O_m,\cM),$
$\WFsc((P-\ev)u)\cap O_m=\emptyset$ and that there exists 
$Q\in\Psisc^{-\infty,0}(O_m)$ elliptic at $\Rp$
that satisfies \eqref{eq:[Q,P]} with
$\WFscp(\tilde F)\cap \WFsc(u)=\emptyset.$ Then $u\in\Isc^{(s),m}(O',\cM)$
where $O'$ is the elliptic set of $Q.$
\end{prop}

\begin{proof}
First consider $u'\in\Isc^{(s),m}(O_m,\cM).$
Let $Au'=(QA_{\alpha,s} u')_{|\alpha|=m}$, regarded as a column vector of
length $|M_m|$. Now consider
\begin{equation}\begin{split}
&\sum_{|\alpha|=m}
K_\alpha
\langle u',i[A_{\alpha,s+1/2}^*Q^*QA_{\alpha,s+1/2},P-\ev]u'\rangle\\
&=\|BAu'\|^2 + \langle Au', GAu' \rangle +\sum_{|\alpha|=m}
\Big(\langle QA_{\alpha,s} u',E_{\alpha,s} u'\rangle
+\langle E_{\alpha,s} u',QA_{\alpha,s} u'\rangle \Big)\\
&+\sum_{|\alpha|=m} \Big(\|\tilde B A_{\alpha,s}u'\|^2
+\langle A_{\alpha,s}u', \tilde F A_{\alpha,s}u'\rangle
+\langle A_{\alpha,s}u', \tilde G A_{\alpha,s}u'\rangle \Big).
\label{HMV.108e}\end{split}\end{equation}
Hence, dropping the term involving $\tilde B$ and
applying the Cauchy-Schwarz inequality
to the terms with $E_{\alpha,s}$, $G$ and $\tilde G$, we have for any $\ep>0$,
\begin{multline}
\|BAu'\|^2\leq \sum_\alpha
\Big|\langle u',i[A_{\alpha,s+1/2}^*Q^*QA_{\alpha,s+1/2},P-\ev]u'\rangle 
\Big|\\
+\ep \Big( \| Au' \|^2 + \sum_\alpha (\|QA_{\alpha,s} u'\|^2 +
\|A_{\alpha,s} u'\|^2 ) \Big) \\
+\ep^{-1} \Big( \| GAu' \|^2 + \sum_\alpha (\|E_{\alpha,s} u'\|^2
+ \|K_\alpha \tilde G A_{\alpha,s} u'\|^2 ) \Big)
+|\langle A_{\alpha,s}u', \tilde F
A_{\alpha,s}u'\rangle|.
\end{multline}
Choosing $\ep>0$ small enough, the second term on the right can
be absorbed in the left hand side (since $B$ is strictly positive), 
and we get 
\begin{equation}\begin{split}
\frac1{2} \|BAu'\|^2\leq &\sum_\alpha
|\langle u',i[A_{\alpha,s+1/2}^*Q^*QA_{\alpha,s+1/2},P-\ev]u'\rangle|\\
&\qquad
+\ep^{-1} \Big( \| GAu' \|^2 + \sum_\alpha (\|E_{\alpha,s} u'\|^2
+ \|K_\alpha \tilde G A_{\alpha,s} u'\|^2 ) \Big)\\
&\qquad +|\langle A_{\alpha,s}u', \tilde F
A_{\alpha,s}u'\rangle|.
\end{split}\label{eq:B'}\end{equation}
We apply this with $u'$ replaced by $u_r=(1+r/x)^{-1}u 
=\frac{x}{x+r}\,u$ $r>0$, where now $u \in \Isc^{(s),m-1}(O_m,\cM)$.
Then letting $r\to 0,$ using the strong convergence of $(1+r/x)^{-1}$ to the
identity as in
the argument of Lemma~\ref{lemma:pos-comm-reg}, and the
assumption that $\WFscp(\tilde F)\cap\WFsc(u) =\emptyset$, shows that
$BAu\in L^2_\scl(X),$ finishing the proof.
\end{proof}

We may arrange that for all $m,$
\begin{equation}\label{eq:Q-ell}
\overline{O_{m+1}}\subset O_{m}\Mand Q=Q_m\ \text{is elliptic on}\ O_{m+1}.
\end{equation}
By induction, the proposition implies that
if $\WFsc^{*,s}(u)\cap O_1=\emptyset$ for some
$s<-1/2,$ and $\WFsc((P-\ev)u)\cap O_1=\emptyset$ then for all $m$,
$u\in\Isc^{(s),m}(O_{m+1}).$ Since we may have
$\cap_{m=1}^\infty O_m=\{\Rp\},$ we need an additional assumption, amounting
to `propagation of regularity', to extend
the result to $\Isc^{(s)}(O,\cM)$ for some neighbourhood $O$ of $\Rp.$
Of course, there is no need for such an assumption if
$\WFsc(u)\cap O\subset\{\Rp\}$ for then $u\in\Isc^{(s),m}(O_{m+1},\cM)$ implies
that $u\in\Isc^{(s),m}(O,\cM).$ We state these two cases separately, but
prove them together.

\begin{prop} \label{HMV.108} Suppose that $s<-1/2,$
$\Rp\in\RP_+(\ev),$ $\ev\nin\Cv(V),$ either
\eqref{eq:A_i-comm} or \eqref{eq:A_i-comm-p} hold, and let
$O$ be a neighbourhood of $\Rp.$
Then for any $u\in\dist(X),$
\begin{equation*}\begin{split}
\WFsc^{*,s}(u)\cap O=\emptyset,\ \WFsc(u)\cap O\subset\{\Rp\}
&\Mand\WFsc((P-\ev)u)\cap O=\emptyset\\
&\Rightarrow
u\in\Isc^{(s)}(O,\cM).
\end{split}\end{equation*}
\end{prop}

In the next proposition we assume `propagation of regularity' to show that
$u\in\Isc^{(s),m}(O,\cM).$

\begin{prop} \label{HMV.108pp} Suppose that $s<-1/2,$
$\Rp\in\RP_+(\ev),$ $\ev\nin\Cv(V),$ either
\eqref{eq:A_i-comm} or \eqref{eq:A_i-comm-p} hold, and let $Q$ and
$O_m$ be as in \eqref{eq:Q-O_m} (or \eqref{eq:Q-O_m-p}), \eqref{eq:[Q,P]}
and \eqref{eq:Q-ell}.
Suppose also that there exists a
neighbourhood $O$ of $\Rp$ such that for all $m,$
\begin{equation}\label{eq:O_m-O}
u\in\Isc^{(s),m}(O_{m+1},\cM)\Mand \WFsc((P-\ev)u)\cap O=\emptyset
\Rightarrow u\in\Isc^{(s),m}(O,\cM).
\end{equation}
Then for any $u\in\dist(X),$ satisfying
\begin{equation*}
\WFsc^{*,s}(u)\cap O=\emptyset
\Mand\WFsc((P-\ev)u)\cap O=\emptyset\Rightarrow
u\in\Isc^{(s)}(O,\cM).
\end{equation*}
\end{prop}

\begin{proof}
By hypothesis, $u\in\Isc^{(s),0}(O,\cM).$ We now show that for $m\geq 1,$
$u\in\Isc^{(s),m-1}(O,\cM)$ implies that $u\in\Isc^{(s),m}(O,\cM).$
But indeed, $u\in\Isc^{(s),m-1}(O,\cM)$ implies that
$u\in\Isc^{(s),m-1}(O_m,\cM),$ hence $u\in\Isc^{(s),m}(O_{m+1},\cM)$ by
Proposition~\ref{HMV.108p} and \eqref{eq:Q-ell}. In Proposition~\ref{HMV.108},
$\WFsc(u)\cap O\subset\{\Rp\},$ so $u\in\Isc^{(s),m}(O_{m+1},\cM)$
implies that $u\in\Isc^{(s),m}(O,\cM)$. In Proposition~\ref{HMV.108pp},
by \eqref{eq:O_m-O}, $u\in\Isc^{(s),m}(O,\cM).$
This completes the inductive step.
\end{proof}

\begin{rem} If $\Rp$ is a sink or center, we use Proposition~\ref{HMV.108}.
If $\Rp$ is a saddle, we use Proposition~\ref{HMV.108pp}. In that case,
we take $O$ to be a $W$-balanced neighbourhood of $\Rp.$ Then
$u\in\Isc^{(s),m}(O,\cM)$ is a finite Legendre regularity condition,
and \eqref{eq:O_m-O} holds by propagation of Legendre regularity outside
critical points of the bicharacteristic flow, which is an immediate
consequence of the parametrix construction of Duistermaat and H\"ormander.
The propagation of Legendre regularity away from critical points
can also be proved by the use of test modules as in this section.
In fact, the use of \eqref{eq:O_m-O} can be avoided altogether
if we construct
$Q$ more carefully and use the
term with $\tilde B$ as well in \eqref{HMV.108e}. Namely, fix a small
neighbourhood $O'$ of $\Rp.$ Next, given $m,$
we first choose $K_\alpha>0$ so that $C'$ is positive definite at $\Rp$.
Then we construct $Q$ such that $\Rp\nin\WFsc(\Id-Q),$ $Q$ elliptic on $O'$
and \eqref{eq:[Q,P]} holds, and in addition
$(\tilde B^*\tilde B)_{|\alpha|=m}+Q^*C'Q$
(the first term denoting a diagonal matrix) can be written as $B^*B+F+G$
with $\WFscp(F)\cap \WFsc(u)=\emptyset,$ $G\in\Psisc^{-\infty,1}(X).$
We can do this since for any constant $N=N_m,$ we can arrange that
$Nq_0\leq-\scH_p q_0$ on a neighborhood of $\WFsc(u)$
outside a small $m$-dependent neighborhood of $\Rp$ where
$q_0=\sigma_\pa(Q).$ Then the proof of Proposition~\ref{HMV.108p} applies
directly and shows that $u\in\Isc^{(s),m}(O',\cM)$ for all $m.$
\end{rem}

\section{Outgoing eigenfunctions at a center}\label{Out.center}

We first analyze the structure of elements of $E_{\mic,\out}(\Rp,\ev)$ when
$\Rp\in\Min_+(\ev)$ is a center for the vector field $W,$ \ie when
\begin{equation}
2V''(\Cp)> \ev-V(\Cp),\ \Cp=\pi(\Rp).
\label{eq:Cp-center}\end{equation}
As shown in Corollary~\ref{cor:min-ident-8}, by
microlocalizing the solution near $\Rp$ we may simply assume
that
\begin{equation}
(P-\ev)u=f\in\dCI(X),\ u\in\CmI(X),\
\WFsc(u)\subset\{\Rp\}.
\label{HMV.71}\end{equation}

We proceed in three steps, first using the commutator methods outlined
above to obtain iterative regularity, with respect to the
$\Psisc^{-\infty,0}(O)$-module 
\begin{equation}
\cM\text{ generated by }\Id,\ \{x^{-\frac12}A,\ A\in
\Psisc^{-\infty,0}(O),\ \sigma_{\pa}(A)(\Rp)=0\}\Mand x^{-1}(P-\ev).
\label{HMV.72}\end{equation}
Then we pass to a blown up space where this regularity becomes conormal
regularity, once a phase is factored out. Finally, using the equation again,
we find the expansion \eqref{HMV.47}.

\begin{prop}\label{HMV.79}
The module $\cM$ in \eqref{HMV.72} is a test module in the sense of
Definition~\ref{Test.module}, which satisfies \eqref{eq:A_i-comm}.
If $u$ satisfies \eqref{HMV.71} at $\Rp$, 
then $u\in\Isc^{(-1/2-\ep)}(O,\cM)$ for all $\ep>0.$
\end{prop}

\begin{proof}
Choose a complex symbol $a\in\Cinf_c(\scT^*X)$ such that $da_j(\Rp)$,
restricted to $\Sigma(\ev)$ at $\Rp,$ is a non-trivial eigenvector for the
linearization of $W$ with eigenvalue $-2\tilde\nu r_1\in\bbC;$ here
$\tilde\nu=\nu(\Rp)$.
Then $\bar a$ has the
same property with eigenvalue $-2\tilde\nu r_2=-2\tilde\nu \overline{r_1}$
and together they
generate, over $\CI_c(\Sigma (\ev))$ and locally near $\Rp,$ the ideal of
functions vanishing at $\Rp.$ Recall from \eqref{HMV.37}
that $\Re r_1=\Re r_2=1/2.$ Certainly $A\in x^{-\frac12}\Psisc^{-\infty,0}(X)$
with symbol $x^{-\frac12}a$ is in the set in \eqref{HMV.72} and indeed
$\cM$ is generated by $A_0=\Id,$ $A_1=A,$ $A_2=A^*$ and $A_3=x^{-1}(P-\ev).$

By \eqref{HMV.46},
the commutator $i[A_j,P-\ev]\in\Psisc^{-\infty,-1/2}(X)$ has principal symbol
\begin{equation}
-\scH_p (x^{-1/2}a_j)=2\tilde\nu x^{-1/2}(i\Im(r_j) a_j+e_j),\ j=1,2,
\label{eq:comm-a_j}\end{equation}
where the $e_j$ vanish to second order at $\Rp.$ Thus
\begin{equation}
ix^{-1}[ A_j,P-\ev]=C_{j0}x^{1/2}+C_{j1}A_1+C_{j2} A_2+x^{1/2}C_{j3}A_3,
\ j=1,2,
\label{HMV.77}\end{equation}
with $C_{jk}\in\Psisc^{-\infty,0}(X),$ $\sigma_\pa(C_{jk})(\Rp)=0$ if $j\neq
k,$ and where $\sigma_\pa(C_{jj})(\Rp)=i2\tilde\nu\Im(r_j)$
is pure imaginary. In
particular $\cM$ is closed under commutators and hence is a test
module satisfying \eqref{eq:A_i-comm}. 

To prove the statement about $u$,
note that $\WFsc^{*,-1/2-\ep}(u) \cap O = \emptyset$ 
for any $W$-balanced neighbourhood $O$ of $\Rp$,
by Corollary~\ref{cor:micro-rough-prop}.
If we take any $Q\in\Psisc^{-\infty,0}(O)$
with $\Rp\nin\WFscp(\Id-Q)$, then \eqref{eq:[Q,P]-min-point} is satisfied,
so Proposition~\ref{HMV.108} yields the result. 
\end{proof}

\begin{rem}\label{rem:cancellation} 
The factor $x^{-1/2}$ in \eqref{HMV.72} is chosen precisely
that the real part of $r_j$ in \eqref{eq:comm-a_j} is cancelled, allowing
the application of Proposition~\ref{HMV.108}. A similar cancellation is
crucial in \eqref{eq:comm-a_j-2} in the next section.
\end{rem}

Next we reinterpret this iterative regularity more geometrically by
introducing the parabolically blown up space 
\begin{equation}
X_{\Cp}=[X;\{\Cp\}]_{1/2},\ \beta:X_{\Cp}\longrightarrow X.
\label{HMV.111}\end{equation}
Such parabolic blow-ups are defined in considerable generality, including
this case, in \cite{MR92i:32016}. We recall a form of the construction,
sufficient for our purposes here, in Section~\ref{sec:blow-up}. Recall
that $X_{\Cp}$ is a compact manifold with corners, of dimension two.
The two boundary curves are the front face $\ff=\beta ^{-1}(\Cp),$ created
by the blowup, and the closure of the pullback of $\pa X \setminus \{ \Cp
\}$ which we shall refer to as the old boundary. We shall denote boundary
defining functions for these boundary curves by 
$\rho_{\new}$ and $\rho_{\old}$, respectively. 
We also have boundary defining
functions for the other boundary curves of $X,$ which lift unchanged to
boundary curves of $X_{\Cp};$ we shall denote the product of these as $\rho
'$ since they play no significant r\^ole here.

One class of natural Sobolev spaces on such a manifold with corners is the
class of weighted $\bl$-Sobolev spaces, which we proceed to define. 
Let $\Vb(X_{\Cp})$ be the space
of smooth vector fields on $X_{\Cp}$ which are tangent to all boundary
faces and let $\mub$ be a b-density on $X_{\Cp}$, \ie, a density of the
form $(\rho'\rho_{\old}\rho_{\new})^{-1} \mu$, where $\mu$ is a smooth,
nonvanishing density. Then we set
\begin{equation}\begin{split}
(\rho')^a\rho _{\old}^b\rho _{\new}^c H^M_{\bl}(X_{\Cp})=
\left\{ u ; u = (\rho')^a\rho _{\old}^b\rho _{\new}^c v, \ v \in 
H^M_{\bl}(X_{\Cp}) \right\} \ , \\
H^M_{\bl}(X_{\Cp}) = 
\left\{u\in L^2(X_{\Cp};\mub);\ \Vb^{k}u\subset L^2(X_{\Cp};\mub) \ \forall \,
k \leq M \right\}. \end{split}
\label{HMV.114}\end{equation}
These are the conormal functions (conormal with
respect to the boundary) and we also use the short-hand notation 
\begin{equation}
(\rho')^{\infty}\rho^{\infty}_{\old}\rho_{\new}^{c}H^{\infty}_{\bl}(X_{\Cp})
=\bigcap\limits_{a,b}
(\rho')^{a}\rho^{b}_{\old}\rho_{\new}^{c}H^{\infty}_{\bl}(X_{\Cp})
\label{HMV.115}\end{equation}
for the subspaces which are rapidly decreasing up to the boundaries other
than $\ff.$

\begin{prop}\label{HMV.112} If $\Rp$ is a center for $W$ in $\Sigma(\ev),$
with $\cM$ the test module given by \eqref{HMV.72} in a $W$-balanced
neighbourhood of $\Rp$ then, for any $s,$ multiplication gives an isomorphism
\begin{multline}
\{u\in\CmI(X);\WFsc(u)\subset\{\Rp\},\ u\in I^{(s)}(O,\cM)\}\ni u\longmapsto\\ 
e^{-i\nu(\Rp)/x}u\in(\rho')^{\infty}\rho_{\old}^{\infty}\rho_{\new}^{2s+3/2}
H^{\infty}_{\bl}(X_{\Cp}).
\label{HMV.113}\end{multline}
\end{prop}

\begin{proof} We shall work in Riemannian normal coordinates in the
boundary based at the minimum under discussion. In terms of these
coordinates the boundary symbol of $P-\ev$ is 
\begin{equation}
p-\ev=\nu^2+\mu ^2+V_0(0)+y^2a(y)-\ev+O(x), \ y^2 a(y) = V_0(y) - V_0(0).
\label{HMV.116}\end{equation}
Here $\mu =0,$ $y=0$ and $\nu=\tilde\nu=\sqrt{\ev-V_0(0)}> 0$ and $a(0) =
V_0''(0)/2 > 0$ at the radial
point. The module $\cM$ is therefore generated by $\Id,$ $x^{\frac12}D_y,$
$x^{-\frac12}y$ and 
\begin{equation}
x^{-1}\left((x^2D_x)^2+x^2D_y^2 - \tilde\nu^2 +y^2a(y)\right).
\label{HMV.117}\end{equation}

Consider the effect of conjugation by $e^{i\tilde\nu/x}.$ This maps each
$\Psi^{m,l}_{\scl}(X)$ isomorphically to itself and hence 
\begin{equation}
\tcM=e^{-i\tilde\nu/x}\cM e^{i\tilde\nu/x}
\label{HMV.118}\end{equation}
is another test module, but in the open set 
\begin{equation}
\tilde O=\{(\nu - \tilde\nu,y,\mu);(\nu,y,\mu)\in
O\}\ni\tilde\Rp=(0,0,0) ,
\label{HMV.119}\end{equation}
$\tilde\Rp$ being the image of $\Rp.$ It is generated by the conjugates of the
generators, namely $\Id,$ $x^{\frac12}D_y,$ $x^{-\frac12}y$ and
\begin{equation}
\frac1x\left((x^2D_x-\tilde\nu)^2-\tilde\nu^2\right)+xD_y^2+x^{-1}y^2a(y).
\label{HMV.120}\end{equation}

The last two terms in \eqref{HMV.120} are in $\tcM^2.$ So, following the
observation \eqref{HMV.110}, we may drop them without changing the spaces
$I^{(s)}(\tilde O,\tcM).$ The new generator can then be written 
\begin{equation*}
\frac1x\left((x^2D_x-\tilde\nu)^2-\tilde\nu^2\right)=(xD_xx-2\tilde\nu)xD_x.
\end{equation*}
The factor $x D_x x - 2\tilde \nu\in\Psisc^{1,0}(X)$
is elliptic at the base point $(0,0,0)$
so can be dropped 
without changing the test module. We can also add the product of the other
two generators, again using \eqref{HMV.110}. Thus we see that 
\begin{multline}
I^{(s)}(O,\cM)=e^{i\tilde\nu/x}\cdot I^{(s)}(\tilde O,\cM')\Mwith\\ \cM'\text{
generated by }\Id,\ x^{\frac12}D_y,\ x^{-\frac12}y\Mand 2xD_x+yD_y.
\label{HMV.121}\end{multline}

Now, consider the effect of the parabolic blow up of $\Cp,$ passing from
$X$ to $X_{\Cp}$ on 
\begin{equation}
v\in x^{s}L^2_{\scl}(X)\Mwith \WFsc(v)\subset\{\tilde\Rp\},\ v\in
I^{(s)}(\tilde O,\cM').
\label{HMV.124}\end{equation}
As a neighbourhood of the front face we may take
\begin{equation}
[0,\delta )_\rho\times\overline{\bbR_Y}
\label{HMV.122}\end{equation}
where $\rho =(x+|y|^2)^{\frac12}=x^{\frac12}(1+|Y|^2)^{-\frac12},$
$Y=x^{-\frac12}y$ and $\overline{\bbR_Y}$ is the radial compactification of
the line (so that on $\overline{\bbR_Y}$, $1/|Y|$ is a defining function
for `infinity'). Thus we may take $\rho_{\ff} = \rho$ and $\rho_{\old} = (1
+ Y^2)^{-1/2}$ for explicit boundary defining functions. In terms
of coordinates $(\rho, Y)$, the generators of $\cM'$ in \eqref{HMV.121} become 
\begin{equation}
\Id, \ x^{\frac12}D_y=D_Y+Y(1+|Y|^2)^{-1}\rho D_\rho,\ 
x^{-\frac12}y=Y, \ 2xD_x+yD_y = \rho D_\rho.  
\label{HMV.123}\end{equation}
We may replace the second generator by just $D_Y$. Computing the 
pull-back of the Riemannian measure we see that 
in region \eqref{HMV.122} we have
\begin{equation}
Y^kD_Y^l(\rho D_\rho )^m \beta ^*v 
\in \rho_{\old}^{2s+2}\rho _{\new}^{2s+3/2}
L^2_{\bl}([0,\delta )_\rho\times\overline{\bbR_Y})\ \forall\ k,l,m\in\bbN_0.
\label{HMV.125}\end{equation}

Full iterative regularity with respect to $Y$ and $D_Y$ in a polynomially
weighted $L^2$ space is equivalent to regularity in the Schwartz
space, and hence to smoothness up to, and rapid vanish at, the old
boundary $Y=\infty.$ Thus \eqref{HMV.125} reduces to the statement that
$\beta ^*v$ is Schwartz in $Y$ with values in the conormal space in
$\rho = \rho_{\ff}$ which is precisely the content of \eqref{HMV.113}.
\end{proof}

\begin{thm}\label{thm:center-smooth}
Let $\Rp \in \Min_+(\ev)$ be a center for $W$,
\ie \eqref{eq:Cp-center} holds, and let $(x,y)$ be local coordinates centered
at $\Cp,$ with $y$ given by arclength along $\pa X$.
Let $X=x^{\frac12}$ and $Y=x^{-\frac12}y,$ let $h=V''(\Cp)$, and let
$Q$ be the operator
\begin{equation}
Q=D_Y^2-\frac{\tilde\nu}{2}(YD_Y+D_YY)+\frac12hY^2,\
\label{eq:Q}\end{equation}
on $L^2(\ff, dY)$, with eigenvalues
$\beta_j$ and normalized eigenfunctions $v_j,$ $j=1,2,\ldots$.  

Suppose that $u$ satisfies \eqref{HMV.71}. 
Then for all $N$ there exist sequences $\gamma_j,$ $\gamma_{j,k}$,
$1\leq k\leq N$ which are rapidly decreasing in $j$ (\ie for all $s$,
there exists $C_s>0$ such that $|\gamma_j|,
|\gamma_{j,k}|\leq C_s j^{-s},$ $j\geq 1$) and such that
\begin{equation}\begin{split}
u=e^{i\tilde\nu/x} v,&\quad v=\sum\limits_{j} X^{1/2}
X^{i(\beta_j+V_1(0))/\tilde \nu}
\left(\gamma _j+\sum\limits_{k=1}^N\gamma_{j,k}X^k+
\gamma'_{j,N}(X)\right)v_j(Y),\\
\text{ where }\ &\qquad \tilde\nu=\nu(\Rp),\ V_1(y) = \dbyd{V}{x} (0,y), \ 
\gamma'_{j,N}\in \cS(\bbR;X^{N}H^{\infty}_{\bl}([0,\delta))).
\end{split}\end{equation}
The map from microfunctions $u$ satisfying \eqref{HMV.71} to rapidly
decreasing sequences $\{\gamma_j\}$ is a bijection.
\end{thm}

\begin{proof}
Combining Propositions~\ref{HMV.79} and \ref{HMV.112} we have now
established the conormal regularity of the outgoing eigenfunctions at a
center on the resolved space $X_{\Cp}.$ To deduce the existence of
expansions as in \eqref{HMV.47} we consider the lift of the operator $P$ to
$X_{\Cp}.$ Using the coordinates $X=x^{\frac12}$ and $Y=x^{-\frac12}y,$
arising from Riemannian normal coordinates in the boundary, we find that
$\tilde P - \ev = e^{-i\tilde\nu/x} (P - \ev) e^{i\tilde\nu/x}$ takes the
form 
\begin{multline}
\tilde P - \ev =e^{-i\tilde\nu/x}\Big( (x^2D_x)^2+ix^3D_x+ x^2D^2_y+
x^2 A(x,y,x^2D_x)+ xB(x,y,xD_y)(xD_y) \\ 
+ x C(x,y)(x^2 D_x) (xD_y) +V(x,y) - \ev \Big) e^{i\tilde\nu/x} \\
= X^2\left(-\tilde\nu (XD_X+i/2)+Q+V_1(0)\right)+
X^3B'(x,y,XD_X,D_Y),
\label{HMV.127}\end{multline}
where $Q$ is as in \eqref{eq:Q},
$B'$ is a differential operator in $XD_X$ and $D_Y$ of order at most two
with coefficients smooth in $(x,y)$, and $V_1(y)=\pa_xV\big|_{x=0}.$

The operator $Q$ is a harmonic oscillator; in fact, conjugation by
$e^{i\tilde\nu Y^2/4}$ gives
\begin{equation}
e^{i\tilde\nu Y^2/4} Q e^{-i\tilde\nu Y^2/4} = D_Y^2 + \alpha^2 Y^2, \quad
\alpha = \sqrt{\frac{V_0''(0)}{2} - \frac{\tilde\nu^2}{4}} > 0,
\label{eq:alpha}\end{equation}
where $\alpha$ is real by \eqref{eq:Cp-center}. Thus, $Q$ has discrete
spectrum with eigenvalues $\beta_j = \alpha(2j+1)$, and eigenvalues $v_j
\in \cS(\bbR)$ which we take to be normalized in $L^2$.

The conclusion of \eqref{HMV.113} is that if $u$ is a microlocal
eigenfunction then $v=e^{-i\tilde\nu/x}u$ has conormal regularity on
$X_{\Cp}$ and vanishes rapidly at the old boundary. Even though
$X=x^{\frac12},$ and $Y=x^{-\frac12}y$ are singular coordinates at the old
boundary, the rapid vanishing there means that 
$v(Y,X)$ is Schwartz in $Y$ with values in a conormal
space $X^{r}H^{\infty}_{\bl}([0,\delta ))$ in $X.$ The remainder term in
\eqref{HMV.127} maps this space into
$\cS(\bbR_Y;X^{r+3}H^{\infty}_{\bl}([0,\delta))).$ Thus the condition
$(P-\ev)u\in\dCI(X)$ becomes the iterative equation 
\begin{equation*}\begin{split}
(-\tilde\nu (XD_X+i/2)+Q+V_1(\Cp))v=
XB'&(X^2,XY,XD_X,D_Y)v\\
&\in \cS(\bbR_Y;X^{r+1}H^{\infty}_{\bl}([0,\delta )).
\end{split}\end{equation*}
Since $v$ is Schwartz in $Y$ with values in a conormal
space $X^{r}H^{\infty}_{\bl}([0,\delta ))$ in $X,$
the eigenfunction expansion
\begin{equation*}
v=\sum\limits_{j}\gamma_j(X)v_j(Y)
\end{equation*}
converges in the space $\cS(\bbR_Y;X^{r}H^{\infty}_{\bl}([0,\delta
))),$ and the coefficients satisfy 
\begin{equation}
(\tilde\nu XD_X+i \tilde\nu /2-\beta _j-V_1(\Cp))\gamma _j(X)=f_j(X)\in
X^{r+1}H^{\infty}_{\bl}([0,\delta).
\label{HMV.131}\end{equation}
It follows that
\begin{equation*}
\gamma_j(X)=X^{1/2}X^{i(\beta_j+V_1(\Cp))/\tilde\nu}\gamma_j+\gamma '_j(X)
\end{equation*}
with
$\gamma _j$ constant and $\gamma'_j\in X^{r+1}H^{\infty}_{\bl}([0,\delta)).$
In fact the sequence $f_j$ is rapidly decreasing in $j$ with values in the
space in \eqref{HMV.131} so the same is true of the terms $\gamma _j'$
which come from integration. The $\gamma _j$ therefore also form a rapidly
decreasing sequence. These arguments can be iterated, giving the asymptotic
expansion, meaning that for any $N$  
\begin{equation}\begin{split}
&v=\sum\limits_{j} X^{1/2}
X^{i(\beta_j+V_1(\Cp))/\tilde\nu}
\left(\gamma _j+\sum\limits_{k=1}^N\gamma_{j,k}X^k
+\gamma'_{j,N}(X)
\right)v_j(Y),\\
&\qquad\qquad \gamma'_{j,N}\in\cS(\bbR;X^{N}H^{\infty}_{\bl}([0,\delta)))
\label{HMV.132}\end{split}\end{equation}
where the series in $j$ converges rapidly. If all the leading
constants $\gamma _j$ vanish then repeated integration shows that
$v\in\dCI(X_{\Cp})=\dCI(X)$ is rapidly decreasing.

Conversely the $\gamma _j$ in \eqref{HMV.132} form an arbitrary rapidly
decreasing sequence since $v$ may be constructed iteratively as in
\eqref{HMV.132} and then asymptotic summation, which can be made uniform in
the Schwartz parameter $j$ gives a corresponding microlocal eigenfunction.

This completes the proof of the theorem.
\end{proof}

\section{Outgoing eigenfunctions at a sink}\label{S.Out.sink}

In this section we consider the structure of microlocally outgoing
eigenfunctions at a radial point $\Rp \in \Min_+(\ev)$ which is a sink for
$W$, which occurs when
\begin{equation*}
\ev>\ev_{\Hess}(\Cp)=\ev_{\Hess}=V(\Cp)+2V''(\Cp), \quad \Cp = \pi(\Rp).
\end{equation*}
In section~\ref{S.Classical} it was shown that 
in this case the linearization of $W$ at $\Rp$ has two negative
eigenvalues, $-2\tilde\nu r_1$ and $-2\tilde\nu r_2$, where $0 < r_1 < 1/2 <
r_2 < 1$, $r_1 + r_2 = 1$ and $\tilde\nu = \nu(\Rp)$. The discussion
here closely parallels that in the previous section, so emphasis is placed
on the differences. At the end of the section we
also discuss the degenerate case, where $\ev=\ev_{\Hess}$ (and $r_1 = r_2 =
1/2$). 

As in the previous section, we may assume the microlocal eigenfunctions
satisfy \eqref{HMV.71}. We again 
use commutator methods to deduce iterative regularity, now
with respect to the
$\Psisc^{-\infty,0}(X)$-module 
\begin{multline}
\cM\text{ generated by }\Id,\ \{x^{-r_1}B,\ B\in
\Psisc^{-\infty,0}(X),\ \sigma_{\pa}(B)(\Rp)=0\},\\
\{x^{-r_2}B,\ B\in
\Psisc^{-\infty,0}(X),\ \sigma_{\pa}(B)|_{L_2}=0\}\Mand x^{-1}(P-\ev).
\label{HMV.72p}
\end{multline}
Recall here that $L_2$ is a smooth curve given by
Proposition~\ref{prop:smooth-Leg}.

We separate the resonant ($r_2/r_1 = 2, 3, 4, \dots$)
and non-resonant cases, although the proof for the
resonant case would go through in the non-resonant case as well, since
the proof in the latter case is more transparent.

\begin{prop}\label{HMV.79p}
Suppose $\Rp$ is a sink for $W$, with eigenvalues $-2\tilde\nu r_i$, as
above, and suppose that $r_2/r_1$ is not an integer. Then 
the module $\cM$ in \eqref{HMV.72p} is a test module in the sense of
Definition~\ref{Test.module}, which satisfies \eqref{eq:A_i-comm}.
If $u$ satisfies \eqref{HMV.71} at $\Rp$, 
then $u\in\Isc^{(-1/2-\ep)}(O,\cM)$ for all $\ep>0.$
\end{prop}

\begin{proof}
Choose real symbols $a_j\in\Cinf_c(\scT^*X)$ such that $da_j(\Rp)$ is,
restricted to $\Sigma(\ev)$ at $\Rp,$ a non-trivial eigenvector for the
linearization of $W$ with eigenvalue $-2\tilde\nu r_j\in\bbR$
and with $a_2=0$ on $L_2$ near $\Rp.$ Then $a_1$,
$a_2$ generate, over $\Cinf_c(\Sigma(\ev))$ and locally near $\Rp$,
the ideal of functions vanishing at $\Rp.$ Similarly, $a_2$ generates, over
$\Cinf_c(\Sigma(\ev))$ and locally near $\Rp,$ the ideal of functions which
vanish on $L_2.$ Then if $A_j\in \Psisc^{-\infty,-r_j}(X)$ have symbols
$x^{-r_j}a_j,$ $j=1,2,$ $\cM$ is generated by $A_0=\Id,$ $A_1,$ $A_2$ and
$A_3=x^{-1}(P-\ev)$ (notice that the operator with symbol $x^{-r_1} a_2$ is
in $\cM$, since $r_1 < r_2$). 

By \eqref{HMV.46},
the commutator $ix^{-1}[A_j,P-\ev]\in\Psisc^{-\infty,-r_j}(X)$ has
principal symbol 
\begin{equation}
- \scH_p (x^{-r_j}a_j)=x^{-r_j}e_j,\ j=1,2,
\label{eq:comm-a_j-2}\end{equation}
where the $e_j$ vanish to second order at $\Rp$ (here we exploit a
cancellation, similar to that in Remark~\ref{rem:cancellation}). 
Since $\scH_p$ is tangent to $L_2,$ $e_2$ in fact vanishes on
$L_2.$ Thus
\begin{equation}\begin{split}
&ix^{-1}[ A_1,P-\ev]=C_{10}x^{1-r_1}+C_{11}A_1+x^{r_2-r_1}C_{12} A_2
+x^{1-r_1}C_{13} A_3,\\
&ix^{-1}[ A_2,P-\ev]=C_{20}x^{1-r_2}+C_{22} A_2+x^{1-r_2}C_{23} A_3,\\
\end{split}\label{HMV.77p}\end{equation}
with $C_{jk}\in\Psisc^{-\infty,0}(X),$ $\sigma_\pa(C_{jk})(\Rp)=0$ for all
$j$ and all $k \geq 1$. 
Note that $[A_1,A_2] \in\Psisc^{-\infty,0}(X)$ as $r_1+r_2=1.$
So $\cM$ is closed under commutators, and hence is a test
module satisfying \eqref{eq:A_i-comm}. 

The statement about $u$ follows exactly as in Proposition~\ref{HMV.79}. 
\end{proof}

Next we deal with the case when $r_2/r_1=N$ is an integer, so $r_1^{-1} =
1/(N+1)$. 
Let $\mathcal{I}$
be the ideal of $\Cinf$ functions on $C_\pa$ vanishing at $\Rp$.
Now consider the
$\Psisc^{-\infty,0}(X)$-module 
\begin{multline}
\cM\text{ generated by }\Id,\ \{x^{-k/(N+1)}B,\ B\in
\Psisc^{-\infty,0}(X),\ \sigma_{\pa}(B)\in\mathcal{I}^k\},\ k=1,\ldots,N\\
\{x^{-N/(N+1)}B,\ B\in
\Psisc^{-\infty,0}(X),\ \sigma_{\pa}(B)|_{L_2}\in\mathcal{I}^{N}\}
\Mand x^{-1}(P-\ev).
\label{HMV.72pp}
\end{multline}

\begin{prop}\label{HMV.79pp}
Suppose $\Rp$ is a sink for $W$, with eigenvalues $-2\tilde\nu r_i$,
and suppose that $r_2/r_1 = N$ is an integer. Then 
the module $\cM$ in \eqref{HMV.72pp} is a test module in the sense of
Definition~\ref{Test.module}, which satisfies \eqref{eq:A_i-comm-p}.
If $u$ satisfies \eqref{HMV.71} at $\Rp$, 
then $u\in\Isc^{(-1/2-\ep)}(O,\cM)$ for all $\ep>0.$
\end{prop}

\begin{proof}
Assume first that $r_1<1/2<r_2$, \ie that $\ev > \ev_{\Hess}$. 
Choose real symbols $a_j\in\Cinf_c(\scT^*X)$ such that $da_j(\Rp)$ is,
restricted to $\Sigma(\ev)$ at $\Rp,$ a non-trivial eigenvector for the
linearization of $W$ with eigenvalue $-2\tilde\nu r_j\in\bbR,$ and $a_2$
vanishes on $L_2$. Then $a_1,$ $a_2$ generate, over
$\Cinf_c(\Sigma(\ev))$ and locally near $\Rp$,
the ideal of functions vanishing at $\Rp.$ Similarly, $a_2$ generates, over
$\Cinf_c(\Sigma(\ev))$ and locally near $\Rp$,
the ideal of functions vanishing on $L_2.$
Certainly $A_k\in \Psisc^{-\infty,-kr_1}(X)$
with symbol $x^{-k r_1}a_1^k,$ $k=1,2,\ldots,N$,
are in the set in \eqref{HMV.72pp}, as is $A_{N+1}$ with symbol
$x^{-r_2}a_2,$ and indeed
$\cM$ is generated by $A_0=\Id,$ $A_1,\ldots,A_{N+1}$ and $A_{N+2}=
x^{-1}(P-\ev).$

By \eqref{HMV.46},
the commutator $ix^{-1}[A_1,P-\ev]\in\Psisc^{-\infty,-r_1}(X)$ has
principal symbol 
\begin{equation*}
-\scH_p (x^{-r_1}a_1)=x^{-r_1}e_1,
\end{equation*}
and $ix^{-1}[A_{N+1},P-\ev]\in\Psisc^{-\infty,-r_2}(X)$ has
principal symbol 
\begin{equation*}
-\scH_p (x^{-r_2}a_2)=x^{-r_2}e_2,
\end{equation*}
where the $e_j$ vanish to second order at $\Rp.$
Since $W$ is tangent to $L_2$ to order $N,$ $e_2=e'_2+e''_2,$ $e'_2$
vanishing on $L_2$ and $e''_2\in\mathcal{I}^N.$ Thus
\begin{equation}
\begin{aligned}
ix^{-1}[ A_1,P-\ev]=&C_{10}x^{1-r_1}+C_{11}A_1+x^{r_2-r_1}C_{1,N+1} A_{N+1}\\
&+x^{1-r_1}C_{1,N+2} A_{N+2},\\
ix^{-1}[ A_k,P-\ev]=&C_{k0}x^{1-kr_1}+C_{kk}A_k+x^{r_2-k r_1}C_{k,N+1} A_{N+1}\\
&+x^{1-k r_1}C_{k,N+2} A_{N+2},\ k=2,3,\ldots,N\\
ix^{-1}[ A_{N+1},P-\ev]=&C_{20}x^{1-r_2}+C_{N+1,N} A_N+C_{N+1,N+1}A_{N+1}\\
&+x^{1-r_2}C_{N+1,N+2} A_{N+2},
\end{aligned}
\label{HMV.77pp}\end{equation}
with $C_{jk}\in\Psisc^{-\infty,0}(X),$
\begin{equation*}
\sigma_\pa(C_{jk})(\Rp)=0\ \text{for all}\ 
j\leq k\Mand k \geq 1,
\end{equation*}
but $\sigma_\pa(C_{N+1,N})(\Rp)$ usually does not vanish.
It is not hard to check that $\cM$ is
closed under commutators, and hence is a test module satisfying
\eqref{eq:A_i-comm-p}. 

The statement about $u$ follows exactly as in Proposition~\ref{HMV.79}. 

If $r_1=r_2=1/2,$ there is only one
eigenvector, which we may arrange to be $da_1(\Rp).$ We choose any
$a_2$ with $da_2(\Rp)$ linearly independent. The argument
as above still works, with $N=1.$
\end{proof}

\begin{rem} In the case that $r_2/r_1$ is an integer, $W$ is not
tangent to the curve $L_2$, and it is necessary to add extra generators
$A_2, \dots, A_N$ to make the module closed under commutators, as the proof
above shows.
Nevertheless, the enveloping algebra of $\cM$ in \eqref{HMV.72pp}
is the same as the enveloping
algebra of the module in \eqref{HMV.72p}, so for regularity considerations,
we work with \eqref{HMV.72p} below, instead of \eqref{HMV.72pp}, even when
$r_2/r_1$ is integral.
\end{rem}

We can again interpret this interative regularity more geometrically by
introducing a blown up space, although with blow-up of different
homogeneities in different variables, as discussed in
Section~\ref{sec:blow-up}. Namely, we consider the blow-up
of $X$ at $p$ along the vector field $\tilde V=r_1^{-1}x\pa_x+y\pa_y$:
\begin{equation}
X_{\Cp,r_1}=[X;\{\Cp\}]_{\tilde V},\ \tilde V=r_1^{-1}x\pa_x+y\pa_y,
\end{equation}
with blow-down map $\beta = \beta_{r_1}$. 
This space is different to that considered in the previous section (unless
$r_1 = 1/2$), but it `looks similar', and we use similar notation to
describe it. It is again a compact manifold with corners, with two boundary
curves, the  front face $\ff=\beta ^{-1}(\Cp)$ and the old boundary .
The two boundary curves are the front face $\ff=\beta ^{-1}(\Cp),$ created
by the blowup, and the closure of the pullback of $\pa X \setminus \{ \Cp
\}$ which we shall refer to as the old boundary. We shall again denote boundary
defining functions for these boundary curves by 
$\rho_{\new}$ and $\rho_{\old}$, and denote the product of 
boundary defining functions for the other boundary curves of $X$ by $\rho'$.

We next describe the appropriate conjugation. Since $L_2$ is Legendre and
has full rank projection to $\pa X$ near $\Rp,$ we can arrange that nearby
$L_2$ is given by $\nu=\Phi_2(y),$ $\mu=\Phi_2'(y),$ so inside
$\Sigma(\ev),$ $\mu-\Phi'_2(y)$ is a defining function for $L_2$.

\begin{prop}\label{HMV.112p} If $\Rp \in \Min_+(\ev)$ is a sink for $W$ in
$\Sigma(\ev),$ 
with $\cM$ the test module given by \eqref{HMV.72p} or \eqref{HMV.72pp}
in a $W$-balanced
neighbourhood of $\Rp$ then for any $s,$ multiplication gives an isomorphism
\begin{multline}
\{u\in\CmI(X);\WFsc(u)\subset\{\Rp\},\ u\in I^{(s)}(O,\cM)\}\ni u\longmapsto\\ 
e^{-i\Phi_2(y)/x}u\in(\rho')^{\infty}\rho_{\old}^{\infty}
\rho_{\new}^{(s+1)/r_1-1/2}
H^{\infty}_{\bl}(X_{\Cp,r_1}).
\label{HMV.113p}\end{multline}
\end{prop}

\begin{rem}
$\Phi_2$ can be replaced by any other smooth function $\tilde\Phi_2$ such
that $\Phi_2-\tilde\Phi_2=\mathcal{O}(y^{1/r_1})$.
\end{rem}

\begin{proof} As in the beginning of the proof of Proposition~\ref{HMV.112} we
have coordinates in which \eqref{HMV.116} holds. In this case the module
$\cM$ is therefore generated by $A_0=\Id,$ $A_1=x^{-r_1}y,$
$A_2=x^{-r_2}(xD_y-\Phi_2'(y))$ and $A_3=x^{-1}(P-\ev).$
Note that $\Phi_2(y)-\nu$ also vanishes on $L_2,$ hence
$x^{-r_2}(x^2D_x+\Phi_2(y))\in \cM$ as well.

Consider the effect of conjugation by $e^{i\Phi_2(y)/x}.$ This maps each
$\Psi^{m,l}_{\scl}(X)$ isomorphically to itself and hence 
\begin{equation}
\tcM=e^{-i\Phi_2(y)/x}\cM e^{i\Phi_2(y)/x}
\label{HMV.118p}\end{equation}
is another test module, but in the open set 
\begin{equation}
\tilde O=\{(\nu-\Phi_2(y),y,\mu-\Phi'_2(y));
(\nu,y,\mu)\in O\}\ni\tilde\Rp=(0,0,0),
\label{HMV.119p}\end{equation}
$\tilde\Rp$ being the image of $\Rp.$ It is generated by the conjugates of the
generators, namely $\Id,$ $x^{-r_1}y,$ $x^{-r_2}(xD_y)$ and
\begin{equation}
x^{-1} \tilde P = x^{-1} e^{-i\Phi_2(y)/x}(P-\ev)e^{i\Phi_2(y)/x}
\label{HMV.120p}\end{equation}

To compute this generator, we write $P - \ev$ as in \eqref{HMV.127} and
express \eqref{HMV.120p} as
\begin{equation}\begin{gathered}
x^{-1} \Big( (x^2 D_x - \Phi_2)^2 + ix(x^2 D_x - \Phi_2) + (x D_y + \pa_y
\Phi_2)^2  + x^2 A(x,y,x^2D_x - \Phi_2) \\
+ xB(x,y,xD_y + \pa_y \Phi_2)(xD_y + \pa_y \Phi_2)\\
+ x C(x,y)(x^2 D_x - \Phi_2)(xD_y + \pa_y
\Phi_2) + V(x,y) - \ev \Big) \\
= x^{-1} \Big( \Phi_2^2 + (\pa_y \Phi_2)^2 + V_0(y) - \ev \Big) \\
-2\Phi_2(xD_x) - i\Phi_2 + 2(\pa_y \Phi_2) D_y - i \pa^2_y \Phi_2 +
B(x,y,\pa_y \Phi_2)\pa_y \Phi_2\\
 -C(x,y)\Phi_2 \pa_y \Phi_2 + V_1(y)
+B_1 y^2D_y+B_2xD_y^2+B_3x^3D_x^2\\
+ B_4x^2D_x
+B_5xD_y+B_6y+B_7x,\ B_j\in\Diffsc^2(X).
\end{gathered}\label{eq:tilde-P}\end{equation}

In the nonresonant case, the term proportional to $x^{-1}$ vanishes
identically, since $\Phi_2$ satisfies the eikonal equation $\Phi^2 + (\pa_y
\Phi)^2 + V - \ev = 0$. Eliminating all terms which are in the span of
the generators $\Id$, $x^{-r_1}y$, $x^{-r_2}(xD_y)$ over
$\Psisc^{*,0}(X)$, as well as (following the
observation \eqref{HMV.110}) elements of $\cM^2$,
we are left with the generator
$$
-2\Phi_2(xD_x) + 2(\pa_y \Phi_2) D_y.
$$
A Taylor series analysis of $\Phi_2$ gives 
\begin{equation}
\Phi_2(y) = \tilde\nu(1 - r_1 y^2/2 + O(y^3)).
\label{eq:Phi_2-Taylor}\end{equation} 
Thus, $2(\pa_y \Phi_2) D_y$ is of the form $a(y) yD_y$. This may be written
$a(y)(x^{-r_1}y)(x^{r_1}D_y)$, the product of two
generators and a smooth function, so we may eliminate this term,
again following the
observation \eqref{HMV.110}. 
Thus, in the nonresonant case \eqref{HMV.120p} may be replaced by the
generator $xD_x$. This is true in the resonant case as well, since,
although the term proportional to $x^{-1}$ in \eqref{eq:tilde-P} does not
vanish, it is $O(y^{1/r_1})$, and $x^{-1} y^{1/r_1} \in \cM^{1/r_1}$. 
In summary, we have shown
\begin{multline}
I^{(s)}(O,\cM)=e^{i\Phi_2(y)/x}\cdot I^{(s)}(\tilde O,\cM')\Mwith\\ \cM'\text{
generated by }\Id,\ x^{-r_1}y,\ x^{-r_2}(xD_y)\Mand xD_x.
\label{HMV.121p}\end{multline}

Now, consider the effect of the inhomogeneous blow up of $\Cp,$ passing from
$X$ to $X_{\Cp}$ on 
\begin{equation}
v\in x^{s}L^2_{\scl}(X)\Mwith \WFsc(v)\subset\{\tilde\Rp\},\ v\in
I^{(s)}(\tilde O,\cM').
\label{HMV.124p}\end{equation}
As a neighbourhood of the front face we may take
\begin{equation}
[0,\delta )_\rho\times\overline{\bbR_Y}
\label{HMV.122p}\end{equation}
where $\rho =(x^{2r_1}+y^2)^{1/2},$
$Y=x^{-r_1}y$ and $\overline{\bbR_Y}$ is the radial compactification of
the line, where now $1/|Y|^{1/r_1}$
is a defining function for `infinity' on $\overline{\bbR_Y}$. In terms
of these coordinates, the generators of $\cM'$ in \eqref{HMV.121p} become 
\begin{equation}
\Id, \ x^{-r_1}y=Y, \ x^{r_1}D_y=D_Y+(1+|Y|^2)^{-1}Y\rho D_\rho, \ 
r_1(1+Y^2)^{-1}\rho D_\rho-r_1 Y D_Y. 
\label{HMV.123p}\end{equation}
A suitable linear combination of the last two generators, namely
$yD_y+r_1^{-1}xD_x$ from \eqref{HMV.121p}, is $\rho D_\rho.$
Using \eqref{HMV.110} as before to simplify the third generator and
computing the pull-back measure we see that in \eqref{HMV.122p}  
\begin{equation}
Y^kD_Y^l(\rho D_\rho )^m \beta ^*v\in \rho_{\old}^{(s+1)/r_1}
\rho _{\new}^{(s+1)/r_1-1/2}
L^2_{\bl}([0,\delta )_\rho\times\overline{\bbR_Y})\ \forall\ k,l,m\in\bbN_0.
\label{HMV.125p}\end{equation}

Full iterative regularity with respect to $Y$ and $D_Y$ in a polynomially
weighted $L^2$ space is equivalent to regularity in the Schwartz
space, and hence to smoothness up to, and rapid vanish at, the old
boundary $Y=\infty.$ Thus \eqref{HMV.125p} reduces to the statement that
$\beta ^*v$ is Schwartz in $Y$ with values in the conormal space in
$\rho_o$ which is precisely the content of \eqref{HMV.113p}.
\end{proof}

We now show that outgoing microlocal eigenfunctions have full asymptotic
expansions on $X_{\Cp, r_1}$. For this we recall from \cite{RBMCalcCon}
that if $I\subset\Real$ is discrete, $I\cap (-\infty,a]$ finite for $a\in\bbR$,
and $I+\bbN_0\subset I$ then
$\mathcal{A}^I_{\text{phg},\ff} (X_{\Cp, r_1})$ is the space of polyhomogeneous
conormal distributions with index set $I$, that vanish to infinite order
everywhere but on the front face $\ff$. That is, $u\in
\mathcal{A}^I_{\text{phg},\ff} (X_{\Cp, r_1})$ means there exist
$\phi_t\in\Cinf(X_{\Cp, r_1})$, $t\in I$, vanishing to infinite order
everywhere but on the front face, such that
\begin{equation*}
u\sim \sum_{t\in I}\phi_t \rho_{\new}^t.
\end{equation*}
Here the summation is understood as asymptotic summation, i.e.\ with the
notation of \eqref{HMV.125p}, the difference
$v=u-\sum_{t\in I,\ t\leq a}\phi_t \rho_{\ff}^t$ satisfies estimates
\begin{equation}
Y^kD_Y^l(\rho D_\rho )^m v\in \rho_{\new}^{a}
L^2_{\bl}([0,\delta )_\rho\times\overline{\bbR_Y})\ \forall\ k,l,m\in\bbN_0.
\label{HMV.125pp}\end{equation}

\begin{thm}\label{thm:sink-smooth}
Suppose that $\Rp \in \Min_+(\ev)$ is a sink of $W$, and
$u$ satisfies \eqref{HMV.71}. Let $\beta = r_2/2 +iV_1(\Cp)/2\tilde\nu$.
If $r_2/r_1$ is not an integer, then
\begin{equation}
u_0=x^{-\beta}e^{-i\Phi_2/x}u\in\mathcal{A}_{\text{phg},\ff}^I(X_{\Cp, r_1}),
\ I=\bbN_0+\bbN_0\frac{1}{r_1}+\bbN_0\frac{2r_2-1}{r_1}+\bbN_0\frac{r_2}{r_1}.
\label{eq:sink-smooth-8}\end{equation}
In particular, $u_0$ vanishes to infinite order off the front face,
and is continuous up to the front face.
The map sending microfunctions 
$u$ satisfying \eqref{HMV.71} to $u_0|_{\ff}$ is a bijection.

If $1/r_1$ is an integer, $r_1<1/2,$ then the same conclusions hold with
\eqref{eq:sink-smooth-8} replaced by
\begin{equation}
u_0=x^{-\beta}x^{icY^{1/r_1}}e^{-i\Phi_2/x}u,\quad Y=y/x^{r_1},
\label{eq:sink-smooth-8p}\end{equation}
where $c$ is a constant determined by the Taylor series, up to order $1/r_1$,
of $V$ and $g$ at $\Cp$.
\end{thm}

\begin{proof} Suppose first that $r_2/r_1$ is not an integer. 
Let $\delta = \min(2r_2 - 1, r_1) > 0$. 
We see from \eqref{eq:tilde-P} and \eqref{eq:Phi_2-Taylor} that
\begin{equation}
x^{-1} \tilde P  = -2\tilde \nu \big(xD_x + r_1yD_y + i\beta\big) + 
\rho_{\ff}^{\delta/r_1} P_2,
\label{eq:effective}\end{equation}
where $P_2$ is a differential operator of degree at most two generated by
vector fields tangent to the boundary of $X_{\Cp, r_1}$. By
Proposition~\ref{HMV.112p}, $v=e^{-i\Phi_2/x}u$ is invariant under such
vector fields. Thus, we deduce that 
\begin{equation*}
x^{-1} \tilde P  v =  -2\tilde\nu(xD_x + r_1yD_y + i\beta) v + 
\rho_{\ff}^{\delta/r_1} P_2 v \in \CIdot(X),
\end{equation*}
so with $x{D_x}_{|Y} = xD_x + r_1 y D_y$ denoting the derivative keeping $Y$,
rather than $y$, fixed, we have by Proposition~\ref{HMV.112p}
\begin{equation}
(x{D_x}_{|Y}+i\beta)v\in x^{-1/2+\delta'}L^2_\scl(X) \quad \forall \, 0 <
\delta' < \delta.
\label{eq:sink-smooth-32}\end{equation}
Writing $v=x^{\beta}\tilde v,$ and noting that $\Re \beta = r_2/2$, 
this yields
\begin{equation*}
x{D_x}_{|Y} \tilde v \in x^{-1/2 + \delta' - r_2/2} L^2_\scl(X) = 
x^{-1/2 + \delta' - r_2/2} L^2\big( \frac{dx dy}{x^3} \big).
\end{equation*}
Taking into account the smoothness in $Y,$ and changing the measure,
this means for each fixed $Y$,
\begin{equation*}
\pa_x \tilde v\in x^{-1/2+\delta'}L^2([0,1)_{x};dx) \subset
x^{\delta''}L^1([0,1)_{x};dx),\ \forall \delta'' < \delta'. 
\end{equation*}
But that implies that $\tilde v$ is continuous to $X=0,$ and after subtracting
$\tilde v(0,Y),$ the result is bounded by $C x^{\delta''}$.
This gives that, modulo $x^{-1/2+\ep'}L^2_\scl(X),$ $\ep'>0$ small,
microlocally near the critical point,
$u$ has the form $x^{\beta}e^{i\Phi_2(y)/x}u_0,$ where $u_0$ is smooth
on the blown-up space and rapidly vanishing off the front face.
A simple asymptotic series construction then yields the
asymptotic series described before, and then the uniqueness
result Proposition~\ref{prop:unique} 
shows that $u$ is actually given by such a series.

Now suppose that $r_2/r_1$ is an integer. We only need minor modifications
of the argument given above; we indicate these below. Equation
\eqref{eq:effective} must be replaced by 
\begin{equation*}
x^{-1} \tilde P  = -2\tilde \nu \big(xD_x + r_1yD_y + c x^{-1}
y^{1/r_1} +  i\beta\big) + x^\delta P_2,
\end{equation*}
where here $\delta = \min(2r_2 - 1, r_1) > 0$. 
Equation \eqref{eq:sink-smooth-32} must be replaced by 
\begin{equation}
(x{D_x}_{|Y}+ cY^{1/r_1} + i\beta)v\in x^{-1/2+\delta'}L^2_\scl(X) \quad
\forall 0 < \delta' < \delta.
\label{eq:replace}\end{equation}
Thus, we reach the same conclusion as before if we redefine
\begin{equation*}
\tilde v=x^{-\beta} x^{icY^{1/r_1}} v.
\end{equation*}
\end{proof}

\begin{rem}
First, we remark that if $1/r_1\nin\Nat$,
the form of the parameterization \eqref{eq:sink-smooth-8} is unchanged
if we replace $\Phi_2,$ parameterizing the 
smooth Legendrian $L_2$, by $\tilde\Phi_2\in\Cinf(X)$
with $\Phi_2-\tilde\Phi_2\in\mathcal{I}^{N},$ $N>1/r_1.$ Indeed, then
\begin{equation*}
(\Phi_2-\tilde\Phi_2)/x=(y/x^{r_1})^N x^{r_1 N-1} a,\quad a\in\Cinf(X),
\end{equation*}
so both the form of \eqref{eq:sink-smooth-8} and $u_0|_{\ff}$ are unaffected
by this change of phase functions. In particular, if $r_1>1/N,$ but we let
$r_1\to 1/N,$ the factor $x^{r_1 N-1}$ decays less and less,
and indeed in the limit it takes a different form \eqref{eq:sink-smooth-8p}.

If now we allow $\Phi_2-\tilde\Phi_2\in\mathcal{I}^{N}$,
$|r_1^{-1}-N|<1$ for some $N\in\Nat,$ then
\eqref{eq:effective} is replaced by
\begin{equation}
x^{-1}e^{-i\tilde\Phi_2/x}(P-\ev)e^{i\tilde\Phi_2/x}=-2\tilde\nu(xD_x+r_1yD_y)
+cy^N/x+i\beta+ \rho_{\ff}^{\delta}P_2,
\label{eq:sink-smooth-24pp}\end{equation}
where $P_2$ is in the enveloping algebra of $\tcM$ and $\delta > 0$.
(Note that the $y^N/x$ term can be included in $P_2$ if $N>1/r_1,$ but not
otherwise, and we are interested in taking $r_1\to 1/N$!)
Proceeding with the argument as in the resonant case,
\eqref{eq:replace} is replaced by
\begin{equation}
(x{D_x}_{|Y}+cx^{r_1N-1}Y^{N}+i\beta)v\in x^{-1/2+\delta}L^2_\scl(X).
\label{eq:sink-smooth-32pp}\end{equation}
If $r_1N-1\neq 1$,
we can remove the $(cx^{r_1N-1}Y^{N}+i\beta)$ terms by introducing
an integrating factor, namely by writing
$v=x^{-\beta}e^{-ic(r_1N-1)^{-1} x^{r_1N-1}Y^N}\tilde v$,
to obtain that
$\tilde v$ is continuous to $x=0,$ and the asymptotics take the form
\begin{equation}
u=x^{\beta}e^{-ic(r_1N-1)^{-1} x^{r_1N-1}Y^N}
e^{i\Phi_2/x}u_0,\quad Y=y/x^{r_1},
\label{eq:sink-smooth-8pp}\end{equation}
with $u_0$ continuous up to the front face of $X_{\Cp,r_1}$,
and \eqref{eq:sink-smooth-8p} is the limiting case as $r_1 N\to 1$.
In particular, the special role played by $1/r_1\in\Nat$
in the statement of the theorem is partly due to its formulation.
Here if $N>1/r_1,$ the factor $e^{-ic(r_1N-1)^{-1} x^{r_1N-1}Y^N}$ tends
to $1$ in the interior of $\ff,$ i.e.\ does not affect the form of
the asymptotics, but if $N<1/r_1,$ it introduces an oscillatory factor with
phase $-c(r_1N-1)^{-1}y^N/x$.
\end{rem}

Finally, we analyze what happens at the Hessian threshold, \ie when
$r_1=r_2=1/2.$ The correct space to describe the asymptotics
is not quite the parabolic blow-up $[X;\{\Cp\}]_{\frac 12},$ on which
$Y=y x^{-1/2}$ is smooth, rather the space on which
\begin{equation*}
\frac{y}{x^{1/2}\log x}=\frac{Y}{\log x}
\end{equation*}
is smooth.
To simplify the statement we only consider the top order
of the asymptotics.

\begin{thm} Suppose that $\Rp \in \Min_+(\ev)$, and $\ev$ is at the Hessian
threshold for $\Rp$, so that $r_1(\Rp)=r_2(\Rp)=\frac 12$. Let 
$\beta$ as in Theorem~\ref{thm:sink-smooth}, let $\Phi_2(y)=
\tilde\nu(1 - y^2/4)$, where $\tilde\nu=\nu(\Rp)$, and let $Y = x^{-1/2}y$.
Suppose that $u$ satisfies \eqref{HMV.71}. Then
\begin{equation}
u_0=x^{-\beta}(\log x)^{1/2}
e^{-i\Phi_2/x}e^{i\frac{\tilde\nu}4\,\frac{Y^2}{\log x}}u
\label{eq:logarithmic}\end{equation}
is such that
\begin{equation}
u_0-g(Y/\log x)\in (\log x)^{-1} L^\infty(X) \text{ for some }g\in\cS(\bbR).
\label{eq:sink-smooth-8q}\end{equation}
The map from microfunctions 
$u$ satisfying \eqref{HMV.71} to $g$ is a bijection.
\end{thm}

\noindent
The full asymptotic expansion, which is not stated above,
is in terms of powers of $\log x$ and it arises from
a stationary phase argument as can be seen from the proof given below.
One could
instead give a description of $u$ as an oscillatory integral, analogous
to how $u$ is described as an oscillatory sum in
Theorem~\ref{thm:center-smooth}, with an error term in
$x^{\frac12-\delta}L^\infty$ for all $\delta>0.$ In particular, not only
are $u,$ hence all terms of the (logarithmic)
asymptotic expansion, determined by $g,$
but the terms only depend on the principal symbol of $P-\ev$ near $\Rp,$
modulo $\mathcal{I}^2,$ and the subprincipal symbol at $\Rp.$

\begin{proof}
We use the computations of Theorem~\ref{thm:center-smooth}. Thus, with
$\tilde P=e^{-i\tilde\nu/x}Pe^{i\tilde\nu/x},$ and $X = x^{1/2}$,
\begin{multline}
\tilde P = X^2\left(-\tilde\nu (XD_X+i/2)+Q+V_1(\Cp)\right)+
X^3B'(X,XD_X,Y,D_Y),\\ 
\text{with }Q=e^{i\tilde \nu Y^2/4} D_Y^2 e^{-i\tilde \nu Y^2/4},
\label{HMV.127p}\end{multline}
since $\alpha$ in \eqref{eq:alpha} is zero at the Hessian threshold.
Thus
\begin{multline}
\tilde P'  =e^{-i\Phi_2/x}(P-\ev)e^{i\Phi_2/x}\\
=X^2\left(-\tilde\nu (XD_X+i/2)+D_Y^2+V_1(\Cp)\right)+
X^3B''(X,XD_X,Y,D_Y).
\end{multline}
We thus deduce that
\begin{equation*}
X^{-3}\left(\tilde P' -X^2(-\tilde\nu XD_X+D_Y^2 -
2i\tilde\nu\beta)\right) 
\end{equation*}
is in the enveloping algebra of $\tcM.$ Thus, with $u=e^{i\Phi_2/x} 
X^{\beta}v,$
\begin{equation*}
(-\tilde\nu XD_X+D_Y^2)v\in x^{1-\delta}L^2_{\scl}(X),\ \forall\delta>0.
\end{equation*}
The operator on the left hand side is the Schr\"odinger operator $-D_t+D_Y^2$
after a logarithmic change of variables $t=\frac1{\tilde\nu}\log X.$
This can now be solved explicitly by taking the Fourier transform in $Y,$
using the fact that $v$ is Schwartz in $Y$ for $X$ bounded away from $0.$
The solution thus has the form of the inverse Fourier transform of
$X^{i\eta^2/\tilde\nu}\tilde g(\eta),$ plus
faster vanishing terms, where $\tilde g$ is Schwartz, and $\eta$ is
the dual variable of $Y$. This gives an expression of the form
\eqref{eq:logarithmic}. 
\end{proof}

\section{Outgoing eigenfunctions at a saddle}\label{sect:saddle}
We next analyze microlocal eigenfunctions $u \in E_{\mic,+}(\Rp,\ev)$ when
$q \in \Max_+(\ev)$, hence $q$ is an outgoing saddle point for $W$. Recall
from section~\ref{S.Classical} that in this case there are always two
smooth Legendre curves through $q$ which are tangent to $W$, $L_1$ and
$L_2$, and there is a local coordinate $v_i$ on $L_i$ with $v_i=0$ at
$q$ such that $W$ takes the form $-2\tilde\nu r_1 v_2 \pa_{v_2}$ on $L_2$ and 
$-2\tilde\nu r_2 v_1 \pa_{v_1}$ on $L_1$. Here $\tilde\nu
= \nu(q)$ is positive, and $r_i$
satisfy $r_1 < 0$ and $r_2 > 1$. Hence $L_2 \subset \Phi_+(\{q\})$ is the
outgoing Legendrian in this case, and by microlocalizing 
$u \in E_{\mic,+}(\Rp,\ev)$ to a $W$-balanced neighbourhood $O$ of $q$, we
may assume that  
\begin{equation}
u\in\CmI(X),\ O\cap\WFsc((P-\ev)u)=\emptyset,\ \WFsc(u)\cap O\subset L_2.
\label{HMV.71r}\end{equation}

As before we first use commutator methods to deduce iterative regularity, now
with respect to the
$\Psisc^{-\infty,0}(X)$-module 
\begin{equation}
\cM\text{ generated by }\Id,\ \{x^{-1}B,\ B\in
\Psisc^{-\infty,0}(X),\ \sigma_{\pa}(B)|_{L_2}=0\}.
\label{HMV.72r}
\end{equation}

In particular, $x^{-1}(P-\ev)\in\cM$, since 
$\sigma_{\pa}(P-\ev)$ vanishes on $L_2$.

\begin{prop}\label{HMV.79r}
The module $\cM$ in \eqref{HMV.72r} is a test module in the sense of
Definition~\ref{Test.module}, which satisfies \eqref{eq:A_i-comm}.
If $u$ satisfies \eqref{HMV.71r} at $\Rp$, 
then $u\in\Isc^{(-1/2-\ep)}(O,\cM)$ for all $\ep>0.$
That is, $u$ is a Legendre distribution microlocally near $\Rp$.
\end{prop}

\begin{proof}
Choose real symbols $a_j\in\Cinf_c(\scT^*X)$ such that $a_j$ vanishes on
$L_j$ in a neighbourhood of $q$, and such that $da_j(\Rp)$ is,
restricted to $\Sigma(\ev)$ at $\Rp,$ a non-trivial eigenvector for the
linearization of $W$ with eigenvalue $-2\tilde\nu r_j\in\bbR$.
Then $a_1,$ $a_2$ generate, over
$\Cinf_c(\Sigma(\ev))$ and locally near $\Rp$,
the ideal of functions vanishing at $\Rp.$ Similarly, $a_i$ generates, over
$\Cinf_c(\Sigma(\ev))$ and locally near $\Rp$,
the ideal of functions also with differentials vanishing on $L_i.$
Certainly $A_1\in x^{-1}\Psisc^{-\infty,0}(X)$
with symbol $x^{-1}a_2,$ is in the set in \eqref{HMV.72r} and indeed
$\cM$ is generated by $A_0=\Id,$ $A_1$ and $A_2=x^{-1}(P-\ev).$

By \eqref{HMV.46},
the commutator $i[A_1,P-\ev]\in\Psisc^{-\infty,-1}(X)$ has principal symbol
\begin{equation*}
-\scH_p (x^{-1}a_2)=x^{-1}(2\tilde\nu(r_2-1)a_2+e_2).
\end{equation*}
where $e_2$ vanishes to second order at $\Rp.$
Since $\scH_p$ is tangent to $L_2,$ $e_2$ in fact vanishes on
$L_2.$ Thus
\begin{equation}
ix^{-1}[ A_1,P-\ev]=C_{10}+C_{11}A_1
+C_{12} A_2,
\label{HMV.77r}\end{equation}
with $C_{jk}\in\Psisc^{-\infty,0}(X),$ $\sigma_\pa(C_{12})(\Rp)=0$,
and $\sigma_\pa(C_{11})(\Rp)=2\tilde\nu(r_2-1)>0$.
In particular $\cM$ is
closed under commutators, and hence is a test module satisfying
\eqref{eq:A_i-comm}.

To prove the statement about $u$,
note that $\WFsc^{*,-1/2-\ep}(u) \cap O = \emptyset$ 
for any $W$-balanced neighbourhood $O$ of $\Rp$,
by Corollary~\ref{cor:micro-rough-prop}. To apply
Proposition~\ref{HMV.108p}, we need to construct a $Q$ satisfying
\eqref{eq:[Q,P]}. We take $Q\in\Psisc^{-\infty,0}(X)$ such that
\begin{equation}
\sigma_{\pa}(Q)=q= \chi_1(a_1^2)\chi_2(a_2)\psi(p),
\end{equation}
where $\chi_1,\chi_2,\psi\in\Cinf_c(\Real),$ $\chi_1,\chi_2\geq 0$ are
supported 
near $0,$ $\psi$ supported near $\sigma,$ $\chi_1,\chi_2\equiv 1$ near $0$
and $\chi_1'\leq 0$ in $[0,\infty).$ Note that $a_2=0$ on $L_2,$ so
$\supp d(\chi_2\circ a_2)\cap L_2=\emptyset$.
On the other hand,
\begin{equation}
\scH_p \chi_1(a_1^2)=2 a_1(\scH_p a_1)\chi_1'(a_1^2)
=-4\tilde\nu a_1(r_1 a_1+e_1)\chi_1'(a_1^2),
\end{equation}
with $e_1$ vanishing quadratically at $\Rp.$ Moreover, on $\supp\chi_1'$,
$a_2$ is bounded away from $0$, and $r_1<0$, so $r_1a_1+e_1\neq 0$,
provided $\supp\chi_1'$
is sufficiently small.
Since $r_1<0$ and $\chi_1'\leq 0$,
$4\tilde\nu r_1 a_1^2\chi_1'(a_1^2)$ is positive. Hence, choosing
$\chi_1$ such that $(-\chi_1\chi_1')^{1/2}$ is $\Cinf,$ and $\supp\chi_1$ is
sufficiently small, we can write
\begin{equation}\begin{split}
&\sigma_{\pa}(i[Q^*Q,P-\ev])=-\scH_p q^2=4\tilde\nu \tilde b^2+\tilde f,\\
&\ \tilde b=(a_1(r_1a_1+e_1)\chi_1'(a_1^2)\chi_1(a_1^2))^{1/2}\chi_2(a_2)
\psi(p),\ \supp\tilde f\cap L_2=\emptyset.
\end{split}\label{eq:Q-saddle}\end{equation}
Now we can apply Proposition~\ref{HMV.108pp}, using the remark that
follows it, to finish the proof.
\end{proof}

Let $\Phi_2$ parameterize $L_2$ near $\Rp,$ as in the previous section.
Multiplication by $e^{-i\Phi_2(y)/x}$ 
maps $L_2$ to the zero section $\nu=0,$ $\mu=0.$ The set of
Legendre distributions
associated to the zero section is exactly that of
distibutions  conormal to the boundary. The test module corresponding
to this class (which is the conjugate of $\cM$ under $e^{-i\Phi_2/x}$),
is $\tcM=\Vb(X),$ generated by $xD_x$ and $D_y$,
microlocalized in $\tilde O,$ the image of $O$ under
the multiplication. We thus deduce the
following corollary.

\begin{cor}\label{HMV.79c}
Suppose that $u$ satisfies \eqref{HMV.71r} for $\Rp \in \Max_+(\ev)$. 
Then $e^{-i\Phi_2(y)/x}u\in\Isc^{(-1/2-\ep)}(\tilde O,\tcM)$ for
all $\ep>0,$ \ie its microlocalization to $\tilde O$ is in
$x^{1/2-\ep}H_{\bl}^\infty(X)$ for all $\ep>0$.
\end{cor}

\begin{rem}\label{rem:conormal}
The $L^2$ conormal space $x^s H_{\bl}^\infty(X)$ is contained in the weighted
$L^\infty$ space with the weight shifted by $\epsilon$; that is,
$v\in x^sH_{\bl}^\infty(X)$
gives $v\in x^{s-\ep}L^\infty(X)$ for all $\ep>0$ (by Sobolev embedding). 
Conversely, if $v\in x^{m}L^\infty(X)\cap x^sH_{\bl}^\infty(X)$
then by interpolation, $v\in x^{m-\ep}H_{\bl}^\infty(X)$ for all $\ep>0$.
\end{rem}

The equation $(P-\ev)u\in\dCI(X)$ for Legendre functions $u$ reduces
to a transport equation. We can obtain this transport equation
by a very explicit conjugation of $P-\ev$ as the following proposition
shows. 

\begin{prop}\label{prop:saddle-L_2-model}
There exists an operator $A$ of the form 
$$
A = e^{i\Phi_2/x} \circ F^* \circ x^{\beta} e^{if},
$$
where $\beta = r_2/2 +iV_1(\Cp)/2\tilde\nu$, $f \in \Cinf(X)$, and $F^*$ is
pullback by a diffeomorphism $F$ on $X$, such that 
\begin{equation}
A^{-1}(P-\ev)A-\tilde P_0\in x^2\tcM^2,\quad \tilde P_0=-2\tilde\nu
((x^2D_x)+r_1 y(xD_y)).
\end{equation}
\end{prop}

\begin{proof}
By \eqref{eq:tilde-P} and \eqref{eq:Phi_2-Taylor}, we have
\begin{multline*}
e^{-i\Phi_2/x} (P - \ev) e^{i\Phi_2/x} = \\
-2\tilde\nu \big( (x^2D_x)+r_1 y(xD_y) + ix\beta(y) \big) +
yF(y)(x^2D_x)+y^2G(y)(xD_y) \ \operatorname{mod}\ x^2\tcM^2,
\end{multline*}
where $\beta(0) = \beta$. 
Since $r_1 < 0$, the vector field $x \pa_x + r_1 y \pa_y$ is nonresonant. 
Hence by a change of coordinates $x'=a(y)x,$ $y'=b(y)y$,
we may arrange that the vector field becomes
$$
-2\tilde\nu(((x')^2D_{x'})+r_1 y(x'D_{y'})),
$$
modulo terms in $x^2\tcM^2$ and subprincipal
terms. We therefore deduce that there is a local diffeomorphism $F$ such
that 
\begin{equation}\begin{split}
&(F^{-1})^* e^{-i\Phi_2/x}(P-\ev)e^{i\Phi_2/x} F^*
-\tilde P_1 \in x^2\tcM^2,\\
&\tilde P_1= -2\tilde\nu((x^2D_x)+r_1 y(xD_y)+ ix\beta(y)).
\end{split}\end{equation}
Now we let $\tilde\beta(y)=(\beta(y)-\beta(0))/y$ and set 
$f(y) = \int_0^y \tilde\beta(z)\,dz$. This yields
\begin{equation}
x^{-\beta} e^{-if} \tilde P_1 e^{if} x^{\beta} -\tilde P_0 \in x^2\tcM^2,
\end{equation}
which completes the proof of the proposition.
\end{proof}

\begin{rem} We have not conjugated the operator $P-\ev$ microlocally to
$\tilde P_0$; we have only conjugated it to the correct form `along $L_2$'
-- this suffices because of the a priori conormal estimates.  Note also
that $e^{-i\Phi_2/x}u\in x^s H_{\bl}^\infty(X)$ is equivalent to
$A^{-1}u\in x^{s+r_2/2}H_\bl^\infty(X).$
\end{rem}

Note that $(yx^{-r_1})^n,$ for $n\geq 0$ integral, satisfies $\tilde
P_0(yx^{-r_1})^n=0,$ and $(yx^{-r_1})^n\in x^s H_{\bl}^\infty(X)$ for all
$s<-r_1 n$ (where $s$ is increasing with $n$, since $r_1 < 0$). 
These can be modified to obtain microlocal
solutions of $A^{-1}(P-\ev)A v\in\dCI(X)$.

\begin{prop}\label{prop:saddle-approx-soln}
For each integer $n\geq 0,$ there exists $v_n\in \cap_{s<-r_1 n}
x^s H_{\bl}^\infty(X)$ such that 
$v_n-(yx^{-r_1})^n\in x^{s'} H_{\bl}^\infty(X)$ for some $s'>-r_1 n$,
and such that $\Rp\nin\WFsc(A^{-1}(P-\ev)A v_n)$.
\end{prop}

\begin{proof}
Applying
$A^{-1}(P-\ev)A$ to $x^\alpha y^n$ yields
a function of the form $x^{\alpha+2} f',$ $f'\in\Cinf(X).$ The equation
$\tilde P_0 x^{\alpha+1} g=x^{\alpha +2}f'$ has a smooth
solution $g$ modulo $x^{\alpha+3} f'',$ $f''\in\Cinf(X),$ provided $\alpha+1$
is not an integer multiple of $r_1$ (and then only a logarithmic factor
in $x$ is needed). Applying this argument with $\alpha = -r_1 n$, etc.,
we deduce that the error terms arising from $(yx^{-r_1})^n$ can
be solved away iteratively.
\end{proof}

This result shows that $A^{-1}(P-\ev)A,$ hence $P-\ev,$ have a series of
approximate generalized eigenfunctions $v_n,$ resp.
$u_n=Av_n,$ where $u_n$ has the form
$u_n=e^{i\Phi_2/x} x^{\beta - r_i n} \tilde v_n,$ with $\tilde v_n$
polyhomogeneous conormal and continuous up to the boundary.

We proceed to show that for any $u$ satisfying \eqref{HMV.71r}
there exist constants
$a_n$ such that $u \sim \sum_{n=0}^\infty a_n u_n.$ The main technical
result is the following lemma.

If $U\subset X$ is open, we denote by $ x^s H_{\bl}^\infty(U)$ the space
of functions consisting of the restrictions of $ x^s H_{\bl}^\infty(X)$
to $U$.

\begin{lemma}\label{lemma:saddle-vf-soln}
Suppose that $U_0$ is a neighbourhood of $\Cp$ in $X$.
There exists a neighbourhood $U\subset U_0$ of $\Cp$ such that
the following hold.

Suppose that
$v|_{U}\in x^{r}H_\bl^\infty(U)$ and
$x^{-1}\tilde P_0 v|_{U}\in x^s H_{\bl}^\infty(U),$ $s> r$. Then

\begin{enumerate}
\item
If $s\leq 0$ then $v|_{U}\in x^{s-\ep}H_\bl^\infty(U)$ for all $\ep>0$.

\item
If $s>0,$ then there exists a constant $a_0$
such that $v|_{U}-a_0\in x^{r'}H_\bl^\infty(U)$ for any
$r'<\min(-r_1,s)$.

\item
If $r>0,$ then $v|_{U}\in x^{r'}H_\bl^\infty(U)$ for all
$r'<\min(-r_1,s)$.

\item
More generally, if $n\geq 0$ is an integer, $r>-(n-1)r_1,$ $s>-nr_1$,
then there exists a constant $a_n$ such that
$v|_{U}-a_n(x^{-r_1}y)^n
\in x^{r'}H_\bl^\infty(U)$ for any $r'<\min(-(n+1)r_1,s)$.

\item
If $s>r,$ $r>(n-1)r_1,$ $n\geq 0$ an integer,
then $v|_{U}\in x^{r'}H_\bl^\infty(U)$ for all $r'<\min(-n r_1,s)$.
\end{enumerate}
\end{lemma}

\begin{proof}
Let $I$ be a small open interval, $0\in I,$ and let $x_0>0$ be small, so
that $U=[0,x_0]_x\times I_y\subset U_0$ is a coordinate neighbourhood of $\Cp$.
Thus, $v|_U\in x^{r}H_\bl^\infty(U)$ and
$x^{-1}\tilde P_0 v|_U\in x^s H_{\bl}^\infty(U),$ $s> r$.
The integral curves of the vector field $\tilde P_0$ are given by
$$
x^{-r_1} y=\text{constant}.
$$
Since $-r_1>0$, if $x$ is increasing along an
integral curve then $|y|$ is decreasing, and vice versa. Hence,
the integral curve $\gamma$ through $(x,y)\in U$ satisfies $\gamma(T)\in U$
provided $x_0\leq x(\gamma(T))\leq x$.

Since $\tilde P_0$ is a vector field, $v-v|_{x=x_0}$ is given by the
integral of $f=x^{-1}\tilde P_0 v$
along integral curves of $\tilde P_0$.
Namely, the solution is given by
$$
v(x,y)=-(2\tilde\nu)^{-1}
\int_{x_0}^x f(t,\left(\frac{x}{t}\right)^{-r_1} y)\,\frac{dt}{t}
+v(x_0,\left(\frac{x}{x_0}\right)^{-r_1} y).
$$
Note that the second term, which solves the homogeneous equation,
is certainly a polyhomogeneous function on $(x,y)$
down to $x=0,$ since it only evaluates $v$ at $x=x_0.$ Thus, it has
a full asymptotic expansion, corresponding to the Taylor series of $v$
at $x=x_0$ around $y=0$ of the form $\sum_{j=0}^\infty c_j (x^{-r_1}y)^j$.

So suppose that $s\leq 0$ first. Then the
first term (when restricted to $U$)
is in $x^{s-\delta}L^\infty$ for all $\delta>0$ since $f\in x^sL^\infty$
along the restriction of $\gamma$ to $x\leq x\circ \gamma\leq x_0.$ Hence,
$v\in x^{s-\delta}L^\infty$ for all $\delta>0,$ so
$v\in x^{s-\ep}H_\bl^\infty(X)$
for all $\ep>0$.

On the other hand, suppose that $f\in x^s L^\infty$ for some $s>0$.
The first term gives a convergent integral, hence is $L^\infty$; indeed,
it is continuous to $x=0$ with limit
$\int_{x_0}^0 f(t,0)\,\frac{dt}{t}$.
Let $a_0 = v(x_0,0)+\int_{x_0}^0 f(t,0)\,\frac{dt}{t}$.
Then $v-a_0\in x^{\min(-r_1,s-\ep)}L^\infty(X)$
for all $\ep>0,$ so
$v-a_0\in x^{\min(-r_1,s)-\ep}H_{\bl}^\infty(X)$ by
Remark~\ref{rem:conormal}. To see this, we write
$$
v - a_0 = \int\limits_{x_0}^x [f(t,0)-f(t,\left(\frac{x}{t}\right)^{-r_1}
y)]\,\frac{dt}{t} 
=\int\limits_{x_0}^x \int\limits_0^1 \left(\frac{x}{t}\right)^{-r_1} y
\frac{\pa f}{\pa y}(t,\rho \left(\frac{x}{t}\right)^{-r_1}
y)\,d\rho\,\frac{dt}{t}, 
$$
and estimate the integrand directly. A similar argument works for arbitrary
$n$.
\end{proof}

\begin{cor}\label{cor:saddle-vf-soln}
Let $v_n,$ $n\geq 0$ integer, be given by
Proposition~\ref{prop:saddle-approx-soln}.
Suppose that $U_0$ is a neighbourhood of $\Cp.$ There exists a neighbourhood
$U\subset U_0$ of $\Cp$ with the following properties.

Suppose that
$v\in x^{r}H_\bl^\infty(U)$ and
$x^{-1}A^{-1}(P-\ev)A v\in \dCI(U)$.
Then there exists a constant $a_0$
such that $v-a_0 v_0\in x^{r'}H_\bl^\infty(X)$ for any $r'<-r_1$.
More generally, if $n\geq 0$ is an integer, $r>-(n-1)r_1$,
then there exists a constant $a_n$ such that
$v-a_n v_n\in x^{r'}H_\bl^\infty(X)$ for any $r'<-(n+1)r_1$.
\end{cor}

\begin{proof}
First, by Proposition~\ref{prop:saddle-L_2-model},
$A^{-1}(P-\ev)A=\tilde P_0+x^2 E,$ $E\in\tcM^2$.

Let $r_0=\sup\{r';v\in x^{r'}H_\bl^\infty(U)\}.$ Thus,
$v\in x^{r'}H_\bl^\infty(U)$ for all $r'<r,$ hence
$Ev\in x^{r'}H_\bl^\infty(U).$ Thus,
$x^{-1}\tilde P_0 v\in x^{r'+1}H_\bl^\infty(U)$ for all $r'<r_0$.

We can now
apply the previous lemma.
Namely, suppose first that $r_0<0,$ so
$x^{-1}\tilde P_0 v\in x^{r_0+1-\ep}H_\bl^\infty(U)$ for all $\ep>0$.
Part (i) of the lemma, applied with $s=\min(0,r_0+1-\ep)\leq 0,$ shows that
$v\in x^{s-\ep}H_\bl^\infty(U)$ for all $\ep>0.$ Thus, either
$v\in x^{-\ep}H_\bl^\infty(U)$ for all $\ep>0,$ but this contradicts
$r_0<0,$ or
$v\in x^{r_0+1-\ep}H_\bl^\infty(U)$ for all $\ep>0,$ which in turn contradicts
the definition of $r_0.$ Thus, $r_0\geq 0$.

By part (ii) of the lemma, there exists a constant $a_0$
such that $v-a_0\in x^{r'}H_\bl^\infty(U)$ for any $r'<\min(-r_1,1)$.

Let $v'=v-a_0 v_0,$ with $v_0$ given by
Proposition~\ref{prop:saddle-approx-soln}. Then $v'
\in x^{r'}H_\bl^\infty(U)$ for some $r'>0$ still and it satisfies
$x^{-1}A^{-1}(P-\ev)A v'\in \dCI(U)$.
Now let $r'_0=\sup\{r';v'\in x^{r'}H_\bl^\infty(U)\}>0$.
As above, this yields 
$x^{-1}\tilde P_0 v'\in x^{r'+1}H_\bl^\infty(U)$ for all $r'<r'_0$.
Suppose that $r'_0<-r_1.$ Then, by the part (iii) of the
lemma, applied with $s=r'_0+1-\ep,$ $\ep>0,$ $v'\in x^{r'}H_\bl^\infty(U)$
for all $r'<\min(r_1,r'_0+1-\ep),$ hence for some $r'>r'_0,$ since we assumed
$r'_0<-r_1.$ Therefore $v'=v-a_0\in x^{r'}H_\bl^\infty(U)$ for any $r'<-r_1$,
finishing the proof of the first part.

The general case, with $n$ arbitrary, is analogous.
\end{proof}

\begin{thm}\label{thm:saddle-smooth}
Suppose $q \in \Max_+(\ev)$, $O$ is a $W$-balanced neighbourhood of $\Rp$
and $u$ satisfies \eqref{HMV.71r}. Let $u_n=Av_n$
be the local approximate solutions constructed in
Proposition~\ref{prop:saddle-approx-soln}.
Then there exist unique constants $a_n$ such that for any $u'$ with
$u'\sim \sum_{n=0}^\infty a_n u_n,$ $O\cap\WFsc(u-u')=\emptyset$.
That is,
\begin{equation}
u\sim \sum_{n=0}^\infty a_n u_n
\label{eq:saddle-asympt}\end{equation}
microlocally near $\Rp.$ The map $u\mapsto \{a_n\}_{n\in\Nat}$ is a
bijection from microfunctions satisfying \eqref{HMV.71r} to complex-valued
sequences.
Thus, $E_{\mic,\out}(\Rp,\ev)$ is isomorphic to the space of arbitrary
complex-valued sequences, \ie to $\Cx[[x]]$.
\end{thm}

\begin{rem}
The theorem implies that microlocally near $\Rp$,
$u\in x^{r_2/2}L^\infty,$ hence
$u\in x^sL^2_{\scl}(X)$ for all $s<r_2/2 - 1.$ Since $r_2>1,$ this is
an improvement over $x^{-1/2-\ep}L^2_{\scl}(X),$ $\ep>0$.
\end{rem}

\begin{proof}
Let $O_1\ni\Rp$ be any open set with $\overline{O_1}\subset O$,
and choose some
$Q\in\Psisc^{-\infty,0}(X)$ such that $\WFsc'(\Id-Q)\cap \overline{O_1}
=\emptyset$,
$\WFsc'(Q)\subset O.$ Let $\tilde u=Qu.$ Then
\begin{equation}\begin{split}
&\WFsc(\tilde u)\subset
O\cap\WFsc(u)\subset O\cap L_2,\\
&\WFsc((P-\ev)\tilde u)\subset
\WFscp([P,Q])\cap\WFsc(u)\subset (O\setminus \overline{O_1})\cap L_2,
\end{split}
\label{eq:saddle-smooth-16}\end{equation}
hence $\Rp\nin\WFsc((P-\ev)\tilde u)$.

By Corollary~\ref{HMV.79c} and Proposition~\ref{prop:saddle-L_2-model},
$v=A^{-1}\tilde u\in x^rH^\infty_\bl(U')$,
for all $r<(1-r_2)/2=r_1/2<0$,
and by \eqref{eq:saddle-smooth-16}, $A^{-1}(P-\ev)\tilde u$ is
in $\dCI(U')$ for a sufficiently small neighbourhood $U'$ of the critical
point $\Cp\in\pa X$.
Hence we can apply Corollary~\ref{cor:saddle-vf-soln} to
$v=A^{-1}\tilde u$.
We deduce that
there exists a constant $a_0$ such that
$A^{-1}\tilde u-a_0 v_0\in x^{r'}H_\bl^\infty(X)$ for all $r'<-r_1$.
Proceeding iteratively now proves \eqref{eq:saddle-asympt}. Surjectivity
follows by asymptotically summing the $u_n$ given by
Proposition~\ref{prop:saddle-approx-soln}, and injectivity follows from
Corollary~\ref{cor:saddle-vf-soln}. 
\end{proof}

\section{Incoming eigenfunctions at an outgoing saddle}\label{S.In.saddle}

Let $\Rp \in \Max_+(\ev)$ as before, but now consider microlocally incoming
functions $u$ at $\Rp$, \ie, $u \in E_{\mic, -}(\Rp, \ev)$. We have
$\Phi_-(\{q\}) \subset L_1$ locally near $\Rp$, so by microlocalizing 
we may assume that  
\begin{equation}
u\in\CmI(X),\ O\cap\WFsc((P-\ev)u)=\emptyset,\ \WFsc(u)\cap O\subset L_1
\label{HMV.71rrrr}\end{equation}
in some $W$-balanced neighbourhood $O$ of $\Rp$. Thus, in
this section $L_1$ plays the
role that $L_2$ played in the previous section. Essentially all the
results go through as before. In particular, let $\Phi_1/x$ parameterize
$L_1,$ and consider
$\Psisc^{-\infty,0}(X)$-module 
\begin{equation}
\cM\text{ generated by }\Id,\ \{x^{-1}B,\ B\in
\Psisc^{-\infty,0}(X),\ \sigma_{\pa}(B)|_{L_1}=0\}.
\label{eq:cM-saddle-L_1}
\end{equation}
This is a test module, which does {\em not} satisfy the positivity
estimates of \eqref{eq:A_i-comm}-\eqref{eq:A_i-comm-p} -- unless we
take $s>-1/2$, which is of no interest since we will {\em not} have
$u\in x^sL^2_\scl(X)$ for such $s$.

The test module conjugate to $\cM$ under
$e^{-i\Phi_1/x}$
is $\tcM=\Vb(X),$ generated by $xD_x$ and $D_y$.
We then deduce the following analog of
Proposition~\ref{prop:saddle-L_2-model}.

\begin{prop}\label{prop:saddle-L_1-model}
There exists an operator $\tilde A$ having the same structure as in
Proposition~\ref{prop:saddle-L_2-model} such that 
\begin{equation}
\tilde A^{-1}(P-\ev)\tilde A
-\tilde P_0\in x^2\tcM^2,\quad \tilde P_0=-2\tilde\nu
((x^2D_x)+r_2 y(xD_y)).
\end{equation}
\end{prop}

Now $r_2>1,$ so the structure of the Legendre solutions of $\tilde P_0 w=0$ is
different from that considered in the previous section.
In particular, $(yx^{-r_2})^n,$ $n\geq 0,$ integer satisfies
$\tilde P_0(yx^{-r_2})^n=0,$ and $(yx^{-r_2})^n\in x^s H_{\bl}^\infty(X)$ for
all $s<-r_2 n.$ Since $r_2>0,$ these become larger at $\pa X$ as $n\to\infty$.
These can still be modified
to obtain microlocal solutions of $\tilde A^{-1}(P-\ev)
\tilde A v\in\dCI(X),$ with the
proof essentially identical to that presented in the previous section.

\begin{prop}\label{prop:saddle-big-approx-soln}
For each integer $n\geq 0,$ there exists $w_n\in \cap_{s<-r_2 n}
x^s H_{\bl}^\infty(X)$ such that 
$w_n-(yx^{-r_2})^n\in x^{s'} H_{\bl}^\infty(X)$ for some $s'>-r_2 n$,
and such that $\Rp\nin\WFsc(\tilde A^{-1}(P-\ev)\tilde A w_n)$.
\end{prop}

Note that only {\em finite} linear combinations of these are tempered
distributions,
since every $u\in\dist(X)$ lies in $x^s H_{\bl}^m(X)$ for some $s,m$.
In fact, the space of these finite linear combinations is isomorphic to
$E_{\mic,\inc}(\Rp,\ev)$ as we prove in Section~\ref{sec:pairing}.

\section{Microlocal Morse decomposition}\label{sec:Morse}

Equipped with the microlocal eigenspaces $E_{\mic,\pm}(q,\ev)$ defined
below \eqref{HMV.42}, we can decompose $R(\ev+i0)(\dCI(X))$ or more
precisely $\RR_{\ess,+}^\infty(\ev)$ defined in \eqref{HMV.137} 
into microlocal 
components. Using the identification \eqref{HMV.60} of the space $\SEEF(\ev)$
of `smooth eigenfunctions' (see \eqref{HMV.10}) with
$\RR_{\ess,+}^\infty(\ev),$ this decomposes 
elements of $\SEEF(\ev)$ into microlocal building blocks, and in particular
proves Theorem~\ref{HMV.52}.

Let us start with the perfect Morse case, with $\Cp_{\min}$ and $\Cp_{\max}$
denoting the minimum and maximum of $\Vy.$ Denote by $\iota$
the identification map
\begin{equation}
\iota:E_{\mic,\out}(\Rp_{\min}^+,\ev)\longrightarrow \RR_{\ess,+}^\infty(\ev).
\label{eq:iota}\end{equation}
Namely, by Corollary~\ref{cor:min-ident-8}, $u\in
E_{\mic,\out}(\Rp_{\min}^+,\ev),$ may be regarded as a distribution with
$\WFsc(u)\subset\{\qminp\}$ and $f=(P-\ev)u\in\dCI(X),$ determined up
to addition of an element of $\dCI(X).$ Then $u=R(\ev+i0)f$ since both
sides are outgoing, and their difference is a generalized eigenfunction of
$P,$ so $u$ can indeed be identified with an element of
$\RR_{\ess,+}^\infty(\ev).$ The map
\begin{equation*}
r:\RR_{\ess,+}^\infty(\ev)\longrightarrow E_{\mic,\out}(\qmaxp, \ev)
\end{equation*}
arises by microlocal restriction. That is, we choose $Q \in
\Psisc^{0,0}(X)$ with $\sigma(Q) = 1$ in a neighbourhood of $\qmaxp$, and
map $u \in \RR_{\ess,+}(\ev)$ to $Qu \in E_{\mic,+}(\qmaxp)$. 

\begin{thm}[Theorem~\ref{HMV.52}] Suppose that $\dim X=2,$ $\pa X$ is
connected and $\Vy$ is perfect Morse. For $\kappa<\ev<K,$ $\iota$ is an
isomorphism. If $K<\ev$ then $r$ is surjective and has null
space which restricts isomorphically onto $E_{\mic,\out}(\qminp,\ev)$
leading to the short exact sequence
\begin{equation*}
0 \longrightarrow E_{\mic,\out}(\qminp,\ev)\longrightarrow
R^\infty_{\ess,+}(\ev) \longrightarrow E_{\mic,\out}(\qmaxp,\ev)
\longrightarrow 0.
\end{equation*}
\end{thm}

\begin{proof}
The injectivity of $\iota$ is clear.

If $\ev < K$, then any $u \in \RR_{\ess,+}^\infty(\ev)$ has scattering
wavefront set in $\{ \nu \geq 0 \} \cap \Sigma(\ev)$, and is invariant
under bicharacteristic flow. Thus, we must have $\WFsc(u) \subset \{ \qminp
\}$. If $Q$ is in
$\Psisc^{0,0}(X)$ with $\sigma(Q) = 1$ in a neighbourhood of $\qminp$, 
then $Qu \in E_{\mic,+}(\qminp,\ev)$ and $\iota(Qu) = u$, which shows that
$\iota$ is an isomorphism in this range. 

Next suppose that $\ev > K$. Then $r\circ\iota=0$ since $\qmaxp\nin\WFsc
(\iota(u))$ for $u\in E_{\mic,\out}(\qminp,\ev).$ Conversely, if
$v\in \RR_{\ess,+}^\infty(\ev)$ maps to $0$ under $r$,
then by definition $\nu\geq\nu(\qmaxp)$ on $\WFsc(v)$,
and $\qmaxp\nin\WFsc(v),$ hence by H\"ormander's theorem on
the propagation of singularities, $\WFsc(v)
\subset\{\qminp\}.$ Thus, $v$ is a representative of an element
of $E_{\mic,\out}(\qminp,\ev),$ and is indeed in the range of $\iota.$
This shows exactness at $\RR_{\ess,+}^\infty(\ev)$.
Finally, surjectivity of $r$ can be seen as follows.
Any element of $E_{\mic,\out}(\qmaxp,\ev)$ has a representative $u$
as in Corollary~\ref{cor:min-ident-8}. In particular,
$f=(P-\ev)u\in\dCI(X),$ $\WFsc(u)\subset\{\nu\geq\nu(\qmaxp)\}$,
so $u=R(\ev+i0)f,$ showing surjectivity of $r$.
\end{proof}

The general situation, \ie when $\dim X=2$
and $\Vy$ is Morse, is not substantially different if none of the
bicharacteristic curves $L_+=L_+(q)$, emanating from $q \in \Max_+(\ev)$, hits
another $q' \in \Max_+(\ev)$. 
For each $q \in \Min_+(\ev)$ there is an identification map
\begin{equation}
\iota_{q}:E_{\mic,\out}(q,\ev)\to \RR_{\ess,+}^\infty(\ev)
\end{equation}
defined as described after \eqref{eq:iota}. 
Similarly, for each $q\in \Max_+(\ev),$ there is a restriction map
\begin{equation*}
R_{q}:\RR_{\ess,+}^\infty(\ev)\to E_{\mic,\out}(q)
\end{equation*}
as described above. Then similar arguments to above yield

\begin{prop}
Suppose that $\dim X=2,$ $\Vy$ is Morse and none of the curves $L_+(q)$, $q
\in \Max_+(\ev)$, hit another $q' \in \Max_+(\ev)$. Then
\begin{equation*}
0 \longrightarrow \!\!\!\! \bigoplus_{q\in\Min_+(\ev)} \!\!\!\!
E_{\mic,\out}(q,\ev)
\longrightarrow R^\infty_{\ess,+}(\ev)
\longrightarrow \!\!\!\! \bigoplus_{q'\in\Max_+(\ev)} \!\!\!\!
E_{\mic,\out}(q',\ev) \longrightarrow 0
\end{equation*}
is a short exact sequence.
\end{prop}
\noindent Note that in case $\Vy$ is perfect Morse, we simply recover
Theorem~\ref{HMV.52}.

In complete generality, \ie for general Morse $\Vy,$ the relationship
between critical points is more complicated. We introduce a partial
order on $\RP_+(\ev)$ corresponding to the flow-out under $W$.

\begin{Def}\label{Def:partial-order} 
If $\Rp,$ $\Rp'\in\RP_+(\ev)$ we say that $\Rp\leq \Rp'$ if
$\Rp'\in\Phi_+(\{\Rp\})$ and $\Rp<\Rp'$ if $\Rp\leq\Rp'$ but $\Rp'\neq\Rp.$
A subset $\Gamma\subset\RP_+(\ev)$ is closed under $\leq$ if for all
$\Rp\in\Gamma,$ we have $\{ q' \in \RP_+(\ev) \mid q \leq q' \} \subset
\Gamma$. We call the set $\{ q' \in \RP_+(\ev) \mid q \leq q' \}$ the
string generated by $q$. 
\end{Def}

This partial order relation between two radial points in $\RP_+(\ev)$
corresponds to the existence of a sequence $\Rp_j\in\RP_+(\ev),$
$j=0,\ldots,k,$ $k\geq 1$, with $\Rp_0=\Rp,$ $\Rp_k=\Rp'$ and such that for
every $j=0,\ldots,k-1$, there is a bicharacteristic $\gamma_j$
with $\lim_{t\to-\infty}\gamma_j=\Rp_j$ and
$\lim_{t\to+\infty}\gamma_j=\Rp_{j+1}.$

We can now define a subspace of $\RR_{\ess,+}^\infty(\ev)$ associated to a
$\leq$-closed subset $\Gamma\subset\RP_+(\ev)$ by setting  
\begin{multline}
\RR_{\ess,+}^\infty(\ev,\Gamma)\\
=\{u\in\dist(X);(P-\ev)u\in\dCI(X),\ \WFsc(u)\cap\RP_+(\ev)
\subset\Gamma\}/\dCI(X).
\label{HMV.133}\end{multline}
Notice that the assumption that $\Gamma\subset\RP_+(\ev)$ implies that
$\WFsc(u)\subset\{\nu>0\}.$ The `trivial' case where $\Gamma =\RP_+(\ev)$
is the space we are ultimately trying to describe:
\begin{equation}
\RR_{\ess,+}^\infty(\ev,\RP_+(\ev))\equiv \RR_{\ess,+}^\infty(\ev),
\ \ev\notin\Cv(V).
\label{HMV.134}\end{equation}

\begin{prop}\label{prop:Gamma}
Suppose that $\Gamma\subset\RP_+(\ev)$ is $\leq$-closed and $\Rp$ is a
$\leq$-minimal element of $\Gamma$. Then with $\Gamma '=\Gamma
\setminus\{\Rp\}$ 
\begin{equation*}\xymatrix{ 0 \ar[r]^{}&
\RR_{\ess,+}^\infty(\ev,\Gamma')\ar[r]^{\iota}&
\RR_{\ess,+}^\infty(\ev,\Gamma)\ar[r]^{r_\Rp}&
E_{\mic,\out}(\ev,\Rp) \ar[r]^{}& 
0}
\end{equation*}
is a short exact sequence.
\end{prop}

\begin{proof} The injectivity of $\iota$ follows from the definitions. The
null space of the microlocal restriction map $r_{\Rp},$ which can be viewed
as restriction to a $W$-balanced neighbourhood of $\Rp,$ is precisely the
subset of $\RR_{\ess,+}^\infty(\ev,\Gamma)$ with wave front set disjoint
from $\{\Rp\},$ and this subset is $\RR_{\ess,+}^\infty(\ev,\Gamma').$ Thus
it only remains to check the surjectivity of $r_{\Rp}.$

We do so first for the strings generated by $q \in \RP_+(\ev)$. For $q \in
\Min_+(\ev)$, the string just consists of $q$ itself and the result follows
trivially. So consider the string $S(q)$ generated by $q \in \Max_+(\ev)$. 
By Corollary~\ref{cor:min-ident-8} any element of
$E_{\mic,+}(q,\ev)$ has a representative $\tilde u$ satisfying
$(P-\ev)\tilde u\in\dCI(X)$ with $\WFsc(\tilde
u)\subset\Phi_+( \{ \Rp \} )$, which immediately gives surjectivity in
this case. 

For any $\leq$-closed set $\Gamma$ and $\leq$-minimal element $\Rp$,
the string $S(q)$ is contained in $\Gamma$, so the surjectivity of $r_{\Rp}$
follows in general.
\end{proof}

Notice that we can always find a sequence 
$\emptyset=\Gamma_0\subset\Gamma_1\subset\ldots\subset\Gamma_n=\RP_+(\ev),$
of $\leq$-closed sets with
$\Gamma_j\setminus\Gamma_{j-1}$ consisting of a single point $\Rp_j$
which is $\leq$-minimal in $\Gamma_j$: we simply order the $q_i \in
\RP_+(\ev)$ so that $\nu(q_1) \geq \nu(q_2) \geq \dots$, and set $\Gamma_i =
\{ q_1, \dots, q_i \}$. Then Proposition~\ref{prop:Gamma} implies the
following

\begin{thm}
Suppose that
$\emptyset=\Gamma_0\subset\Gamma_1\subset\ldots\subset\Gamma_n=\RP_+(\ev),$
is as described in the previous paragraph. Then
\begin{equation}
\{0\}\longrightarrow R^\infty_{\ess,+}(\ev,\Gamma_1)\hookrightarrow \ldots
\hookrightarrow
R^\infty_{\ess,+}(\ev,\Gamma_{n-1})\hookrightarrow R^\infty_{\ess,+}(\ev),
\end{equation}
with
\begin{equation}
R^\infty_{\ess,+}(\ev,\Gamma_{j})/R^\infty_{\ess,+}(\ev,\Gamma_{j-1})
=E_{\ess,+}(\Rp_j,\ev),\ j=1,2,\ldots,n.
\end{equation}
\end{thm}

This theorem shows that there are elements of $E^{\infty}_{\ess,+}(\ev)$
associated to each radial point $\Rp.$ The parameterization of
the microlocal eigenspaces shows that in some sense, there are `more'
eigenfunctions associated to the minima than to the maxima. One way
to make this precise is to show that the former class is dense
in $E^{\infty}_{\ess}(\ev),$ in the topology of $\dist(X),$ or indeed
in the topology of $E^s_{\ess}(\ev)$ for $s>0$ small.
We define
\begin{equation}
E^{\infty}_{\Min,+}(\ev)=\left\{u\in E^{\infty}_{\ess}(\ev);
\WFsc(u)\cap\RP_+(\ev)\subset\Min_+(\ev)\right\}.
\label{HMV.152p}\end{equation}

\begin{prop}\label{HMV.164} Let $s_0=-1/2+\min
\{r_2(\Rp)/2;\ \Rp\in\Max_+(\ev)\}>0$. Then for any $\ev\notin\Cv(C)$ and
$0<s<s_0,$ $E^{\infty}_{\Min,+}(\ev)$ is dense in $E^{\infty}_{\ess}(\ev)$ in
the topology of $E^{s}_{\ess}(\ev)$ given by
Proposition~\ref{prop:topologies-s}. 
\end{prop}

\begin{proof} We first suppose that 
\begin{equation}
\ev \text{ is not in the point spectrum of } P.
\label{eq:spec-ass}\end{equation}
Let $u\in E^{\infty}_{\ess}(\ev).$ We know from
Theorem~\ref{thm:saddle-smooth}
(see the remark following it) that
\begin{equation*}
\WFsc^{*,s-1/2}(u)\subset \Min_-(\ev)\cup\Min_+(\ev),\ s<s_0,
\end{equation*}
\ie the wave front set of $u$ relative to $x^{s-1/2}\Hsc^\infty(X)$ is
localized at the radial points over the minima, provided that $s<s_0$; in
particular we can take $s>0$.

Now choose $Q\in\Psisc^{-\infty,0}(X)$ which has $\WFscp(Q)$ localized near
$\Min_+(\ev)$ and $\WFscp(\Id-Q)$ disjoint from $\Min_+(\ev).$
Then
\begin{equation*}
\WFsc^{*,s-1/2}(u-Qu)\subset \Min_-(\ev),\ s<s_0.
\end{equation*}

Let $\Pt=\Lap+\Vt$
be the operator constructed in Lemma~\ref{lemma:Vt} with $\tilde\nu$
large, so that $\ev-\tilde\nu^2<\inf V.$ Then every radial
point $\Rp\in\widetilde{\RP_+}(\ev)$
of $\Pt$ lies above a minimum of $\Vt$ (\ie $\widetilde{\RP_+}(\ev)
=\widetilde{\Min_+}(\ev)$). Let $\tilde\Pi$ be the orthogonal projection
off the $L^2$-nullspace of $\Pt-ev$, so $\Id-\tilde\Pi$ is a finite
rank projection to a subspace of $\dCI(X)$.
Since $\WFscp([P,Q])\subset\WFscp(Q)\cap\WFscp(\Id-Q)$, we have 
\begin{equation*}
\WFsc^{*,s-1/2}(u)\cap \WFscp([P,Q])=\emptyset,\ s<s_0.
\end{equation*}
Moreover,
\begin{equation*}
\WFscp(P-\tilde P)\cap\WFscp(Q)=\emptyset,
\end{equation*}
and $(P-\ev)Qu=[P,Q]u+Q(P-\ev)u,$ so
\begin{equation*}
(P-\ev)Qu,(\tilde P-\ev)Qu\in x^{s+1/2}\Hsc^\infty(X),\ s<s_0.
\end{equation*}
Here we can take $s>0$, hence we can apply $R(\ev-i0)$ to these functions, and
\begin{equation*}
\WFsc^{*,s-1/2}(R(\ev-i0)(P-\ev)Qu)\subset \Min_-(\ev).
\end{equation*}

We claim that
\begin{equation}
u=v-R(\ev-i0)(P-\ev)v,\ v=Qu.
\label{HMV.168}\end{equation}
To see this, let $\tilde u=u-(v-R(\ev-i0)(P-\ev)v).$ For $s<s_0$,
$\WFsc^{*,s-1/2}(\tilde u)$ lies in $\Min_-(\ev),$ since
this is true both for $u-v$ and for $R(\ev-i0)(P-\ev)v.$ Moreover,
$(P-\ev)\tilde u=0,$ so the uniqueness theorem, Theorem~\ref{thm:unique},
shows that $\tilde u=0$.

Since $v=Qu$ has wavefront set confined to
$\{\nu>0\}$ and hence is outgoing,
\begin{equation*}
\tilde\Pi v=\tilde R(\ev+i0)(\tilde P-\ev)v.
\label{HMV.167}\end{equation*}
Thus,
\begin{equation*}\begin{split}
u=(\Id-\tilde\Pi) Qu+\tilde R(\ev+i0)w&-R(\ev-i0)(P-\ev)\tilde R(\ev+i0)w,\\
 &\Mwhere
w=(\tilde P-\ev)Qu.
\label{HMV.169}\end{split}\end{equation*}

Now consider $u'_j=\phi(x/r_j) u,$ where $\phi\in\Cinf(\bbR)$ is identically
$1$ on $[2,+\infty),$ $0$ on $[0,1],$ and $r_j>0$ satisfy
$\lim_{j\to\infty} r_j=0.$ Thus, $u'_j\in\dCI(X),$ and 
\begin{equation}\begin{aligned}
{}&u_j=(\Id-\tilde\Pi) Q u'_j+\tilde R(\ev+i0)w_j-R(\ev-i0)(P-\ev)\tilde R(\ev+i0)w_j,\\
{}&w_j=(\tilde P-\ev)Qu'_j=(\tilde P-P)Qu'_j+[P,Q\phi(x/r_j)]u
\end{aligned}\label{HMV.171}\end{equation}
satisfy $(P-\ev)u_j=0$ and $\WFsc(u_j)\cap\{\nu\geq 0\}\subset \Min_+(\ev),$
so $u_j\in E^{\infty}_{\Min,+}(\ev)\subset\SEEF(\ev).$ 
Since $w_j\to w$ and $(P-\ev)Qu'_j \to (P-\ev)Qu$ in
$x^{s+1/2}L^2_{\scl}(X)$,
we have $u_j \to u$ in $x^{r}L^2_{\scl}(X)$, for any $r < -1/2$ by
Theorem~\ref{thm:lim-abs}. Moreover, if $B$ is as in \eqref{eq:B} then
since $\WFsc'(B) \cap \RP(\ev) = \emptyset$ and $\WFsc(w) \cap \RP(\ev) =
\emptyset$, we have $Bu_j \to Bu$ in $x^{s-1/2}L^2_{\scl}(X)$ by
Theorem~\ref{thm:prop-sing}, proving
convergence in $E^s_{\ess}(\ev)$.

We now indicate the modifications necessary if \eqref{eq:spec-ass} is not
satisfied. Let $\Pi$ denote the orthogonal projection off the $L^2$
$\ev$-eigenspace. Equation \eqref{HMV.168} must be replaced by 
$u = \Pi v - R(\ev-i0)(P-\ev)v$. 
Then if we define $u_j$ by
\begin{equation*}
u_j= \Pi \tilde R(\ev+i0)w_j-R(\ev-i0)(P-\ev)\tilde R(\ev+i0)w_j,
\end{equation*}
instead of \eqref{HMV.171}, the argument goes through. 
\end{proof}

\section{Pairing and duality}\label{sec:pairing}

We define a basic version of `Green's pairing' in this context on 
\begin{equation}
\Et^{-\infty}(\ev)=\left\{u\in\CmI(X);(P-\ev)u\in\dCI(X)\right\}.
\label{HMV.142}\end{equation}
Namely 
\begin{equation}
\Bt(u_1,u_2)=i\left(
\langle (P-\ev)u_1,u_2\rangle-\langle u_1,(P-\ev)u_2\rangle\right),\
u_i\in\Et^{-\infty}(\ev).
\label{HMV.145}\end{equation}
The imaginary factor is inserted only to make the pairing sesquilinear, 
\begin{equation*}
\overline{\Bt(u_1,u_2)}=\Bt(u_2,u_1).
\label{HMV.143}\end{equation*}
Observe that 
\begin{equation*}
\Bt(u,v)=0\Mif u\in\dCI(X)
\label{HMV.151}\end{equation*}
since then integration by parts is permitted.

More significantly we define a sesquilinear form on $E^{\infty}_{\ess}(\ev)$
using $\Bt.$ Assuming that $\ev\notin\Cv(V),$ choose 
\begin{multline}
A\in\Psisc^{0,0}(X)\Msuchthat
\WFscp(A)\cap\Sigma(\ev)\subset \{ \nu > -a(\ev) \}  \Mand\\
\WFscp(\Id-A)\cap\Sigma(\ev)\subset \{ \nu < a(\ev) \}=\emptyset
\label{HMV.147}\end{multline}
where $a(\ev)$ is as in Remark~\ref{rem:Phi(RP)}. 
If $u\in E^{\infty}_{\ess}(\ev)$ then $\WFsc((P-\ev)Au)=\emptyset$ so
$Au\in\Et^{-\infty}(\ev).$ This allows us to define
\begin{equation}
B:E^{\infty}_{\ess}(\ev)\times E^{\infty}_{\ess}(\ev)\longrightarrow \bbC,\
B(u_1,u_2)=\Bt(Au_1,Au_2).
\label{HMV.148}\end{equation}

\begin{lemma}\label{HMV.149} If $\ev\notin\Cv(V)$ and $A$ satisfies
\eqref{HMV.147} then $B$ in \eqref{HMV.148} is independent of the choice of
$A$ and defines a non-degenerate sesquilinear form.
\end{lemma}

\begin{proof} The space of cut-off operators satisfying \eqref{HMV.147} is
convex. If we consider the linear homotopy between two of them and denote
the bilinear form $B_t$ for $t\in[0,1]$ then 
\begin{equation*}
\frac d{dt}B_t(u_1,u_2)=\Bt(\frac d{dt}A_tu_1,A_tu_2)+ \Bt(A_tu_1,\frac
d{dt}A_tu_2)=0
\label{HMV.150}\end{equation*}
since $\frac{d}{dt}A_tu_i\in\dCI(X)$ for $i=1,2.$

The non-degeneracy of $B$ follows by rewriting it 
\begin{equation}
B(u_1,u_2)=i\langle (P-\ev)Au_1,u_2\rangle,\ u_1,\ u_2\in\SEEF(\ev).
\label{HMV.160}\end{equation}
Indeed the difference between the right side of \eqref{HMV.160} and
\eqref{HMV.148} is  
\begin{equation*}
i\left(\langle(P-\ev)Au_1,(\Id-A)u_2\rangle
+\langle Au_1,(P-\ev)Au_2\rangle\right).
\label{HMV.161}\end{equation*}
In the first term the scattering wavefront sets of $Au_1$ and $(\Id-A)u_2$
do not meet, so integration by parts is permitted and, $u_2$ being an
eigenfunction, this difference vanishes. If \eqref{HMV.160} vanishes for
all $u_1$ it follows that $\langle f,u_2\rangle =0$ for all $f\in\dCI(X)$
since it vanishes for $f\in(P-\ev)\dCI(X) \oplus E_{\pp}(\ev)$ which, by
Propositions~\ref{prop:Sp} and \ref{HMV.192} is a complement to
$(P-\ev)A\SEEF(\ev)$ in $\dCI(X).$ This proves the non-degeneracy.
\end{proof}

\begin{lemma}\label{HMV.163} The pairing $B$ extends to a non-degenerate
pairing and topological duality
\begin{equation}
B:\SEEF(\ev)\times \EEF(\ev)\longrightarrow \bbC,\ B(u,v)=\Bt(Au, Av).
\label{HMV.162}\end{equation}
Moreover, for all $s\in\bbR,$
$B$ extends to a nondegenerate pairing
\begin{equation}
B:\FEEF s(\ev)\times \FEEF{-s}(\ev)\longrightarrow \bbC,
\ B(u,v)=\Bt(Au, Av).
\label{HMV.162p}\end{equation}
\end{lemma}

\begin{proof} The extension of the pairing to \eqref{HMV.162} follows
directly from \eqref{HMV.160}, since $Au_1\in\dCI(X)$ when
$u_1\in\SEEF(\ev).$ The same argument shows that it is non-degenerate. To
see that it is a topological pairing observe that the Fr\'echet topology on
$\SEEF(\ev)$ arises from its identification with the quotient of
$\dCI(X)$ by the closed subspace $(\Ham-\ev)\dCI(X).$
Thus an element of the dual space may be identified with a distribution
$u\in\CmI(X)$ which vanishes on $(\Ham-\ev)\dCI(X),$ \ie is an element of
$\EEF(\ev).$ The pairing between $\EEF(\ev)$ and $\SEEF(\ev)$ defined in
this way is precisely $B.$

A similar argument applies to the spaces of eigenfunctions of finite
regularity. The square of the Hilbert norm on $\FEEF{s}(\ev)$ may be taken
to be
\begin{equation}
\|Au\|^2_{x^{s-1/2}L^2}+\|u\|^2_{x^{-N}L^2}
\label{HMV.199}\end{equation}
where $A$ is as in \eqref{HMV.147} and $N$ is sufficiently large, in
particular $N>\max\{-s+\frac12,\frac12\}.$ Writing the pairing in the form 
\begin{equation}
B(u_1,u_2)=i\langle [P,A]u_1,u_2\rangle 
\label{HMV.200}\end{equation}
shows its continuity with respect to \eqref{HMV.199} (for a different $A$
which is microlocally the identity on the essential support of the operator
in \eqref{HMV.200}.) The non-degeneracy follows much as before.
\end{proof}

As well as these global forms of the pairing, we can define versions of it on the
components of the microlocal Morse decomposition. Namely, for
$\Rp\in\RP_+(\ev),$ choose $Q\in\Psisc^{0,0}(X)$ with $\Rp\notin\WFscp(\Id-Q)$ and
$\WFscp(Q)\cap\Sigma(\ev)$ supported in a $W$-balanced neighbourhood $O$ of
$\Rp.$ Then consider 
\begin{equation}
\Et_{\mic,+}(\Rp,\ev)\times\Et_{\mic,-}(\Rp,\ev)\ni
(u,v)\longmapsto \Bt(Qu,Qv).
\label{HMV.174}\end{equation}
The pairings on the right are well defined. Indeed, for such
$(u,v)$,
\begin{equation}
\WFsc(u)\cap\WFsc(v)\cap O\subset\{\Rp\},
\end{equation}
since $\WFsc(u)\cap O\subset\Phi_+(\{\Rp\})$ and
$\WFsc(v)\cap O\subset\Phi_-(\{\Rp\}).$
Moreover, $\Rp\nin\WFsc\left((P-\ev)Qu\right)$
and $\Rp\nin\WFsc\left((P-\ev)Qv\right),$ so
\begin{equation}
\WFsc\left((P-\ev)Qu\right)\cap\WFsc(Qv)=
\WFsc\left((P-\ev)Qv\right)\cap\WFsc(Qu)=\emptyset,
\label{HMV.175}\end{equation}
so $B$ in \eqref{HMV.162} is well-defined.

If $Q'$ is any other microlocal cutoff with the same properties then
$\Rp\nin\WFscp(Q-Q'),$ hence $\WFsc((Q-Q')u)\cap\WFsc(Qv)=\emptyset$,
so $\tilde B(Qu-Q'u,Qv)=0.$ Similarly, $\tilde B(Qu,Qv-Q'v)=0$,
so $\tilde B(Qu,Qv)=\tilde B(Q'u,Q'v)$.

If we define a new pairing $B'$ as in \eqref{HMV.162}, but using the
operator $\tilde P$ from Lemma~\ref{lemma:Vt} in place of $P$, then
$B'(Qu,Qv)=\Bt(Qu,Qv)$ provided that $O\cap\WFscp(P-\tilde P)=\emptyset$.

\begin{prop}\label{HMV.172} For each $\Rp\in\RP_+(\ev)$ the sesquilinear
form \eqref{HMV.174} descends to a non-degenerate pairing
\begin{equation}
\Bt_{\Rp}:E_{\mic,+}(\Rp,\ev)\times
E_{\mic,-}(\Rp,\ev)\longrightarrow \bbC.
\label{HMV.173}\end{equation}
\end{prop}

\begin{proof} The terms involved in the quotients all have wavefront
set disjoint from $\WFscp(Q),$ so integration by parts shows that
\eqref{HMV.173} is independent of choices.

To show non-degeneracy,
suppose that $\Rp$ is a minimal element of $\RP_+(\ev)$ with respect
to the partial order $\leq$ of Definition~\ref{Def:partial-order}. 
Then every
element $v$ of $\Et_{\mic,-}(\Rp,\ev)$ has a representative
in $\EEF(\ev)$ with $\WFsc(v)\subset\{\nu\leq\nu(\Rp)\}$ and
$\WFsc(v)\cap \{\nu=\nu(\Rp)\}=\{\Rp\}.$ By
Lemma~\ref{HMV.163}, there exists $u\in\SEEF(\ev)$ such that $\Bt(u,v)
\neq 0.$ Since $\Rp$ is $\leq$-minimal, we have
$\WFsc(u)\cap\WFsc(v)\subset\{\Rp\}$. Hence, $\Bt(u,v)=
\Bt(Qu,Qv).$ Thus, $Qu$ is an element of $E_{\mic,+}(\Rp,\ev)$
with $\Bt(Qu,v) \neq 0$. 

If $\Rp$ is not $\leq$-minimal, then we use Lemma~\ref{lemma:Vt} with
$\tilde\nu$ smaller than $\nu(\Rp)$, but larger than $\nu(\Rp')$ for any $\Rp' \in
\RP_+(\ev)$ with $\nu(\Rp')$ smaller than $\nu(\Rp)$. Then $\Rp \in
\widetilde{\RP_+(\ev)}$, and $\tilde V = V$ in a neighbourhood of
$\pi(\Rp)$, so $E_{\mic,\pm}(\Rp,\ev)$ is the same space irrespective
of whether we consider $P$ or $\tilde P$. However, 
$\Rp$ is $\leq$-minimal for $\tilde P$, so the result follows from the
previous paragraph. 
\end{proof}

\begin{prop}
Suppose that $\Rp \in \Max_+(\ev)$ and $O$ is a $W$-balanced
neighbourhood of $\Rp$ in $\scT^*_YX.$ Let $u_n=Av_n\in E_{\mic,+}(\Rp,\ev)$
and $\tilde u_n=\tilde Aw_n\in E_{\mic,+}(\Rp,\ev)$
be the local approximate outgoing, resp\@. incoming, solutions constructed in
Proposition~\ref{prop:saddle-approx-soln},
resp\@. Proposition~\ref{prop:saddle-big-approx-soln}. Then
\begin{equation*}\begin{split}
&\Bt_{\Rp}(u_n,\tilde u_m)=0\ \Mif n>m,\\
&\Bt_{\Rp}(u_n,\tilde u_n)\neq 0\ \forall\ n.
\end{split}\end{equation*}
\end{prop}

\begin{proof}
This follows from a straightforward calculation, using the asymptotic
expansions, similarly to the proof of Proposition~\ref{HMV.153} below.
\end{proof}

We can thus renormalize $\tilde u_m$ inductively by letting
\begin{equation*}
\tilde u_m'=\frac{1}{\Bt_{\Rp}(u_m,\tilde u_m)}
\left(\tilde u_m-\sum_{n<m}\Bt_{\Rp}(u_n,\tilde u_m)
\,\tilde u'_n\right).
\end{equation*}
Then $\Bt_{\Rp}(u_n,\tilde u'_m)=0$ if $n\neq m$ and
$\Bt_{\Rp}(u_m,\tilde u'_m)=1$.

An immediate consequence of non-degeneracy is the following theorem.

\begin{thm}\label{thm:saddle-big}
Suppose that $\Rp \in \Max_+(\ev)$ and $O$ is a
$W$-balanced neighbourhood of $\Rp$
in $\scT^*_YX.$ If $u\in E_{\mic,-}(q,\ev)$ is a microlocally incoming
eigenfunction at $q$,
then there exist unique constants $a_n$ such that $a_n\neq 0$ for
only finitely many $n,$ and
\begin{equation}
O\cap\WFsc(u-\sum a_n \tilde u_n)=\emptyset.
\end{equation}
\end{thm}
\noindent 
Thus, $E_{\mic,\inc}(\Rp,\ev)$ is isomorphic to the space of finite
complex-valued sequences, $\Cx[x]$.

\begin{proof}
Choose $s>0$ such that
$u\in W'_s=E_{\mic,-}(\Rp,\ev)\cap x^{-s-1/2}L^2(X)$.
Since $W'_s$ automatically annihilates $E_{\mic,+}(\Rp,\ev)\cap
x^{s-1/2} L^2(X),$ we deduce that $\Bt_{\Rp}(u_n,u)=0$ for $n$ sufficiently
large. Now
$u'=u-\sum_n \Bt_{\Rp}(u_n,u) \tilde u'_n$  satisfies $\Bt_{\Rp}(u_n,u')=0$
for all $n.$ The non-degeneracy of $\Bt_{\Rp}$ and the structure theorem
for $E_{\mic,+}(\Rp,\ev),$ Theorem~\ref{thm:saddle-smooth},
thus imply that $u'=0$.
Since the $\tilde u'_n$ are finite linear combinations of $\tilde u_n$,
the proof is complete.
\end{proof}

We can compute the form of $B$ explicitly on the subspace
$E^\infty_{\Min,+}(\ev)$. 
Recall that for each $\Rp\in\Min_+(\ev)$ there is a map 
\begin{equation*}
M_+(\Rp,\ev):E^{\infty}_{\Min,+}(\ev)\longrightarrow \cS(\bbR_{r(\Rp)})
\label{HMV.156}\end{equation*}
given by \eqref{HMV.48} or \eqref{HMV.50}, which is an isomorphism on
$E_{\mic,+}(\Rp)$. 

\begin{prop}\label{HMV.153} If $\ev\notin\Cv(V)$ then for
$u_1,u_2\in E^{\infty}_{\Min,+}(\ev),$
\begin{equation}
B(u_1,u_2)=2\sum\limits_{\Rp\in\Min_+(\ev)}
\sqrt{\ev-V(\pi(\Rp))}\int_{\bbR_{r(\Rp)}} M_+(\Rp,\ev)(u_1)
\overline{M_+(\Rp,\ev)(u_2)}\,\omega_{\Rp},
\label{HMV.154}\end{equation}
where $\omega_q$ is the density on the front face of the blow-up induced
by the Riemannian density of $h.$
\end{prop}

\begin{proof}
Let $\chi\in\Cinf(\Real)$ be identically $0$ near $0$ and identically
$1$ on $[1,+\infty),$ and let $\chi_r(x)=\chi(x/r).$ Then, with $A$
as in \eqref{HMV.147},
\begin{equation}\begin{split}
B(u_1,u_2)&=i\lim_{r\to 0} (\langle \chi_r (P-\ev)Au_1,Au_2\rangle
-\langle \chi_rAu_1,(P-\ev)Au_2\rangle)\\
&=-i\lim_{r\to 0} \langle [P,\chi_r] Au_1,Au_2\rangle=-i\lim_{r\to 0}
\int_X ([P,\chi_r] Au_1)
\overline{Au_2}\,dg,
\end{split}
\label{eq:B-pairing-8}\end{equation}
since for $r>0$,
$$
\langle (P-\ev)\chi_r Au_1,Au_2\rangle =\langle \chi_r Au_1,(P-\ev)Au_2\rangle.
$$
Since $[P,\chi_r]\to 0$ strongly as a map $\Hsc^{1,s-1/2}(X)\longrightarrow 
\Hsc^{0,s+1/2}(X),$ for $s\in\Real,$
it follows that only microlocal regions where at least one of $Au_1$ and
$Au_2$ is not in $\Hsc^{0,-1/2}(X)$ contribute. That is, we can insert
microlocal cut-offs near the minimal radial points and thereby localize the
computation to a single radial point. 

So, let $q \in \Min_+(\ev)$, let $z = \pi(q)$, and 
assume that $\ev$ is above the Hessian threshold for $q$. 
Using the asymptotic form of $Au_1$ and $Au_2$ given by
Theorem~\ref{thm:sink-smooth}, we need to compute
$$
-i\lim_{r\to 0}\int [\Delta,\chi_r](e^{i\Phi_2/x}x^{\beta}v_1)
e^{-i\Phi_2/x}x^{\overline{\beta}}\overline{v_2}\,\frac{dx}{x^2}
\,\frac{\omega}{x^{r_2}},
$$
where $\omega$ is the density on the front face of the blow-up induced
by the Riemannian density of $h.$ Thus, $\omega=|dh(\Cp)|/x^{r_1}$.
Here, $|dx|$ is fixed at the boundary by the requirement that the
scattering metric takes the form \eqref{HMV11}, so 
this is a well-defined density on the front face of $X_{z,r_1}$. 
In the computation the commutator $[\Delta,\chi_r]$
may be replaced by $2(x^2D_x\chi_r)(x^2D_x)$ modulo terms that vanish
in the limit, and $v_1,$ $v_2$ may be restricted to the front face.
In addition, $x^2D_x$ must fall on $e^{i\Phi_2/x},$ or the limit vanishes.
Finally, $\re \beta = r_2/2$. Hence, the limit is
\begin{equation*}
2i\Phi_2(\Cp)(\int v_1|_{\ff}\overline{v_2|_{\ff}}\,\omega)\,
(\lim_{r\to 0}\int_0^\infty (D_x\chi_r)\,dx)
=2\sqrt{\ev-V(\Cp)}(\int v_1|_{\ff}\overline{v_2|_{\ff}}\,\omega).
\label{HMV.155}\end{equation*}
Since $v_j|_{\ff}=M_+(q,\ev)(u_j),$ this gives the stated result if $\ev$ is
above the Hessian threshold of $\Cp$.

Now suppose that the $\ev$ is below the Hessian threshold of $q$. 
Then the asymptotics of $Au_i$ are given by Theorem~\ref{thm:center-smooth} 
and take the form
\begin{equation*}
e^{i\tilde\nu/x}\sum\limits_{j} X^{1/2}X^{i((2j+1)\alpha + \gamma)}
\gamma_{j} v_j(Y), \ \tilde\nu=\nu(\Cp)=\sqrt{\ev-V(\Cp)},
\end{equation*}
modulo terms whose contribution vanishes in the limit $r\to 0$,
and where the $v_j,$ $j=1,2,\ldots$
are orthonormal eigenfunctions of a harmonic oscillator. 
Here $\gamma_{j} = \gamma_{1,j}$ (for $Au_1$) or $\gamma_{2,j}$ (for
$Au_2$) is a Schwartz sequence in $j,$ hence interchanging
the order of various integrals and sums is permitted. Again,
the commutator $[\Delta,\chi_r]$
may be replaced by $2(x^2D_x\chi_r)(x^2D_x)$ and $x^2D_x$ must fall on
$e^{i\tilde\nu/x},$ to yield terms that do not vanish as $r\to 0$.
Now we interchange the integral in $Y$ and the summations, and use
the fact that the $v_k$ are orthonormal, to conclude that the limit is
\begin{equation*}
2\sqrt{\ev-V(\Cp)}\sum_{k=1}^\infty \gamma_{1,k}\,\overline{\gamma_{2,k}}.
\end{equation*}
Since this is also equal to
\begin{equation*}
2\sqrt{\ev-V(\Cp)}\int M_+(q,\ev)(u_1)\,\overline{M_+(q,\ev)(u_2)},
\quad M_+(q,\ev)(u_i)=\sum_{j=1}^\infty \gamma_{i,j} v_j,
\end{equation*}
the proof is complete.
\end{proof}

\begin{cor}
For $\ev \notin \Cv(V)$, the maps $M_+(\Rp,\ev)$ extend continously to
\begin{equation}
M_+(\Rp,\ev):E^0_{\ess}(\ev)\longrightarrow L^2(\bbR_{r(\Rp)}).
\label{HMV.156p}\end{equation}
The pairing $B$ is positive definite, hence it defines a norm on
$E^0_{\ess}(\ev)$ satisfying
\begin{equation}
0\leq \|u\|^2_0=2 \!\!\! \sum\limits_{\Rp\in\Min_+(\ev)} \!\!\!
\sqrt{\ev-V(\pi(\Rp))}\|M_+(\Rp,\ev)u\|^2=B(u,u)=i\langle[P-\ev,A]u,u\rangle
\label{eq:B-[P,A]}\end{equation}
for any $A\in\Psisc^{-\infty,0}(X)$ as in \eqref{HMV.147}. 
\label{cor:B-pos-def}\end{cor}

\begin{proof}
Equation \eqref{eq:B-[P,A]} holds for $u\in E^\infty_{\Min,+}(\ev),$
with the positivity coming from Proposition~\ref{HMV.153}. Since
$B$ extends to a pairing on $E^0_{\ess}(\ev),$ which is continuous
in the topology $\cT_3^0,$ \eqref{HMV.156p}
follows from the density statements of Proposition~\ref{HMV.164}
and Corollary~\ref{prop:density} (in $\cT_3^0$ in the latter case).
As all other expressions extend to continuous bilinear maps on
$E^0_{\ess}(\ev),$
the density statements of Proposition~\ref{HMV.164} and
Corollary~\ref{prop:density} show that
\eqref{eq:B-[P,A]} remains true. Since $B$ is non-degenerate and it is
a semi-norm, it follows that $B$ is in fact positive definite, hence a norm.
\end{proof}

\begin{lemma}
Suppose that $\ev\nin\Cv(V).$ For $r<-1/2$ there exists $C_r>0$ such that for
$u\in E^0_{\ess}(\ev),$
$\|u\|^2_{x^rL^2_{\scl}(X)}\leq C_r B(u,u).$
\label{lemma:x^rL^2<B}\end{lemma}

\begin{proof}
For simplicity of notation, suppose first that $\ev\nin\sigma_{\pp}(P).$

For $f\in x^s L^2_{\scl}(X),$ $s>1/2,$ let
\begin{equation*}
\tilde u=R(\ev+i0)f-R(\ev-i0)f.
\end{equation*}
For $u\in E^0_{\ess}(\ev)$ then
\begin{equation*}
i\langle u,f\rangle=\Bt(u,R(\ev+i0)f)=\Bt(Au,R(\ev+i0)f)=\Bt(Au,\tilde u)=
B(u,\tilde u).
\end{equation*}
By Cauchy-Schwarz and Corollary~\ref{cor:B-pos-def},
\begin{equation*}
B(u,\tilde u)\leq B(u,u)^{1/2} B(\tilde u,\tilde u)^{1/2}.
\end{equation*}
But
\begin{equation*}
B(\tilde u,\tilde u)=\Bt(R(\ev+i0)f,R(\ev+i0)f)
=i(\langle f,R(\ev+i0)f\rangle-\langle R(\ev+i0)f,f\rangle),
\end{equation*}
so by the limiting absorption principle,
\begin{equation*}
B(\tilde u,\tilde u)\leq \tilde C_s \|f\|_{x^s L^2}^2.
\end{equation*}
Combining these results gives that
\begin{equation}
|\langle u,f\rangle|\leq \tilde C_s^{1/2} B(u,u)^{1/2}\|f\|_{x^s L^2}.
\label{eq:pairing<norm}\end{equation}

Now for $r<-1/2,$ $u\in E^0_{\ess}(\ev),$
let $f=x^{-2r}u\in x^s L^2_{\scl}(X)$ for all $s<-2r-1/2,$
in particular for $s=-r>1/2.$ Applying \eqref{eq:pairing<norm} yields
\begin{equation*}
\|u\|_{x^r L^2}^2\leq
\tilde C_{-r}^{1/2} B(u,u)^{1/2}\|x^{-2r} u\|_{x^{-r} L^2}
=\tilde C_{-r}^{1/2} B(u,u)^{1/2}\|u\|_{x^{r} L^2}.
\end{equation*}
Cancelling a factor of $\|u\|_{x^{r} L^2}$ from both sides proves
the lemma if $\ev\nin\sigma_{\pp}(P).$

If $\ev\in\sigma_{\pp}(P),$
the calculations up to and including \eqref{eq:pairing<norm}
work with $f$ replaced by $\Pi f,$ $\Pi$ the orthogonal projection off
$E_{\pp}(\ev).$ Take $f=x^{-2r}u$ again and use the identity $\Pi u=u$ to
finish the proof in this case.
\end{proof}

An immediate consequence is the following proposition.

\begin{prop}
Suppose that $\ev\nin\Cv(V).$
The norm \eqref{eq:B-[P,A]} on $E^0_{\ess}(\ev)$ is equivalent
to the the norm
\begin{equation*}
\|u\|_{\cT_3^0}=\|u\|_3=\left(\|u\|^2_{x^rL^2_{\scl}(X)}
+\|B'u\|^2_{x^{-1/2}L^2_{\scl(X)}}\right)^{1/2},\ r<-1/2,
\end{equation*}
where $B'$ is as in \eqref{eq:B}. 
\label{prop:B-equiv-norm}\end{prop}

\begin{proof}
We may choose $A$ as above such that, in addition,
$i[P,A]\in\Psisc^{-\infty,1}(X)$ is
of the form $G^2+E,$ $G\in\Psisc^{-\infty,1/2}(X),$
$E\in \Psisc^{-\infty,2}(X).$ We can take, for example, $A$ with
$\sigma(A)=\chi(\nu)\phi^2(p),$ $p=\sigma(P),$ with $\phi\in\Cinf_c(\bbR)$
identically $1$ near $\ev,$ $\chi\in\Cinf(\bbR)$ identically $0$ on
$(-\infty,-a(\ev)/2),$ $1$ on $(a(\ev)/2,+\infty),$ and $\chi'\geq 0$ with
$(\chi')^{1/2}\in\Cinf(\bbR).$ Then the principal symbol of $i[P,A]$ is
$W\sigma(A),$ and $(W\sigma(A))^{1/2}$ is  real and $\Cinf$ by
\eqref{HMV.64}, so $i[P,A]=G^2+E$ as above. Now
\begin{equation*}
B(u,u)=|\langle [P-\ev,A] u,u\rangle|
\geq \|Gu\|^2-|\langle u,Eu\rangle|\geq \|Gu\|^2-C\|u\|^2
_{x^{-1}L^2},
\end{equation*}
so by Lemma~\ref{lemma:x^rL^2<B}, $\|Gu\|^2+\|u\|^2_{x^{-1}L^2}
\leq C'B(u,u).$

Conversely, 
\begin{equation*}
B(u,u)=|\langle [P-\ev,A] u,u\rangle|\leq C''(\|Gu\|_{x^{-1/2}L^2}+
\|u\|_{x^{-1}L^2})^2,
\end{equation*}
and $G$ may be chosen so as to satisfy the conditions of \eqref{eq:B}. 
\end{proof}

\begin{cor}
For $\ev\nin\Cv(V),$ $\Rp\in\Min_+(\ev),$ the Poisson operators
\begin{equation*}\begin{aligned}
P_+&(\Rp,\ev):\cS(\bbR)\to \SEEF(\ev) \\
P_+&(q,\ev)(g) = u \Longleftrightarrow u \in E^\infty_{\ess}(\ev), \
M_+(q,\ev)u = g, \ M_+(q',\ev)u = 0 \text{ for } q' \neq q
\end{aligned}\end{equation*}
extend to continuous linear maps
\begin{equation*}
P_+(\Rp,\ev):L^2(\bbR)\to E^0_{\ess}(\ev).
\end{equation*}
\end{cor}

\begin{proof}
By Proposition~\ref{prop:B-equiv-norm}, for $g\in\cS(\bbR),$
\begin{equation*}\begin{split}
&\|P(\Rp,\ev)g\|_3^2\leq C |B(P(\Rp,\ev)g,P(\Rp,\ev)g)|\\
&\qquad=2\sqrt{\ev-V(\pi(\Rp))}\int_{\bbR} |M_+(\Rp,\ev)(P(\Rp,\ev)g)|^2
\,\omega_{\Rp}
=2\sqrt{\ev-V(\pi(\Rp))}\int_{\bbR}|g|^2\,\omega_{\Rp},
\end{split}\end{equation*}
so the density of $\cS(\bbR)$ in $L^2(\bbR)$ finishes the proof.
\end{proof}

The combination of these results now yields the following theorem.

\begin{thm}[Theorem~\ref{thm:AC-ev}]\label{thm:AC2}
For $\ev \notin \Cv(V)$, the joint map
\begin{equation}
M_+(\ev)=E^0_{\ess}(\ev)\to\oplus_{\Rp\in\Min_+(\ev)}L^2(\bbR),
\ M_+(\ev)u= (M_+(\Rp,\ev)(u))_{\Rp\in\Min_+(\ev)}
\end{equation}
is an isomorphism with inverse
\begin{equation*}
P_+(\ev):\oplus_{\Rp\in\Min_+(\ev)}L^2(\bbR)\to E^0_{\ess}(\ev),\
P_+(\ev)(g_q)_{q \in \Min_+(\ev)} = \sum\limits_q P_+(q,\ev)g_q
\end{equation*}
\end{thm}

We can also express the pairing \eqref{HMV.154} in terms of the expansions
of microlocal eigenfunctions $u \in E_{\mic,-}(q,\ev)$, $q \in
\Min_-(\ev)$, \ie in the incoming region. That leads to maps $M_-(q,\ev)$
and $P_-(q,\ev)$, for $q \in \Min_-(\ev)$ and Theorem~\ref{thm:AC2} holds
with all plus signs changed to minus signs. It is also convenient to
change the notation slightly and identify the spaces $L^2(\bbR)$ for
the incoming and outgoing radial points over $z\in\Min$.

\begin{cor} For $\ev \notin \Cv(V)$, 
the S-matrix may be identified as the unitary operator
$S(\ev)=M_+(\ev)P_-(\ev)$ on $\oplus_{z\in\Min}L^2(\bbR).$
\end{cor}

This theorem is essentially a pointwise version of asymptotic completeness
in $\ev.$ Integrating in $\ev$ gives a version of the usual statement.
Namely, let $I\subset\bbR\setminus\Cv(V)$ be a compact interval.
For $f\in\dCI(X)$ orthogonal to $E_{\pp}(I),$
$u=u(\ev)=R(\ev+i0)f,$ as discussed above,
\begin{equation}\begin{split}
B(u,u)=B_{\ev}(u(\ev),u(\ev))&=i\left(\langle f,R(\ev+i0)f\rangle
-\langle R(\ev+i0)f,f\rangle\right)\\
&=i\langle f,[R(\ev+i0)-R(\ev-i0)]f\rangle
=2\pi \langle f,\Sp(\ev)f\rangle.
\end{split}\end{equation}
Integrating over $\ev$ in $I,$ denoting
the spectral projection of $P$ to $I$ by $\Pi_I,$ and writing $M_+(q,\ev)=0$
if $V(z)>\ev,$ $z = \pi(q)$, we deduce that
\begin{equation}
\|\Pi_I f\|^2=\sum_{\Rp\in\Min_+(\ev)}4\pi\sqrt{\ev-V(\Cp)}
\int\|M_+(q,\ev)R(\ev+i0)f\|^2_{L^2(\bbR),
\omega_{\Rp}}\,d\ev,
\label{eq:Pi_I-isometry}\end{equation}
so $M_+\circ R(.+i0)$ is an isometry on the orthocomplement of the
finite dimensional space $E_{\pp}(I)$ in the range of $\Pi_I$.

\begin{thm}[Asymptotic completeness]\label{thm:AC}
If $I\subset\bbR\setminus\Cv(V)$ is compact then
\begin{equation*}
M_+\circ R(.+i0):{\operatorname{Ran}}( \Pi_I)\ominus E_{\pp}(I)\to
L^2(I\times \bbR)
\end{equation*}
is unitary.
\end{thm}

\def\cprime{$'$} \def\cprime{$'$} \def\cprime{$'$}
  \def\polhk#1{\setbox0=\hbox{#1}{\ooalign{\hidewidth
  \lower1.5ex\hbox{`}\hidewidth\crcr\unhbox0}}} \def\cprime{$'$}
  \def\cprime{$'$} \def\cprime{$'$} \def\cprime{$'$} \def\cprime{$'$}
  \def\cprime{$'$} \def\cprime{$'$} \def\cprime{$'$} \def\cprime{$'$}
  \def\cprime{$'$} \def\cprime{$'$} \def\cprime{$'$} \def\cprime{$'$}
  \def\cprime{$'$} \def\cprime{$'$} \def\cprime{$'$} \def\cprime{$'$}
  \def\cprime{$'$} \def\cprime{$'$} \def\cprime{$'$} \def\cprime{$'$}
  \def\cprime{$'$} \def\cprime{$'$} \def\cprime{$'$} \def\cprime{$'$}
  \def\cprime{$'$} \def\cprime{$'$}
\providecommand{\bysame}{\leavevmode\hbox to3em{\hrulefill}\thinspace}
\providecommand{\MR}{\relax\ifhmode\unskip\space\fi MR }
\providecommand{\MRhref}[2]{%
  \href{http://www.ams.org/mathscinet-getitem?mr=#1}{#2}
}
\providecommand{\href}[2]{#2}


\end{document}